\documentclass[12pt]{article}
\usepackage{amssymb,amsmath,amsopn,amsfonts}
\usepackage[dvips]{color}
\usepackage{colordvi,multicol}
\usepackage{epsf}
\newfam\msbfam
\font\tenmsb=msbm10 \textfont\msbfam=\tenmsb \font\sevenmsb=msbm7
\scriptfont\msbfam=\sevenmsb \font\fivemsb=msbm5
\scriptscriptfont\msbfam=\fivemsb

\def\th#1{\vspace{1mm}\noindent{\bf #1}\quad}

\def\proof{\vspace{1mm}\noindent{\it Proof}\quad}

\numberwithin{equation}{section}

\def\bc{\begin{center}}
\def\ec{\end{center}}
\def\no{\noindent}
\def\hang{\hangindent\parindent}
\def\textindent#1{\indent\llap{\qquad #1\ \ \enspace}\ignorespaces}
\def\ref{\par\hang\textindent}

\begin{document}
\centerline{\Large\bf BSDE and generalized Dirichlet forms:}
\centerline{\Large\bf the finite
dimensional case}

\vspace{0.5 true cm}

 \centerline{RONG-CHAN ZHU$^{\textrm{a,b}}$} \centerline{\small
a. Institute of Applied Mathematics, Academy of Mathematics and Systems
Science,} \centerline{\small Chinese Academy of Sciences, Beijing
100190, China}
\centerline{\small b. Department of Mathematics, University of Bielefeld, D-33615 Bielefeld, Germany}
\centerline{\small E-mail: zhurongchan@126.com}
\footnotetext{\footnotesize Research supported by 973 project, NSFC, key Lab of CAS, the DFG through IRTG 1132 and CRC 701 and the I.Newton Institute, Cambridge, UK}
 \vskip
1.4cm
\begin{abstract}
\vskip 0.1cm \noindent We consider the following
quasi-linear parabolic system of backward partial differential equations
 $$(\partial_t+L)u+f(\cdot,\cdot,u, \nabla u\sigma)=0 \textrm{ on } [0,T]\times \mathbb{R}^d\qquad u_T=\phi,$$
 where $L$ is a possibly degenerate  second order differential operator with merely measurable coefficients. We solve this system in the framework of generalized Dirichlet forms and employ the stochastic calculus
 associated to the Markov process with generator $L$ to obtain a probabilistic representation of the solution $u$ by solving the corresponding backward stochastic differential equation.
The solution satisfies the corresponding mild equation which is equivalent to being a generalized solution of the PDE. A further main result is the generalization of  the martingale representation theorem
using  the stochastic calculus associated to the generalized Dirichlet form given by $L$. The nonlinear term $f$ satisfies a monotonicity condition with respect to $u$ and a Lipschitz condition with respect to $\nabla u$.
\end{abstract}
\no{\footnotesize{\bf Keywords}:\hspace{2mm}  backward stochastic differential equations, quasi-linear parabolic partial differential equations, Dirichlet
forms, generalized Dirichlet forms, Markov processes, martingale representation}

\section{Introduction}

 Consider the following quasi-linear parabolic system of backward partial differential equations
 $$(\partial_t+L)u+f(\cdot,\cdot,u, \nabla u\sigma)=0 \textrm{ on } [0,T]\times \mathbb{R}^d\qquad u_T=\phi,\eqno(1.1)$$
 where $L$ is a second order linear differential operator and $f$ is monotone in $u$ and Lipschitz in $\nabla u$. If $L$ has sufficiently regular coefficients there is a well-known theory,  to obtain a probabilistic representation of the solutions to
(1.1), using corresponding backward stochastic differential equations (BSDE) and also to solve BSDE with the help of (1.1), originally due to E. Pardoux and S. Peng (\cite{PP}).
 The main aim of this paper is to implement this approach for a very general class of linear operators
 $L$, which are possibly degenerate, have merely measurable cofficients and are in general not symmetric. Solving (1.1) for such general $L$ is the first main task of this paper.
The second main contribution is to prove the martingale representation theorem
for the underlying reference diffusion generated by such general operators $L$.

If $f$ and the coefficients of the second-order differential operator $L$ are sufficiently smooth, the PDE  has a classical solution $u$. Consider $Y_t^{s,x}:=u(t,X_t^{s,x})$, $Z_t^{s,x}:=\nabla u\sigma(t,X_t^{s,x})$ where $X_t^{s,x}, s\leq t\leq T$, is the diffusion process with infinitesimal generator $L$ which starts from $x$ at time $s$. Then, using It\^{o}'s formula one checks that $(Y_t^{s,x}, Z_t^{s,x})_{s\leq t\leq T}$ solves the BSDEs
$$Y_t^{s,x}=\phi(X_T^{s,x})+\int_t^Tf(r,X_r^{s,x},Y_r^{s,x},Z_r^{s,x})dr-\int_t^TZ_r^{s,x}dB_r.\eqno(1.2)$$
Conversely, by standard methods one can prove that (1.2) has a unique solution $(Y_t^{s,x}, Z_t^{s,x})_{s\leq t\leq T}$ and then $u(s,x):=Y_s^{s,x}$ is a solution to PDE (1.1).
 BSDEs have been introduced
(in their actual form) by Pardoux and Peng \cite{PP} and found applications in
stochastic control and mathematical finance. If $f$ and the coefficients of $L$ are  Lipschitz continuous then a series of papers (e.g. \cite{BPS05}, \cite{P} and the reference therein) prove that the above relation between PDE (1.1) and BSDE (1.2) remains true, if one considers viscosity solutions to PDE (1.1). In both these approaches, since the coefficients are Lipschitz continuous, the Markov process $X$ with infinitesimal operator $L$ is a diffusion process which satisfies an SDE and so one may use its associated stochastic calculus.

In \cite{BPS05} Bally, Pardoux and Stoica consider a semi-elliptic symmetric second-order differential operator $L$ ( which is written in divergence form ) with measurable coefficients. They prove that the above system of PDE has a unique solution $u$ in some functional space. Then using the theory of \emph{symmetric} Dirichlet forms and its associated stochastic calculus, they prove that the solution $Y^{s.x}$ of the BSDE yields a precised version of the  solution $u$ so that, moreover, one has $Y_t^{s,x}=u(t,X_{t-s}), P^x$-a.s. In \cite{S}, the analytic part of \cite{BPS05} has been generalized to a non-symmetric case with $L$ satisfying the weak sector condition. Here the weak sector condition means
$$((1-L)u,v)\leq K((1-L)u,u)^{1/2}((1-L)v,v)^{1/2}, \textrm{ for } u,v\in \mathcal{D}(L),$$ for some constant $K>0$.  A.Lejay (\cite{L}) and L. Stoica (\cite{Sto}) consider the generator $L=\frac{1}{2}\sum_{i,j=1}^d\frac{\partial}{\partial x_i}(a_{ij}\frac{\partial}{\partial x_j})+\sum_{i=1}^db_i(x)\frac{\partial}{\partial x_i}$ for bounded $a,b$. In their case, the operator $L$  generates a semigroup $(P_t)$, which possesses continuous densities satisfying Aronson's estimate. In \cite{ZR},  T.S. Zhang and Q.K.Ran (see also \cite{Z}) consider  $L$ of a more general form, but $a=(a_{ij})$ is required to be uniformly elliptic  and $b\in L^p$ for $p>d$. Anyway, since $L$ satisfies the weak sector condition in this case, it generates a sectorial ( i.e. a small perturbation of a symmetric)  Dirichlet form, so the theory of Dirichlet forms from \cite{MR} can be applied in \cite{Z}, \cite{ZR}.

In \cite{St} Stannat extends the known framework of Dirichlet forms to the class of generalized Dirichlet forms. By this we can analyze differential operators where the second order part may be degenerate and at the same time the first order part may be unbounded satisfying no global $L^p$-condition for $p\geq d$. The motivation for this paper is to extend the results in \cite{BPS05} to the case, where $L$ generates a generalized Dirichlet form so that we can allow the coefficients of $L$ to be more general.

In this paper, we consider PDE (1.1) for a non-symmetric second order differential operator $L$, which is associated to the bilinear form
$$\aligned\mathcal{E}(u,v):=&\sum_{i,j=1}^d\int a_{ij}(x)\frac{\partial u}{\partial x^i}(x)\frac{\partial v}{\partial x^j}(x)m(dx)+\int c(x)u(x)v(x)\\&+\sum_{i=1}^d\int\sum_{j=1}^d a_{ij}(x)(b_j(x)+\hat{b}_j(x)) \frac{\partial u}{\partial x^i} v(x)m(dx)\textrm{ }\forall u,v \in C_0^\infty(\mathbb{R}^d).\endaligned\eqno(1.3)$$
for $u,v\in C_0^\infty(\mathbb{R}^d)$, where $C_0^\infty(\mathbb{R}^d)$ denotes the space of infinitely differentiable functions with compact support. We stress that $(a_{ij})$ is not necessarily assumed to be (locally) strictly positive definite, but may be degenerate in general. When $b\equiv0$, the bilinear form $\mathcal{E}$ satisfies the weak sector condition. In the analytic part, the condition we assume for the perturbation term given by $b$ is $b\sigma \in L^2(\mathbb{R}^d;\mathbb{R}^d,m)$. Here $\sigma\sigma^*=a$ and $\sigma^*$ is the transpose of the matrix of $\sigma$. In the probabilistic part, we have more general condition on $b$ (i.e. $b\in L^p_{\rm{loc}}$, for $p>d$ see Example 5.6).  That implies that we do not have the weak sector condition for the bilinear form. We use the theory of generalized Dirichlet forms and its associated stochastic calculus ( cf [22, 23] [25, 26]) to generalize the results in \cite{BPS05}. Here $m$ is a finite measure or Lebesgue measure on $\mathbb{R}^d$. If $D$ is a bounded open domain, we choose $m$ as $1_D(x)dx$. Then in certain cases the solution of  PDE (1.1) satisfies  the Neumann boundary condition. If we replace $C_0^\infty(\mathbb{R}^d)$ by $C^\infty_0(D)$, the solution of PDE (1.1) satisfies the Dirichlet boundary condition.

In the analytic part of our paper, we do not need $\mathcal{E}$ to be a generalized Dirichlet form. We start from a semigroup $(P_t)$ satisfying conditions (A1)-(A4), specified in Section 2 below. Such a semigroup can, however, be constructed from a generalized Dirichlet form. It can also be constructed by other methods (see e.g. \cite{DR}). Under conditions (A1)-(A4), the coefficients of $L$ may be quite singular and only very broad assumptions on $a$ and $b$ are needed (see the examples in Section 4 and Section 5, which to the best of our knowledge could not be covered by
other results).

The paper is organized as follows. In Sections 2 and 3, we use functional analytical methods to solve PDE (1.1) in the sense of Definition 2.3, i.e.
there are sequences  $\{u^n\}$ which are strong solutions with data $(\phi^n,f^n) $ such that $$\|u^n-u\|_T\rightarrow0,\|\phi^n-\phi\|_2\rightarrow0,\lim_{n\rightarrow\infty} f^n=f \textrm{ in } L^1([0,T];L^2).$$
Here $\|\cdot \|_T:=(\sup_{t\leq T}\|\cdot\|_2^2+\int_0^T\mathcal{E}^{a,\hat{b}}_{c_2+1}(\cdot)dt)^{1/2},$ where $\mathcal{E}^{a,\hat{b}}$ is the  summand in the left hand side of (1.3) with $b\equiv0$.
The above definition for the solution is equivalent to that of the following mild equation in $L^2$-sense
$$u(t,x)=P_{T-t}\phi(x)+\int_t^TP_{s-t}f(s,\cdot,u_s,D_\sigma u_s)(x)ds.$$
 If we use the definition of  weak solution to define our solution as in [2],  uniqueness of the solution can not be obtained since only $|b\sigma|\in L^2(\mathbb{R}^d;m)$.
Furthermore, the function $f$ in PDE (1.1) need not be Lipschitz continuous with respect to the third variable; monotonicity suffices. And $\mu$ which appears in the monotonicity  conditions (see condition (H2) in Section 3.2 below) can depend on $t$. $f$ is, however, assumed to be Lipschitz continuous with respect to the last variable. We emphasize that the first order term of $L$ cannot be incorporated into $f$ without the condition that $b$ is bounded. Hence we are forced to take it as part of $L$ and hence have to consider a diffusion process $X$ in (1.2) which is generated by an operator $L$ which is the generator of a (in general non-sectorial) generalized Dirichlet form. We also emphasize that under our conditions, PDE (1.1) cannot be tackled by standard monotonicity methods (see e.g. \cite{Ba10}) because of the lack of a suitable Gelfand triple $V\subset H\subset V^*$ with $V$ being a reflexive Banach space.

In Section 4, we extend the stochastic calculus of generalized Dirichlet forms in order to generalize the martingale representation theorem. In order to treat BSDE, we show in Theorem 4.8 that there exists a set of null capacity $\mathcal{N}$ outside of which the following representation theorem holds : for  every bounded $\mathcal{F}_\infty$-measurable random variable $\xi$, there  exists a predictable process $(\phi_1,...,\phi_d):[0,\infty)\times \Omega\rightarrow \mathbb{R}^d$, such that for each probability measure $\nu$, supported by $\mathbb{R}^d\setminus\mathcal{N}$, one has
$$\xi=E^\nu(\xi|\mathcal{F}_0)+\sum_{i=0}^d\int_0^\infty\phi_s^idM_s^{(i)}\qquad P^\nu-a.s..$$
As a result, one can choose the exceptional set $\mathcal{N}$ such that if the process $X$ starts from a point of $\mathcal{N}^c$, it remains always in this set. As a consequence we deduce the existence of solutions for the BSDE using the existence for PDE (1.1) in the usual way, however, only under $P^m$, because of our general coefficients of $L$ (c.f. Theorem 4.12).

In Section 5, we employ the martingale representation to deduce existence and uniqueness for the solutions of BSDE (1.2). Moreover, we give Example 5.6 satisfying Assumption (A5), which to the best of our knowledge could not be treated by
other techniques. As a consequence, in Theorem 5.7, existence of solutions for PDE (1.1), not covered by our analytic results in Section 2, is obtained by  $u(s,x)=Y_s^s$, where $Y_t^s$ is the solution of the BSDE. Moreover we have, $Y^{s}_t=u(t,X_{t-s}),P^x$-a.s., $x\in R^d\backslash \mathcal{N}$.

\section{Preliminaries}

Let $\sigma:\mathbb{R}^d\mapsto \mathbb{R}^d\otimes \mathbb{R}^k$ be a measurable map. Then  there exists a measurable map $\tau:\mathbb{R}^d\mapsto \mathbb{R}^k\otimes \mathbb{R}^d$ such that
$$\sigma\tau=\tau^*\sigma^*,\qquad \tau\sigma=\sigma^*\tau^*,\qquad \sigma\tau\sigma=\sigma.$$ (see e.g. [2, Lemma A.1]). Here $\sigma^*$ is the transpose of the matrix of $\sigma$. Then $a:=\sigma \sigma^*=(a_{ij})_{1\leq i,j\leq d}$ takes values in the space of symmetric non-negative
definite matrices. Let also $b:\mathbb{R}^d\rightarrow \mathbb{R}^d$ be measurable.
Assume that the basic measure $m(dx)$ for the generalized Dirichlet form, to be defined below, is a finite measure or Lebesgue measure on $\mathbb{R}^d$.

Denote  the Euclidean norm and  the scalar product in $\mathbb{R}^d$ by $|\cdot|$, $\langle\cdot,\cdot\rangle$ respectively,  while on the space of matrices $\mathbb{R}^d\otimes \mathbb{R}^k$ we use the trace scalar product and its associated norm, i.e., for $z=(z_{ij})\in \mathbb{R}^d\otimes \mathbb{R}^k$,
 $\langle z_1,z_2\rangle=\rm{trace}(z_1z_2^*), |z|=(\sum_{i=1}^d\sum_{j=1}^kz_{ij}^2)^{1/2}.$ Let $L^2$, $L^2(\mathbb{R}^d;\mathbb{R}^k)$ denote $L^2(\mathbb{R}^d,m)$, $L^2(\mathbb{R}^d,m;\mathbb{R}^k)$ respectively. And $(\cdot,\cdot)$ denotes the $L^2$-inner product. And for $1\leq p\leq\infty$, $\|\cdot\|_p$ denotes the usual norm in $L^p(\mathbb{R}^d;m)$. If $W$ is a function space, we will use $bW$ to denote the bounded function in $W$.

Furthermore, let $a_{ij},\sum_{j=1}^da_{ij}b_j,\sum_{j=1}^da_{ij}\hat{b}_j\in L_{\rm{loc}}^1(\mathbb{R}^d,m)$ and $c\in L_{\rm{loc}}^1(\mathbb{R}^d,\mathbb{R}^+;m) $. We introduce the bilinear form
$$\aligned\mathcal{E}(u,v):=&\sum_{i,j=1}^d\int a_{ij}(x)\frac{\partial u}{\partial x^i}(x)\frac{\partial v}{\partial x^j}(x)m(dx)+\int c(x)u(x)v(x)\\&+\sum_{i=1}^d\int\sum_{j=1}^d a_{ij}(x)(b_j(x)+\hat{b}_j(x)) \frac{\partial u}{\partial x^i} v(x)m(dx)\textrm{ }\forall u,v \in C_0^\infty(\mathbb{R}^d).\endaligned$$
Consider the following conditions:
\vskip.10in
  \noindent(A1) The bilinear form
 $$\mathcal{E}^a(u,v)=\sum_{i,j=1}^d\int a_{ij}(x)\frac{\partial u}{\partial x^i}(x)\frac{\partial v}{\partial x^j}(x)m(dx)\textrm{ }\forall u,v\in C_0^\infty(\mathbb{R}^d),$$ is closable on $L^2(\mathbb{R}^d,m)$.
  \vskip.10in

 Define $\mathcal{E}^a_1(\cdot,\cdot):=\mathcal{E}^a(\cdot,\cdot)+(\cdot,\cdot)$. The closure of $C_0^\infty(\mathbb{R}^d)$ with respect to $\mathcal{E}^a_1$ is denoted by $F^a$. Then $(\mathcal{E}^a,F^a)$
is a well-defined symmetric Dirichlet form on $L^2(\mathbb{R}^d,m)$.

 \vskip.10in

  For the bilinear form $$\aligned\mathcal{E}^{a,\hat{b}}(u,v):&=\sum_{i,j=1}^d\int a_{ij}(x)\frac{\partial u}{\partial x^i}(x)\frac{\partial v}{\partial x^j}(x)m(dx)+\int c(x)u(x)v(x)\\&+\sum_{i=1}^d\int\sum_{j=1}^d a_{ij}(x)\hat{b}_j(x) \frac{\partial u}{\partial x^i} v(x)m(dx),\endaligned$$ we consider the following conditions:

 \vskip.10in

  \noindent(A2) There exists a constant $c_2\geq0$ such that $\mathcal{E}^{a,\hat{b}}_{c_2}(\cdot,\cdot):=\mathcal{E}^{a,\hat{b}}(\cdot,\cdot)+c_2(\cdot,\cdot)$ is a coercive closed form (see e.g. \cite{MR}), and there exist constants $c_1,c_3>0$ such that for $u\in C_0^\infty(\mathbb{R}^d)$
 \begin{equation}c_1\mathcal{E}^a(u,u)\leq \mathcal{E}_{c_2}^{a,\hat{b}}(u,u),\end{equation}
and \begin{equation}\int cu^2 dm \leq c_3 (\mathcal{E}^a(u,u)+\|u\|_2^2).\end{equation}
 \vskip.10in

Denote $\tilde{\mathcal{E}}^{a,\hat{b}}(u,v):=\mathcal{E}^{a,\hat{b}}(u,v)+\mathcal{E}^{a,\hat{b}}(v,u)$. The closure of $C_0^\infty(\mathbb{R}^d)$ with respect to $\tilde{\mathcal{E}}_{c_2+1}^{a,\hat{b}}$ is denoted  by $F$. By (2.1) we have $F\subset F^a$. And  for $u\in F$ (2.1) and (2.2) are satisfied.
  \vskip.10in

\noindent(A3) $| b\sigma|\in L^2(\mathbb{R}^d;m)$ and there exists $\alpha\geq 0$ such that
\begin{equation}\int\langle b\sigma,(\nabla u^2)\sigma\rangle dm \geq -\alpha \|u\|_2^2 \textrm{ } u\in C_0^\infty(\mathbb{R}^d). \end{equation}
\vskip.10in

\noindent(A4) There exists a positivity preserving $C_0$-semigroup $P_t$ on $L^1(\mathbb{R}^d;m)$ such that for any $t\in[0,T], \exists C_T>0$ such that $$\|P_t f\|_\infty \leq C_T \|f\|_\infty.$$ Then for $0\leq t\leq T$,  $P_t$ extends to a semigroup on $L^p(\mathbb{R}^d;m)$ for all $p\in[1,\infty)$ by the Riesz-Thorin Interpolation Theorem (denoted by $P_t$ for simplicity) which is strongly continuous on $L^p(\mathbb{R}^d;m)$. We denote its $L^2$-generator by  $(L,\mathcal{D}(L))$ and assume that $b\mathcal{D}(L)\subset bF$, and for any $u\in bF$ there exists uniformly bounded $u_n\in \mathcal{D}(L)$ such that  $\tilde{\mathcal{E}}_{c_2+1}^{a,\hat{b}}(u_n-u)\rightarrow0$ and that it is associated with the bilinear form in the sense that $\mathcal{E}(u,v)=-(Lu,v)$ for $u,v\in b\mathcal{D}(L)$.
 \vskip.10in

We emphasize that in contrast to previous work $P_t$ in (A4) is no longer analytic on $L^2(\mathbb{R}^d;m)$.
By (A4) there exist constants $M_0,c_0$ such that
\begin{equation}\|P_tf\|_2\leq M_0e^{c_0 t}\|f\|_2, \textrm{   } \forall f\in L^2(\mathbb{R}^d;m).\end{equation}

\vskip.10in

To obtain a semigroup $P_t$ satisfying the above conditions, we can use generalized Dirichlet form.
Let us recall the definition of a generalized Dirichlet form from \cite{St}.
 Let $E$ be a Hausdorff topological space and assume that its Borel
$\sigma$-algebra $\mathcal{B}(E)$ is generated by the set $C(E)$ of
all continuous functions on $E$. Let $m$ be a $\sigma$-finite
measure on $(E,\mathcal{B}(E))$ such that $\mathcal{H}:=L^2(E,m)$ is
a separable (real) Hilbert space. Let $(\mathcal{A},\mathcal{V})$ be
a coercive closed form on $\mathcal{H}$ in the sense of \cite{MR}. We will always denote the
corresponding norm by $\|\cdot\|_{\mathcal{V}}$. Identifying
$\mathcal{H}$ with its dual $\mathcal{H}'$ we obtain that $
\mathcal{V}\rightarrow\mathcal{H}\cong\mathcal{H}'\rightarrow\mathcal{V}'$ densely and continuously.

Let $(\Lambda,D(\Lambda,\mathcal{H}))$ be a linear operator on
$\mathcal{H}$ satisfying the following assumptions:

(i) $(\Lambda,D(\Lambda,\mathcal{H}))$ generates a $C_0$-semigroup
of contractions $(U_t)_{t\geq0}$ on $\mathcal{H}$.

(ii)$\mathcal{V}$ is $\Lambda$-admissible, i.e. if $(U_t)_{t\geq}$
can be restricted to a $C_0$-semigroup on $\mathcal{V}$.

Let $(\Lambda,\mathcal{F})$ with corresponding norm $\|\cdot\|_{\mathcal{F}}$ be the closure of
$\Lambda:D(\Lambda,\mathcal{H})\cap\mathcal{V}\rightarrow\mathcal{V}'$ as an operator from $\mathcal{V}$ to $\mathcal{V}'$
 and $(\hat{\Lambda},\hat{\mathcal{F}})$ its dual operator.

 Let\begin{equation*}
\mathcal{E}(u,v)=\left\{\begin{array}{ll} \mathcal{A}(u,v)-\langle\Lambda u,v\rangle&\ \ \ \textrm{if} \ u\in \mathcal{F},v\in \mathcal{V}\\
\mathcal{A}(u,v)-\langle\hat{\Lambda} v,u\rangle&\ \ \ \ \textrm{if}
\ u\in \mathcal{V},v\in \hat{\mathcal{F}},
\end{array}\right.
\end{equation*}
Here  $\langle\cdot,\cdot\rangle$ denotes the dualization between
$\mathcal{V}'$ and $\mathcal{V}$ and $\langle\cdot,\cdot\rangle$ coincides with the inner product $(\cdot,\cdot)_H$ in $H$ when restricted
to $H\times \mathcal{V}$. And
$\mathcal{E}_{\alpha}(u,v):=\mathcal{E}(u,v)+\alpha(u,v)_{\mathcal{H}}$
for $\alpha>0$. We call $\mathcal{E}$ the bilinear form associated
with $(\mathcal{A},\mathcal{V})$ and
$(\Lambda,D(\Lambda,\mathcal{H}))$. If $$u\in \mathcal{F}\Rightarrow
u^+\wedge1\in \mathcal{V} \textrm{ and }
\mathcal{E}(u,u-u^+\wedge1)\geq0$$ then the bilinear form is called
a generalized Dirichlet form. If the adjoint semigroup $(\hat{U}_t)_{t\geq0}$ of $(U_t)_{t\geq 0}$ can also be restricted to a $C_0$-semigroup on $\mathcal{V}$. Let $(\hat{\Lambda},D(\hat{\Lambda},\mathcal{H}))$ denote the generator of $(\hat{U}_t)_{t\geq0}$ on $\mathcal{H}$, $\hat{\mathcal{A}}(u,v):=\mathcal{A}(v,u), u,v\in \mathcal{V}$ and let the coform $\hat{\mathcal{E}}$ be defined as the bilinear form associated with $(\hat{\mathcal{A}},\mathcal{V})$ and
$(\hat{\Lambda},D(\hat{\Lambda},\mathcal{H}))$.

 \vskip.10in

\th{Remark 2.1} (i) Some general criteria imposing conditions on $a$ in order that $\mathcal{E}^a$ be closable
are e.g. given in  [12, Section 3.1] and [14, Chap II, Section 2].

(ii) There are examples considered in [14, Chap. II, Subsection 2d] satisfying (A2). Assume the Sobolev inequality
$$\|u\|_q\leq C(\mathcal{E}^a(u,u)+\|u\|_2^2)^{1/2}, \textrm{  } \forall u\in C_0^\infty(\mathbb{R}^d),$$ is satisfied,
where $\frac{1}{q}+\frac{1}{d}=\frac{1}{2}$ and $\|\cdot\|_q$ denotes the usual norm in $L^q$. If $|\hat{b}\sigma|\in L^d(\mathbb{R}^d;m)+L^\infty(\mathbb{R}^d;m)$ and $c\in L^{d/2}(\mathbb{R}^d;m)+L^\infty(\mathbb{R}^d;m)$, then (A2) is satisfied. In \cite{ZR} they consider the bilinear form $Q(u,v)=\mathcal{E}^{a,\hat{b}}(u,v)+\int\langle d_1(x),\nabla v(x)\rangle u(x)dm$, here $d_1\in L^q(\mathbb{R}^d), q>d$. In their case, the result for the existence of the nonlinear PDE can be obtained by [18, Theorem 4.2.4]  since  the nonlinear part is Lipschitz in $u$ and $\nabla u$. In our case, we have more general conditions on $b$ and $f$, so that  we can not find  a suitable Gelfand triple $V\subset H\subset V^*$ with $V$ being a reflexive Banach space and use monotonicity methods as in \cite{PR}.

(iii) We can construct a semigroup $P_t$ satisfying (A4) by the theory of generalized Dirichlet form. More precisely, suppose  there exists a constant $\hat{c}\geq0$ such that $\mathcal{E}_{\hat{c}}(\cdot,\cdot):=\mathcal{E}(\cdot,\cdot)+\hat{c}(\cdot,\cdot)$ is a generalized Dirichlet form with domain $\mathcal{F}\times \mathcal{V}$ in one of the following three senses:

(a)$(E,\mathcal{B}(E),m)=(\mathbb{R}^d,\mathcal{B}(\mathbb{R}^d),m)$,

 $(\mathcal{E}^{a,\hat{b}}_{c_2},F)=(\mathcal{A},\mathcal{V})$,

  $-\langle\Lambda u,v\rangle-(\hat{c}-c_2)(u,v)=\sum_{i}^d\int \sum_{j=1}^da_{ij}(x)b_j(x) \frac{\partial u}{\partial x^i} v(x)m(dx)$  for $u,v\in C_0^\infty(R^d)$;

(b)$(E,\mathcal{B}(E),m)=(\mathbb{R}^d,\mathcal{B}(\mathbb{R}^d),m)$,

$\mathcal{A}\equiv0$ and $\mathcal{V}=L^2(\mathbb{R}^d,m)$,

$-\langle\Lambda u,v\rangle=\mathcal{E}_{\hat{c}}(u,v)$ for $u,v\in C_0^\infty(\mathbb{R}^d)$ and $C_0^\infty(\mathbb{R}^d)\subset \mathcal{D}(L)$;

(c) $\mathcal{E}_{\hat{c}}=\mathcal{A}$, $\Lambda\equiv0$ (In this case $(\mathcal{E}_{\hat{c}},\mathcal{V})$ is a sectorial Dirichlet form in the sense of \cite{MR}).

Then there exists a sub-Markovian $C_0$-semigroup of contractions $P_t^{\hat{c}}$ associated with the generalized Dirichlet form $\mathcal{E}_{\hat{c}}$. Define $P_t:=e^{{\hat{c}}t}P_t^{\hat{c}}$.  If it is a $C_0$-semigroup on $L^1$ then it satisfies (A4).  Then we have $$\mathcal{D}(L)\subset \mathcal{F}\subset F.$$

(iv) The semigroup can be also constructed by other methods. (see e.g. \cite{DR}, \cite{BBR}, \cite{BDR}).

(v)  By (A3) we have that $\mathcal{E}$ is positivity preserving i.e.
$$\mathcal{E}(u,u^+)\geq 0\textrm{  }\forall u\in \mathcal{D}(L), $$
which can be obtained by the same arguments as [23, Proposition 4.4].

(vi) The condition that for any $u\in bF$ there exists uniformly bounded $u_n\in \mathcal{D}(L)$ such that  $\tilde{\mathcal{E}}_{c_2+1}^{a,\hat{b}}(u_n-u)\rightarrow0$ is satisfied if $C_0^\infty(\mathbb{R}^d)\subset \mathcal{D}(L)$. It can also be satisfied in the case of  (iii) by the theory of generalized Dirichlet form.

(vii) All the conditions are satisfied by  the bilinear form considered in \cite{DR}, \cite{L}, [22, Section 1 (a)]  and the following example which is considered in \cite{St}.

(viii) The notion of quasi-regularity for generalized Dirichlet forms analogously to \cite{MR} has been introduced in \cite{St}. By this and a technical assumption an associated $m$-tight special standard process can be constructed. We will use stochastic calculus associated with this process to conclude our probabilistic results (see Section 4 below).
\vskip.10in

\th{Example 2.2} Let $b_i\in L^2(\mathbb{R}^d;dx)$, $1\leq i\leq d$. Consider the bilinear form $$\mathcal{E}(u,v):=\sum_{i,j=1}^d\int_{\mathbb{R}^d} \frac{\partial u}{\partial x_i}\frac{\partial v}{\partial x_j}dx-\sum_{i=1}^d\int b^i\frac{\partial u}{\partial x_i}vdx; u,v\in C_0^\infty(\mathbb{R}^d)$$
Assume there exist constants $c,L\geq0$ such that
$$\int\langle b, \nabla u\rangle dx\leq 2c\|u\|_1 \textrm{ for all } u\in C_0^\infty(\mathbb{R}^d),u\geq0,$$
$$-\sum_{i,j=1}^d\int b_i\frac{\partial u}{\partial x_j}dxh_ih_j\leq L\|u\|_1|h|^2,$$
 $$\textrm{ for all } u\in C_0^\infty(\mathbb{R}^d),u\geq0, h\in \mathbb{R}^d,$$
 (or equivalently, $b$ is quasi-monotone, i.e. $$\langle b(x)-b(y),x-y\rangle\leq L|x-y|^2, \forall x,y \in \mathbb{R}^d,)$$
 and for some continuous, monotone increasing function $f:[0,\infty)\rightarrow[1,\infty)$ with $\int_0^\infty \frac{dr}{f(r)}=\infty$ we have that $$|b(x)|\leq f(|x|), x\in \mathbb{R}^d.$$
 Then in [23, Subsection II.2] it is proved that there exists a generalized Dirichlet form in $L^2(\mathbb{R}^d)$ extending $\mathcal{E}_c$. We denote the semigroup associated with $\mathcal{E}_c$ by $P_t^c$. If we define $P_t:=e^{ct}P_t^c$, then it is the semigroup associated with $\mathcal{E}$. By the computation in [23, Subsection II.2], $P_t$ is sub-Markovian. So it satisfies the conditions (A1)-(A4).
 \vskip.10in

Further examples  are presented in Section 4 (see Examples 4.2 and 4.3) and Section 6.

Then we use the  same notation $\hat{F}, \mathcal{C}_T,\|\cdot\|_T$ associated with $\mathcal{E}^{a,\hat{b}}$ as in \cite{BPS05}:
$\mathcal{C}_T=C^1([0,T];L^2)\cap L^2([0,T];F)$, which turns out to be the appropriate space of test functions, i.e.
$$\aligned\mathcal{C}_T=\{&\varphi:[0,T]\times \mathbb{R}^d\rightarrow \mathbb{R}|\varphi_t\in F \textrm{ for almost each } t, \int_0^T\mathcal{E}^{a,\hat{b}}(\varphi_t)dt<\infty,\\&t\rightarrow\varphi_t \textrm{ is differentiable in } L^2 \textrm {and } t\rightarrow\partial_t\varphi_t \textrm{ is } L^2-\textrm{continuous on }[0,T]\}.\endaligned$$
Here and below we set $\mathcal{E}^{a,\hat{b}}(u)$ for $\mathcal{E}^{a,\hat{b}}(u,u)$.
We also set $\mathcal{C}_{[a,b]}=C^1([a,b];L^2)\cap L^2([a,b];F)$.
For $\varphi\in\mathcal{C}_T$, we define
$$\|\varphi\|_T:=(\sup_{t\leq T}\|\varphi_t\|_2^2+\int_0^T\mathcal{E}^{a,\hat{b}}_{c_2}(\varphi_t)dt)^{1/2}.$$
$\hat{F}$ is the completion of $\mathcal{C}_T$ with respect to $\|\cdot\|_T$. By \cite{BPS05}, $\hat{F}=C([0,T];L^2)\cap L^2(0,T;F)$. We define the space $\hat{F}^a$ w.r.t. $\mathcal{E}_1^a$ analogous to $\hat{F}$. Then we have $\hat{F}\subset \hat{F}^a$. We also introduce the following space
$$W^{1,2}([0,T];L^2(\mathbb{R}^d))=\{u\in L^2([0,T];L^2);\partial_t u\in L^2([0,T];L^2)\},$$
where $\partial_t u$ is the derivative of $u$ in the weak sense (see e.g. \cite{Ba10}).

\vskip.10in

\subsection{Linear Equations}

Consider the linear equation
\begin{equation}\aligned (\partial_t+L)u+f&=0,\qquad  0\leq t\leq T,\\
u_T(x)&=\phi(x), \qquad x\in \mathbb{R}^d,\endaligned\end{equation}
where $f\in L^1([0,T];L^2),\phi\in L^2$.

As in \cite{BPS05} we set $D_\sigma \varphi:=(\nabla \varphi) \sigma$ for any $\varphi\in C_0^\infty(\mathbb{R}^d)$, define $V_0=\{D_\sigma\varphi:\varphi\in C_0^\infty(\mathbb{R}^d)\}$, and let $V$ be the closure of $V_0$ in $L^2(\mathbb{R}^d;\mathbb{R}^k)$. Then we have the following results:
\vskip.10in
\th{Proposition 2.3} Assume (A1)-(A3) hold. Then:

(i) For every $u\in F^a$ there is a unique element of $V$, which we denote by $D_\sigma u$, such that
$$\mathcal{E}^a(u)=\int\langle D_\sigma u(x),D_\sigma u(x)\rangle m(dx).$$

(ii) Furthermore, if $u\in \hat{F}^a$, then there exists a measurable function $\phi:[0,T]\times \mathbb{R}^d\mapsto \mathbb{R}^d$ such that $|\phi \sigma|\in L^2((0,T)\times \mathbb{R}^d)$ and $D_\sigma u_t=\phi_t\sigma$ for almost all $t\in[0,T]$.

(iii)Let $u^n,u\in \hat{F}^a$ be such that $u^n\rightarrow u$ in $L^2((0,T)\times \mathbb{R}^d)$ and $(D_\sigma u^n)_n$ is Cauchy in $L^2([0,T]\times \mathbb{R}^d; \mathbb{R}^k)$. Then $D_\sigma u^n\rightarrow D_\sigma u$ in $L^2((0,T)\times \mathbb{R}^d;\mathbb{R}^k)$, i.e. $D_\sigma$ is closed as an operator from $\hat{F}^a$ into $L^2((0,T)\times \mathbb{R}^d)$.

\proof See [2, Proposition 2.3].
$\hfill\Box$

\vskip.10in
For $u\in F,v\in bF$, we define $$\mathcal{E}(u,v):=\mathcal{E}^{a,\hat{b}}(u,v)+\int \langle b\sigma, D_\sigma u\rangle v m(dx).$$

\vskip.10in
\th{Notation}We denote by $\tilde{\nabla}u$ the set of all measurable functions $\phi:\mathbb{R}^d\rightarrow \mathbb{R}^d$ such that $\phi\sigma=D_\sigma u$ as elements of $L^2(\mathbb{R}^d,\mathbb{R}^k)$.

\subsection{Solution of the Linear Equation}
\vskip.10in
We recall the following standard notions.

\th{Definition 2.4} (\emph{strong solutions}) A function $u\in \hat{F}\cap L^1((0,T);\mathcal{D}(L))$ is called a strong solution of equation (2.5) with data $\phi,f$, if $t\mapsto u_t=u(t,\cdot)$ is $L^2$-differentiable on $[0,T],\partial _tu_t\in L^1((0,T);L^2)$ and the equalities in (2.5) hold $m$-a.e..
\vskip.10in
\th{Definition 2.5}(\emph{generalized solutions}) A function $u\in \hat{F}$ is called a generalized solution of equation (2.5), if there are sequences  $\{u^n\}$ which are strong solutions with data $(\phi^n,f^n) $ such that $$\|u^n-u\|_T\rightarrow0,\|\phi^n-\phi\|_2\rightarrow0,\lim_{n\rightarrow\infty} f^n=f \textrm{ in } L^1([0,T];L^2).$$
\vskip.10in

\th{Remark}  Here we use the definition of a generalized solution and prove it is equivalent to a mild solution in Proposition 2.9.  If we use the definition of  weak solution to define our solution as in [2],  uniqueness of the solution in Proposition 2.9 can not be obtained since only $|b\sigma|\in L^2(\mathbb{R}^d;m)$.
\vskip.10in

\th{Proposition 2.6} Assume (A3)-(A4) hold.

(i) Let $f\in C^1([0,T];L^p)$ for $p\in[1,\infty)$. Then
$w_t:=\int_t^TP_{s-t}f_sds\in C^1([0,T];L^p),$
 and
$\partial_t w_t=-P_{T-t}f_T+\int_t^TP_{s-t}\partial_sf_sds.$

(ii) Assume that $\phi\in \mathcal{D}(L)$, $f\in C^1([0,T];L^2)$ and for each $t\in [0,T]$, $f_t\in \mathcal{D}(L)$ . Define
$u_t:=P_{T-t}\phi+\int_t^TP_{s-t}f_sds.$
Then $u$ is a strong solution of (2.5) and, moreover, $u\in C^1([0,T];L^2)$.

\proof By the same arguments as in [2, Proposition 2.6].  $\hfill\Box$

\vskip.10in
\th{Remark 2.7} Here in (ii) we add the assumption
 $\phi\in \mathcal{D}(L)$ and   $f_t\in \mathcal{D}(L)$, $t\in [0,T]$, as  we can not deduce $P_t\phi\in \mathcal{D}(L)$ for $\phi\in L^2$, since $(P_t)$ might not be analytic.

\vskip.10in
\th{Proposition 2.8} Suppose (A4) holds. If $u$ is a strong solution for (2.5), it is a mild solution for (2.5) i.e.
$u_t=P_{T-t}\phi+\int_t^TP_{s-t}f_sds.$

\proof For fixed $t$, $\varphi\in \mathcal{D}(\hat{L})$, 
$(u_T,\hat{P}_{T-t}\varphi)-(u_t,\varphi)=\int_t^T(-Lu_s-f_s,\hat{P}_{s-t}\varphi)ds+\int_t^T(u_s,\hat{L}\hat{P}_{s-t}\varphi)ds.$
Here $\hat{L}, \hat{P}_t$ denote the adjoints on $L^2(\mathbb{R}^d,m)$ of $L, P_t$ respectively. As $u$ is a strong solution, we can deduce that
$(u_t,\varphi)=(P_{T-t}\phi+\int_t^TP_{s-t}f_sds,\varphi).$
Since $\mathcal{D}(\hat{L})$ is dense in $L^2$, the result follows.
$\hfill\Box$

\vskip.10in
\th{Proposition 2.9} Assume that conditions (A1)-(A4) hold, $f\in L^1([0,T];L^2)$ and $\phi\in L^2$. Then equation (2.5) has a unique generalized solution  $u\in \hat{F}$ and
\begin{equation} u_t=P_{T-t}\phi+\int_t^TP_{s-t}f_sds.\end{equation}
The solution satisfies the three relations:
\begin{equation}\|u_t\|_2^2+2\int_t^T\mathcal{E}^{a,\hat{b}}(u_s)ds\leq 2\int_t^T(f_s,u_s)ds+\|\phi\|_2^2+2\alpha\int_t^T\|u_s\|_2^2ds, \qquad 0\leq t\leq T.\end{equation}
\begin{equation}\|u\|_T^2\leq M_T(\|\phi\|_2^2+(\int_0^T\|f_t\|_2dt)^2).\end{equation}
\begin{equation}\int_{t_0}^T((u_t,\partial_t\varphi_t)+\mathcal{E}^{a,\hat{b}}(u_t,\varphi_t)+\int\langle b\sigma, D_\sigma u_t
\rangle \varphi_tdm)dt=\int_{t_0}^T(f_t,\varphi_t)dt+(\phi,\varphi_T)-(u_{t_0},\varphi_{t_0}),\end{equation}
for any $\varphi\in b \mathcal{C}_T, t_0\in[0,T]$. $M_T$ is a constant depending on $T$. (2.9) can be extended easily for $\varphi\in bW^{1,2}([0,T];L^2)\cap L^2([0,T];F)$.

Moreover, if $u\in\hat{F}$ is bounded and satisfies (2.9) for any $\varphi\in b \mathcal{C}_T$ with bounded $f,\phi$, then $u$ is a generalized solution given by (2.6).

\proof [Existence] Define $u$ by (2.6).
First assume that $\phi,f$ are bounded and satisfy the conditions of Proposition 2.6 (ii). Then, since $u$ is bounded and by Proposition 2.6 we know that $u$ is a strong solution of (2.5), hence it obviously satisfies (2.9). Furthermore, $u\in C^1([0,T];L^2)$.  Hence, actually $u\in b\mathcal{C}_T$ and consequently,
 $$\int_{t}^T((u_s,\partial_tu_s)+\mathcal{E}^{a,\hat{b}}(u_s,u_s)+\int\langle b\sigma, D_\sigma u_s
\rangle u_sdm)ds=\int_{t}^T(f_s,u_s)ds+(\phi,u_T)-(u_{t},u_{t}).$$
By (2.3) we have $\int\langle b\sigma, D_\sigma u_s
\rangle u_sdm\geq-\alpha\|u_s\|_2^2$. Hence
$$\|u_t\|_2^2+2\int_t^T\mathcal{E}^{a,\hat{b}}(u_s)ds\leq 2\int_t^T(f_s,u_s)ds+\|\phi\|_2^2+2\alpha\int_t^T\|u_s\|_2^2ds, \qquad 0\leq t\leq T.$$
As
$$\aligned\int_t^T(f_s,u_s)ds=&\int_t^T((f_s,P_{T-s}\phi)+(f_s,\int_s^TP_{r-s}f_rdr))ds
\\ \leq & M_0e^{T-t}(\|\phi\|_2\int_t^T\|f_s\|_2ds+\int_t^T(\|f_s\|_2\int_s^T\|f_r\|_2dr)ds),\endaligned$$
and
$\int_t^T\|u_s\|_2^2ds\leq M_{T-t}(\|\phi\|_2^2+(\int_0^T\|f_t\|_2dt)^2),$
 we obtain, 
$\|u_t\|_2^2+\int_t^T\mathcal{E}^{a,\hat{b}}(u_s)ds\leq M_{T-t}(\|\phi\|_2^2+(\int_0^T\|f_t\|_2dt)^2).$
Hence,  it follows that
\begin{equation}\|u\|_T^2\leq M_T(\|\phi\|_2^2+(\int_0^T\|f_t\|_2dt)^2).\end{equation}
Here $M_{T-t}$ can change from line to line and is independent of $f,\phi$.
Now we will obtain the result for general data $\phi$ and $f$. Let $(f^n)_{n\in N}$ be a sequence of bounded function in $C^1([0,T];L^2)$ such that $f_t\in \mathcal{D}(L)$ for a.e. $t\in [0,T]$ and $\int_0^T\|f_t^n-f_t\|_2dt\rightarrow0$ (This sequence can be obtained since $\{\alpha_t g(x);\alpha_t\in C_0^\infty[0,T], g\in b\mathcal{D}(L)\}$ is dense in  $L^1([0,T];L^2)$). Take bounded functions $(\phi^n)_{n\in N}\subset \mathcal{D}(L)$ such that $\phi^n\rightarrow \phi$ in $L^2$. Let $u^n$ denote
the solution given by (2.6) with $f=f^n, \phi=\phi^n$.

By linearity, $u^n-u^m$ is associated with $(\phi^n-\phi^m, f^n-f^m)$. Since by (2.10)
$$\|u^n-u^m\|_T^2\leq M_T(\|\phi^n-\phi^m\|_2^2+(\int_0^T\|f_t^n-f_t^m\|_2dt)^2),$$
we can deduce that $(u^n)_{n\in N}$ is a Cauchy sequence in $\hat{F}$. Hence $u=\lim_{n\rightarrow\infty}u^n$ in $\|\cdot\|_T$ is the generalized solution of (2.5) and (2.6) follows.

 Next we prove (2.7)(2.8) (2.9)  for $u$. We have (2.9) for $u^n$ with $f^n,\phi^n$ and $\varphi\in b\mathcal{C}_T$, i.e.
 $$\int_0^T((u_t^n,\partial_t\varphi_t)+\mathcal{E}^{a,\hat{b}}(u_t^n,\varphi_t)+\int\langle b\sigma, D_\sigma u_t^n
\rangle \varphi_tdm)dt=\int_0^T(f_t^n,\varphi_t)dt+(\phi^n,\varphi_T)-(u_0^n,\varphi_0).$$
We have
$|\int_0^T\mathcal{E}^{a,\hat{b}}(u_t^n-u_t,\varphi_t)dt|\rightarrow0,$
and
 $$\aligned |\int_0^T\int\langle b\sigma, D_\sigma (u_t^n-u_t)
\rangle \varphi_tdmdt|&\leq \|\varphi\|_\infty (\int_0^T\int|b\sigma|^2dmdt)^{\frac{1}{2}}(\int_0^T\int|D_\sigma (u_t^n-u_t)|^2dmdt)^{\frac{1}{2}}
\rightarrow0.\endaligned$$
Hence we deduce (2.9) for any $\varphi\in b \mathcal{C}_T.$

Since $\|u_t^n\|_T\rightarrow\|u_t\|_T$, we conclude
$\lim_{n\rightarrow\infty}\int_0^T\mathcal{E}^{a,\hat{b}}(u_t^n)dt=\int_0^T\mathcal{E}^{a,\hat{b}}(u_t)dt.$ As the relations (2.7), (2.8) hold for the approximating functions, by passing to the limit,  (2.7) and (2.8) follows for $u$.

[Uniqueness] Let $v\in \hat{F}$ be another generalized solution of (2.5) and $(v^n)_{n\in N},(\tilde{\phi}^n)_{n\in N}, (\tilde{f}^n)_{n\in N}$ be the corresponding approximating sequences in the definition of  generalized solutions.
By Proposition 2.8,  $\sup_{t\in[0,T]}\|u^n_t-v^n_t\|_2^2\leq M_T(\|\phi^n-\tilde{\phi}^n\|_2^2+(\int_0^T\|f_t^n-\tilde{f}_t^n\|_2dt)^2).$
Letting $n\rightarrow\infty$, this implies $u=v$.

For the last result  we note that $\forall t_0\geq 0,\varphi\in b\mathcal{C}_T$
\begin{equation}\int_{t_0}^T((u_t,\partial_t\varphi_t)+\mathcal{E}^{a,\hat{b}}(u_t,\varphi_t)+\int\langle b\sigma, D_\sigma u_t
\rangle \varphi_tdm)dt=\int_{t_0}^T(f_t,\varphi_t)dt+(\phi,\varphi_T)-(u_{t_0},\varphi_{t_0}). \end{equation}
For $t\geq \frac{1}{n} $,  define
$u_t^n:=n\int_0^{\frac{1}{n}}u_{t-s}ds, f_t^n:=n\int_0^{\frac{1}{n}}f_{t-s}ds,\phi^n:=n\int_0^{\frac{1}{n}}u_{T-s}ds.$
Let us check that each $u^n$ also fulfills (2.11) with $f^n,\phi^n$.
We set $\varphi_r^s:=\varphi_{r+s}$ for $0\leq s+ r\leq T$. Then for fixed $t_0\in (0,T],$ and $n\geq \frac{1}{t_0}$,
$$\aligned &\int_{t_0}^T((u_t^n,\partial_t\varphi_t)+\mathcal{E}^{a,\hat{b}}(u_t^n,\varphi_t)+\int\langle b\sigma, D_\sigma u_t^n\rangle \varphi_tdm)dt
=\int_{t_0}^T(f_t^n,\varphi_t)dt+(\phi^n,\varphi_T)-(u_{t_0}^n,\varphi_{t_0})dt.\endaligned$$
For the mild solution $v$ associated with $f,\phi$, the above relation also holds with $v^n$ replacing $u^n$. Hence we have
$$\int_{t_0}^T(((u-v)_t^n,\partial_t\varphi_t)+\mathcal{E}^{a,\hat{b}}((u-v)_t^n,\varphi_t)+\int\langle b\sigma, D_\sigma (u-v)_t^n\rangle \varphi_tdm)dt=-((u-v)_{t_0}^n,\varphi_{t_0}).$$
And we have $(u-v)_t^n\in b\mathcal{C}_{[\frac{1}{n},T]}$.  Hence, the above equation holds with $(u-v)_t^n$ as a test function, i.e. for $n\geq\frac{1}{t_0}$
$$\int_{t_0}^T(((u-v)_t^n,\partial_t(u-v)_t^n)+\mathcal{E}^{a,\hat{b}}((u-v)_t^n,
(u-v)_t^n)+\int\langle b\sigma, D_\sigma (u-v)_t^n
\rangle (u-v)_t^n dm)dt=-((u-v)_{t_0}^n,(u-v)_{t_0}^n).$$
So we have $\|(u-v)^n_{t_0}\|_2^2+2\int_t^T\mathcal{E}^{a,\hat{b}}((u-v)_t^n,
(u-v)_t^n)dt\leq 2\alpha\int_{t_0}^T\|(u-v)^n_t\|_2^2dt.$
By Gronwall's Lemma it follows that $\|(u-v)^n_{t_0}\|_2^2=0.$
Letting $n\rightarrow\infty$, we have $\|u_{t_0}-v_{t_0}\|_2=0$. Then letting $t_0\rightarrow0$, we have $\|u_0-v_0\|=0$. Then
$u_t=P_{T-t}\phi+\int_t^TP_{s-t}f_sds$ is a generalized solution for (2.5).
$\hfill\Box$
\vskip.10in

 We can prove the following basic relations for the linear equation which is essential to the following section. The basic idea of the proof comes from [2]. As the definition of the solution is different from [2] and our bilinear form is not symmetric, we need to apply some results in Proposition 2.9 and some  properties of the bilinear form $\mathcal{E}$ from (A1)-(A4) to conclude the following proposition. And here  we also use some new estimates to deal with the non-symmetric part  of the bilinear form $\mathcal{E}$. The proof is included in Appendix C.
\vskip.10in

\th{Proposition 2.10} Assume (A1)-(A4) hold. Let $u=(u^1,...,u^l)$ be a vector-valued function, where each component is a generalized solution of the linear equation (2.5) associated to certain data $f^i\in L^1([0,T];L^2),\phi^i\in L^2$    for $i=1,...,l$. Denote by $\phi,f$ the vectors $\phi=(\phi^1,...,\phi^l),f=(f^1,...,f^l)$ and by $D_\sigma u$ the matrix whose rows consist of the row vectors $D_\sigma u^i$. Then the following relations hold $m$-almost everywhere
\begin{equation}
 |u_t|^2+2\int_t^TP_{s-t}(|D_\sigma u|^2+\frac{1}{2}c|u_s|^2)ds=P_{T-t}|\phi|^2+2\int_t^TP_{s-t}\langle u_s,f_s\rangle ds.
\end{equation}
\begin{equation}
 |u_t|\leq P_{T-t}|\phi|+\int_t^TP_{s-t}\langle \hat{u}_s,f_s\rangle ds.
\end{equation}

\section{The Non-linear Equation}
In the case of non-linear equations, we are going to consider systems of equations, with the unknown functions and their first-order derivatives mixed in the non-linear term of the equation. The non-linear term is a given measurable function $f:[0,T]\times \mathbb{R}^d\times \mathbb{R}^l\times \mathbb{R}^l\otimes \mathbb{R}^k\rightarrow \mathbb{R}^l$, $l\in N$.
We are going to treat the following system of equations.
\begin{equation}
 (\partial_t+L)u+f(\cdot,\cdot,u, D_\sigma u)=0 \qquad u_T=\phi.
\end{equation}
 Here $\phi\in L^2(\mathbb{R}^d,dm;\mathbb{R}^l)$.
\vskip.10in
\th{Definition 3.1} (\emph{Generalized solutions of the nonlinear equation})
A generalized solution of equation (3.1) is a system $u=(u^1,u^2,...,u^l)$ of $l$ elements in $\hat{F}$ with the property that $f^i(\cdot,\cdot, u, D_\sigma u)$ belongs to $L^1([0,T];L^2)$ and
there are sequences  $\{u_n\}$ which are strong solutions of (3.1) with data $(\phi_n,f_n) $ such that $$\|u_n-u\|_T\rightarrow0,\|\phi_n-\phi\|_2\rightarrow0,\textrm{ and }\lim_{n\rightarrow\infty} f_n(\cdot,\cdot,u_n,D_\sigma u_n)=f(\cdot,\cdot,u,D_\sigma u) \textrm{ in } L^1([0,T];L^2).$$
\vskip.10in
\th{Definition 3.2} (\emph{Mild equation})
For every $i\in \{1,...,l\}$ we define the mild equation as
\begin{equation}u^i(t,x)=P_{T-t}\phi^i(x)+\int_t^TP_{s-t}f^i(s,\cdot,u_s,D_\sigma u_s)(x)ds, m-a.e..\end{equation}
\vskip.10in
\th{Lemma 3.3} $u$ is a generalized solution of the nonlinear equation (3.1) if and only if it solves the mild equation (3.2).

\proof The assertion follows by Proposition 2.9.$\hfill\Box$
\vskip.10in
We will use the following notation
$\|\phi\|_2^2=\sum_{i=1}^l\|\phi^i\|_2^2, \phi\in L^2(\mathbb{R}^d;\mathbb{R}^l),$
$\mathcal{E}(u,v)=\sum_{i=1}^l\mathcal{E}(u^i,v^i)$, 

$\mathcal{E}^a(u,v)=\sum_{i=1}^l\mathcal{E}^a(u^i,v^i),u,v\in F^l,\|u\|_T^2:=\sup_{t\leq T}\|u_t\|_2^2+\int_0^T\mathcal{E}^{a,\hat{b}}_{c_2+1}(u_t)dt, u\in \hat{F}^l.$
\subsection{The Case of Lipschitz Conditions}
In this subsection we consider a measurable function $f:[0,T]\times \mathbb{R}^d\times \mathbb{R}^l\times \mathbb{R}^l\otimes \mathbb{R}^k\rightarrow \mathbb{R}^l$ such that
\begin{equation}|f(t,x,y,z)-f(t,x,y',z')|\leq C(|y-y'|+|z-z'|),\end{equation}
with $t,x,y,y',z,z'$ arbitrary and $C$ is a constant independent of $t,x$. Set $f^0(t,x):=f(t,x,0,0)$.
\vskip.10in
\th{Proposition 3.4} Assume that the conditions (A1)-(A4) hold and that  $f$ satisfies condition (3.3), $f^0\in L^2([0,T]\times \mathbb{R}^d,dt\times dm;\mathbb{R}^l)$ and $\phi\in L^2(\mathbb{R}^d;\mathbb{R}^l)$. Then
the equation (3.1) admits a unique generalized solution $u\in\hat{F}^l$ and it satisfies the following estimate
$$\|u\|_T^2\leq e^{T(1+2C+\frac{C^2}{c_1}+2\alpha+c_2)}(\|\phi\|_2^2+\|f^0\|^2_{L^2([0,T]\times \mathbb{R}^d)}).$$

\proof If $u\in \hat{F}^l$, then by relation (3.3) we have
$$\aligned |f(\cdot,\cdot,u,D_\sigma u)|&\leq |f(\cdot,\cdot,u,D_\sigma u)-f(\cdot,\cdot,0,0)|+|f(\cdot,\cdot,0,0)|
\\& \leq C(|u|+|D_\sigma u|)+|f^0|.\endaligned$$
As $f^0 \in L^2([0,T]\times \mathbb{R}^d,dt\times dm;\mathbb{R}^l)$ and $|D_\sigma u|$ is an element of $L^2([0,T]\times \mathbb{R}^d)$, we get $f(\cdot,\cdot,u,D_\sigma u)\in  L^2([0,T]\times \mathbb{R}^d;\mathbb{R}^l)$.

Now we define the operator $A:\hat{F}^l\rightarrow \hat{F}^l$ by
$(Au)^i_t:=P_{T-t}\phi^i+\int_t^T P_{s-t}f^i(s,\cdot,u_s,D_\sigma u_s)ds, i=1,...,l.$
Then Proposition 2.9 implies that $Au\in \hat{F}^l$. In the following we write $f_{u,s}^i:=f^i(s,\cdot,u_s,D_\sigma u_s).$ Since $(Au)^i_t-(Av)^i_t=\int_t^TP_{s-t}(f_{u,s}^i-f_{v,s}^i)ds$ is
 the mild solution with data $(f_u^i-f_v^i,0)$, by the same arguments as in Proposition 2.9  we have
$$\aligned&\|\int_t^TP_{s-t}(f_{u,s}^i-f_{v,s}^i)ds\|_{[t,T]}^2\leq M_T(\int_t^T\|f_{u,s}-f_{v,s}\|_2ds)^2
\\\leq&  M_T(T-t) \int_t^T(\|u_s-v_s\|_2^2+\|D_\sigma u_s-D_\sigma v_s\|^2_2) ds
 \leq  M_T(T-t) \|u-v\|_{[t,T]}^2,\endaligned$$
where $M_T$ can change from line to line.
Here $\|u\|_{[T_a,T_b]}:=(\sup_{t\in[T_a,T_b]}\|u_t\|_2^2+\int_{T_a}^{T_b}\mathcal{E}^{a,\hat{b}}_{c_2+1}(u_t)dt)^{\frac{1}{2}},$
where $0\leq T_a\leq T_b\leq T$. Fix $T_1$ sufficiently small such that $\gamma:=M_T(T-T_1)<1$. Then we have
:$\|Au-Av\|_{[T_1,T]}^2\leq \gamma\|u-v\|^2_{[T_1,T]}.$
Then there exists a unique $u_1\in \hat{F}_{[T_1,T]}$ such that $Au_1=u_1$ where $\hat{F}_{[T_a,T_b]}:=C([T_a,T_b];L^2)\cap L^2((T_a,T_b);F)$ for $T_a\in [0,T]$ and $T_b\in[T_a,T]$.

 We define the operator $A^1:\hat{F}^l\rightarrow \hat{F}^l$ by
$(A^1u)^i_t:=P_{T_1-t}u_{1,T_1}^i+\int_t^{T_1} P_{s-t}f^i(s,\cdot,u_s,D_\sigma u_s)ds, i=1,...,l.$
Then by the same method as above,  we get
$\|A^1u-A^1v\|_{[t,T_1]}^2\leq M_T(T_1-t)\|u-v\|^2_{[t,T_1]}.$
Now we choose $T_2<T_1$ such that $M_T(T_1-t)<1$. Hence
there exists a unique $u_2\in \hat{F}_{[T_2,T_1]}$ such that $A^1u_2=u_2$. Define $u:=u_11_{[T_1,T]}+u_21_{[T_2,T_1)}$. Then we have
if $T_2\leq t\leq T_1$, that $$\aligned&P_{T-t}\phi^i+\int_t^{T} P_{s-t}f^i(s,\cdot,u_s,D_\sigma u_{s})ds
\\&=P_{T-t}\phi^i+\int_t^{T_1} P_{s-t}f^i(s,\cdot,u_s,D_\sigma u_s)ds+P_{T_1-t}(u_{1,T_1}^i-P_{T-T_1}\phi^i)
=u_{2,t}^i.\endaligned$$ If $t> T_1$,
$(Au)_t^i=
P_{T-t}\phi^i+\int_t^{T} P_{s-t}f^i(s,\cdot,u_{1,s},D_\sigma u_{1,s})ds
=u_{1,t}^i.$
Therefore, we can construct a solution over the interval $[T_2,T]$. Clearly there exists $n\in N$ such that $T<n(T-T_1)$. Hence, the construction is done after $n$ steps.

In order to obtain the estimate in the statement, we write
$$\aligned |\int_t^T(f_{u,s},u_s)ds|
\leq &\int_t^T\|f_s^0\|_2\|u_s\|_2ds+C\int_t^T\|u_s\|_2^2ds+C\int_t^T\|D_\sigma u_s\|_2\|u_s\|_2ds
\\ \leq &\frac{1}{2}\int_t^T\|f_s^0\|_2^2ds+(\frac{1}{2}+C+\frac{1}{2c_1}C^2)\int_t^T\|u_s\|_2^2ds+\frac{c_1}{2}\int_t^T\mathcal{E}^a(u_s)ds.\endaligned$$
By relation (2.7) of Proposition 2.9 it follows that
$$\aligned&\|u_t\|_2^2+2\int_t^T\mathcal{E}^{a,\hat{b}}(u_s)ds\leq2\int_t^T(f_{u,s},u_s)ds+\|\phi\|_2^2+2\alpha\int_t^T\|u_s\|_2^2ds\\
\leq&\|\phi\|_2^2+\int_t^T\|f_s^0\|_2^2ds+(1+2C+\frac{C^2}{c_1}+2\alpha+c_2)\int_t^T\|u_s\|_2^2ds+\int_t^T\mathcal{E}^{a,\hat{b}}(u_s)ds.
\endaligned$$
Now by Gronwall's lemma the desired estimate follows.

[Uniqueness] Let $u_1$ and $u_2$ be two solutions of equation (3.1). By using (2.7) for the difference $u_1-u_2$ we get
$$\aligned&\|u_{1,t}-u_{2,t}\|_2^2+2\int_t^T\mathcal{E}^{a,\hat{b}}(u_{1,s}-u_{2,s})ds\\\leq& 2\int_t^T(f(s,\cdot,u_{1,s},D_\sigma u_{1,s})-f(s,\cdot,u_{2,s},D_\sigma u_{2,s}),u_{1,s}-u_{2,s})ds+2\alpha\int _t^T\|u_{1,s}-u_{2,s}\|_2^2ds\\
\leq &2\int_t^T C(|D_\sigma u_{1,s}-D_\sigma u_{2,s}|,|u_{1,s}-u_{2,s}|)ds+(2\alpha+C)\int _t^T\|u_{1,s}-u_{2,s}\|_2^2ds
\\ \leq & (\frac{C^2}{c_1}+c_2+2\alpha+C) \int _t^T\|u_{1,s}-u_{2,s}\|_2^2ds+\int_t^T\mathcal{E}^{a,\hat{b}}(u_{1,s}-u_{2,s})ds.\endaligned$$
By Gronwall's lemma it follows that $u_1=u_2$.
$\hfill\Box$
\vskip.10in
\subsection{The Case of Monotonicity Conditions}
Let $f:[0,T]\times \mathbb{R}^d\times \mathbb{R}^l\times \mathbb{R}^l\otimes \mathbb{R}^k\rightarrow \mathbb{R}^l$ be a measurable function and $\phi\in L^2(\mathbb{R}^d,m;\mathbb{R}^l)$ be the final condition of
 (3.1). We impose the following conditions:

\noindent(H1) (\emph{Lipschitz condition in $z$}) There exists a fixed constant $C>0$ such that for $t,x,y,z,z'$ arbitrary, $|f(t,x,y,z)-f(t,x,y,z')|\leq C|z-z'|.$

\noindent(H2) (\emph{Monotonicity condition in $y$})
For $x,y,y',z$ arbitrary, there exists a function $\mu_t\in L^1([0,T];\mathbb{R})$ such that
$\langle y-y',f(t,x,y,z)-f(t,x,y',z)\rangle
\leq \mu_t|y-y'|^2,$
and $ \alpha_t:=\int_0^t\mu_sds.$

\noindent(H3) (\emph{Continuity condition in $y$})
For $t,x$ and $z$ fixed, the map
$y\mapsto f(t,x,y,z)$ is continuous.

We need the following notation:
$f^0(t,x):=f(t,x,0,0), f'(t,x,y):=f(t,x,y,0)-f(t,x,0,0),$
$f^{',r}(t,x):=\sup_{|y|\leq r}|f'(t,x,y)|.$

\noindent(H4)For each $r>0$, $f^{',r}\in L^1([0,T];L^2).$

\noindent(H5) $\|\phi\|_\infty<\infty,\|f^0\|_\infty<\infty,|\phi|\in L^2, |f^0|\in L^2([0,T];L^2).$

If $m(\mathbb{R}^d)<\infty$ the last two conditions in (H5) are ensured by the boundedness of $\phi$ and $f^0$. The condition (H1), (H4), and (H5) imply that if $u\in b\hat{F}$ , then $|f(u,D_\sigma u)|\in L^1([0,T];L^2)$.  Under the above conditions, even if $E$ is equal to a Hilbert space, it seems impossible to apply general monotonicity methods to the map $\mathcal{V}\ni u\mapsto f(t,\cdot,u(\cdot),D_\sigma u)\in \mathcal{V}'$ because of lack of a suitable reflexive Banach space $\mathcal{V}$ such that $\mathcal{V}\subset\mathcal{H}\subset \mathcal{V}'$. Therefore, also here we proceed developing a hands-on approach to prove existence and uniqueness of solutions for equation (3.1) as done in [2].

\vskip.10in
\th{Lemma 3.5} In (H2) without loss of generality we can assume that $\mu_t\equiv0$.

\proof Let us make the change $u_t^*=\exp(\alpha_t)u_t$ and $$\phi^*=\exp(\alpha_T)\phi, \qquad f_t^*(y,z)=\exp(\alpha_t) f_t(\exp(-\alpha_t)y, \exp(-\alpha_t)z)-\mu_t y,$$ for the data. Next we will prove that $u$ is a generalized solution associated to the data $(\phi,f)$ if and only if $u^*$ is a solution associated to the data $(\phi^*,f^*)$.
Hence we can write $u^{i}_t=P_{T-t}\phi^{i}(x)+\int_t^TP_{s-t}f^i(s,\cdot,u_s,D_\sigma u_s)(x)ds,$ equivalently as
$$\aligned u^{i,*}_t=&\exp(\alpha_t)P_{T-t}\phi^{i}(x)+\exp(\alpha_t)\int_t^TP_{s-t}f^i(s,\cdot,u_s,D_\sigma u_s)(x)ds\\
=&\exp(\alpha_T)P_{T-t}\phi^{i}(x)+(\exp(\alpha_t)-\exp(\alpha_T))P_{T-t}\phi^{i}(x)\\&+\int_t^T(\exp(\alpha_t)-\exp(\alpha_s)) P_{s-t}f^i(s,\cdot,u_s,D_\sigma u_s)(x)ds+\int_t^T\exp(\alpha_s)P_{s-t}f^i(s,\cdot,u_s,D_\sigma u_s)(x)ds\\
=&P_{T-t}\phi^{i,*}(x)+\int_t^TP_{s-t}(\exp(\alpha_s)f^{i}(s,\cdot,\exp(\alpha_s)u_s^*,\exp(\alpha_s)D_\sigma u_s^*)(x))ds\\&-\int_t^TP_{s-t}(\mu_s\exp(\alpha_s)P_{T-s}\phi^i(x)+\int_s^T\mu_s \exp(\alpha_s)P_{l-s}f^i(l,\cdot,u_l,D_\sigma u_l)dl)ds\\=
& P_{T-t}\phi^{i,*}(x)+\int_t^TP_{s-t}f^{i,*}(s,\cdot,u_s^*,D_\sigma u_s^*)(x)ds.\endaligned$$
Next we prove $f^*$ satisfies (H1)-(H5). It is obvious that (H1)-(H5) are satisfied. $\hfill\Box$
\vskip.10in

\th{Lemma 3.6} Assume  that conditions (A1)-(A4), (H1) and the following weaker form of condition (H2) (with $\mu_t\equiv0$) hold:
$$\textrm{(H2')   } \langle y, f'(t,x,y)\rangle\leq 0 \textrm{ for all }t,x,y.$$  If $u$ is a generalized solution of (3.1), then there exists a constant $K$  depending on $C,\mu_t, T,\alpha$ such that
\begin{equation}\|u\|_T^2\leq K(\|\phi\|_2^2+\int_0^T\|f_t^0\|_2^2dt).\end{equation}
\proof Since $u$ is a solution of (3.1),  by Proposition 2.9 we have
$\|u_t\|_2^2+2\int_t^T\mathcal{E}^{a,\hat{b}}(u_s)ds\leq 2\int_t^T(f_s,u_s)ds+\|u_T\|^2_2+2\alpha\int_t^T\|u_s\|_2^2ds.$
Conditions (H1) and (H2') yield
$$\aligned \langle f_s(u_s,D_\sigma u_s),u_s\rangle=&\langle f_s(u_s,D_\sigma u_s)-f_s(u_s,0)+f_s'(u_s)+f_s^0,u_s\rangle\\
\leq& (C|D_\sigma u_s|+|f_s^0|)|u_s|.\endaligned$$
Hence, it follows that
$$\aligned\|u_t\|_2^2+2\int_t^T\mathcal{E}^{a,\hat{b}}(u_s)ds\leq& 2\int_t^T\int(C|D_\sigma u_s|+|f_s^0|)|u_s|dmds+\|u_T\|_2^2+2\alpha\int_t^T\|u_s\|_2^2ds
\\\leq &\int_t^T\mathcal{E}^{a,\hat{b}}(u_s)ds+(\frac{C^2}{c_1} +1+2\alpha+c_2)\int_t^T\|u_s\|_2^2ds+\int_t^T\|f_s^0\|_2^2ds+\|u_T\|^2_2.\endaligned$$
Then by Gronwall's lemma, the assertion follows.$\hfill\Box$
\vskip.10in

By a modification of the arguments  in [2, Lemma 3.3] we have the following estimates.
\vskip.10in

\th{Lemma 3.7} Assume  that conditions (A1)-(A4), (H1) and (H2') hold. If $u$ is a generalized solution of (3.1),
 there exists a constant $K$, which depends on $C,\mu$ and $T$, such that
\begin{equation}
 \|u\|_\infty\leq K(\|\phi\|_\infty+\|f^0\|_\infty).
\end{equation}
\proof By Proposition 2.10 we have
\begin{equation}
 |u_t|^2+2\int_t^TP_{s-t}(|D_\sigma u|^2)\leq P_{T-t}|\phi|^2+2\int_t^TP_{s-t}\langle u_s,f_s(u_s,D_\sigma u_s)\rangle ds.
\end{equation}
Following the same arguments as in the proof of  Lemma 3.6  we deduce
$ \langle f_s(u_s,D_\sigma u_s),u_s\rangle
\leq (C|D_\sigma u_s|+|f_s^0|)|u_s|.$
And by Proposition 2.10 (2.13) we get
$
 |u_s|\leq P_{T-s}|\phi|+\int_s^TP_{r-s}(C|D_\sigma u_r|+|f_r^0|) dr.
$
Then we have
$$\aligned&\int_t^TP_{s-t}\langle f_s(u_s,D_\sigma u_s),u_s\rangle ds\\ \leq&\int_t^TP_{s-t}[(P_{T-s}|\phi|+\int_s^TP_{r-s}(C|D_\sigma u_r|+|f_r^0|) dr)(C|D_\sigma u_s|+|f_s^0|)]ds.\endaligned$$
So, by (3.6) and Lemma C.6 we obtain
$$\aligned &|u_t|^2+2\int_t^TP_{s-t}(|D_\sigma u_s |^2)ds\\
\leq & P_{T-t}|\phi|^2+2(\int_t^TP_{s-t}[(P_{T-s}|\phi|+\int_s^TP_{r-s}(C|D_\sigma u_r|+|f_r^0|) dr)(C|D_\sigma u_s|+|f_s^0|)]ds)
\\ \leq &3 P_{T-t}|\phi|^2+2C^2\int_t^T\int_s^TP_{s-t}(|D_\sigma u_s|P_{r-s}|D_\sigma u_r|)drds+2\int_t^T\int_s^TP_{s-t}(|f_s^0|P_{r-s}|f_r^0|)drds\\
&+2\int_t^T\int_s^TP_{s-t}[P_{r-s}(C|D_\sigma u_r|+|f_r^0|)(C|D_\sigma u_s|+|f_s^0|)]drds.\endaligned$$
Since
$$\aligned&\int_t^T\int_s^TP_{s-t}[P_{r-s}(C|D_\sigma u_r|+|f_r^0|)(C|D_\sigma u_s|+|f_s^0|)]drds\\\leq&
\frac{1}{2}\int_t^T\int_s^T[P_{s-t}(C|D_\sigma u_s|+|f_s^0|)^2]+P_{s-t}[(P_{r-s}(C|D_\sigma u_r|+|f_r^0|))^2]drds\\\leq&
\int_t^T\int_s^TC^2P_{s-t}|D_\sigma u_s|^2+P_{s-t}|f_s^0|^2+\frac{1}{2}P_{r-t}(C|D_\sigma u_r|+|f_r^0|)^2drds\\\leq&2C^2(T-t)\int_t^TP_{s-t}|D_\sigma u_s|^2ds+2(T-t)\int_t^TP_{s-t}|f_s^0|^2ds,
\endaligned$$
and since by Schwartz's inequality one has
$$\aligned&\int_t^T\int_s^TP_{s-t}(|D_\sigma u_s|P_{r-s}|D_\sigma u_r|)drds
\leq \int_t^T\int_s^T\frac{1}{2}(P_{s-t}|D_\sigma u_s|^2)drds+\int_t^T\int_s^T\frac{1}{2}(P_{r-t}|D_\sigma u_r|^2)drds\\ \leq
&(T-t)\int_t^TP_{s-t}|D_\sigma u_s|^2ds,\endaligned$$
we conclude
$$\aligned &|u_t|^2+2\int_t^TP_{s-t}(|D_\sigma u_s |^2)ds
\leq 3 P_{T-t}|\phi|^2+6C^2(T-t)\int_t^TP_{s-t}|D_\sigma u_s|^2ds+6(T-t)\int_t^TP_{s-t}|f_s^0|^2ds.
\endaligned$$
Hence, the  estimate of this lemma holds on the interval $[T-\varepsilon,T]$ where $\varepsilon>0$ such that $6 C^2\varepsilon=1$. So we can deduce by iteration the estimate over the interval $[0,T]$.
We obtain from the first estimate
$$\aligned |u_t|^2
\leq &\sup_{t\in [0,T]}\sup_{x\in \mathbb{R}^d}\tilde{K}(P_{T-t}|\phi|^2+(T-t)\int_t^TP_{s-t}|f_s^0|^2ds)\\
\leq &\sup_{t\in [0,T]}\tilde{K}(\|\phi^2\|_\infty+T^2\|f^0\|_\infty^2)\\
\leq &K^2(\|\phi\|_\infty^2+\|f^0\|_\infty^2)
,\endaligned$$
which implies (3.5).$\hfill\Box$
\vskip.10in

By the same methods as in [2, Theorem 3.2], we obtain the following results. As the method is similar as in the proof of [2, Theorem 3.2], we will give the proof  in the Appendix A.
\vskip.10in
\th{Theorem 3.8} Suppose that $m(dx)$ is a finite measure and that conditions (A1)-(A4), (H1)-(H5) hold. Then there exists a unique generalized solution of equation (3.1) and it satisfies the following estimates for some $K_1$ and $K_2$ independent of $u, \phi,f$
$$\|u\|_T^2\leq K_1(\|\phi\|_2^2+\int_0^T\|f_t^0\|_2^2dt).$$
$$ \|u\|_\infty\leq K_2(\|\phi\|_\infty+\|f^0\|_\infty).$$

\vskip.10in
\th{Lemma 3.9} Assume  conditions (A1)-(A4),(H1)-(H5) hold. If $u\in\hat{F}$ is bounded and  for $\varphi\in b\mathcal{C}_T$ satisfies $\int_0^T\mathcal{E}(u_t,\varphi_t)+(u_t,\partial_t\varphi_t) dt=\int_0^T(f_t(u_t,D_\sigma u_t),\varphi_t) dt+(u_T,\varphi_T)-(u_0,\varphi_0).$
Then we have
$$\|u_t\|_2^2+2\int_t^T\mathcal{E}^{a,\hat{b}}(u_s)ds\leq 2\int_t^T(f_s(u_s,D_\sigma u_s),u_s)ds+\|\phi\|_2^2+2\alpha\int _t^T\|u_s\|_2^2ds, \qquad 0\leq t\leq T.$$
\proof Define $u_t^h=\frac{1}{2h}\int_{t-h}^{t+h} u_sds$. Choose $\phi_t=u_t^h$, then we have
$$\int_t^T\mathcal{E}(u_t,u_t^h)+\frac{1}{2h}(u_t,u_{t+h})-\frac{1}{2h}(u_t,u_{t-h}) dt=\int_t^T(f_t(u_t,D_\sigma u_t),u_t^h) dt+(u_T,u_T^h)-(u_t,u_t^h).$$
That is to say,
\begin{equation}\aligned&\frac{1}{2h}\int_{T-h}^T(u_t,u_{t+h})dt-\frac{1}{2h}\int_{t-h}^t(u_t,u_{t+h})dt+\int_t^T\mathcal{E}(u_t,u_t^h)dt \\=&\int_0^T(f_t(u_t,D_\sigma u_t),u_t^h) dt+(u_T,u_T^h)-(u_t,u_t^h).\endaligned
\end{equation}

Letting $h\rightarrow0$ in (3.7), the assertion follows .$\hfill\Box$
\vskip.10in
For the case $m(dx)=dx$,
we will use a weight function of the form $\pi(x)=\exp[-\rho\theta(x)]$, with $\theta\in C^1(\mathbb{R}^d)$ being a fixed function such that $0\leq \theta(x)\leq |x|$, and $\theta(x)=|x|$ if $|x|\geq 1$, and $\rho\in R^+$.  If one chooses $\rho>0$, then clearly one has $m(\mathbb{R}^d)<\infty$. We  denote the generalized Dirichlet form, function spaces and the generator associated with $\rho>0$ by $\mathcal{E}^\rho$, $\hat{F}^\rho$, $\mathcal{C}_T^\rho$, $L_\rho$ respectively.
In the case $\rho=0$, we drop $\rho$ in the notation, i.e.  $\mathcal{E}=\mathcal{E}^0$. And for the case $\rho=0$, we need   the following condition.

\vskip.10in
\no(A2') (\emph{Sobolev inequality}) If $\rho=0$, then $\sigma$ is a bounded measurable field in $\mathbb{R}^d$ and
$$\|u\|_q\leq C\mathcal{E}^a(u,u)^{1/2}, \textrm{  } \forall u\in C_0^\infty(\mathbb{R}^d),$$
where $\frac{1}{q}+\frac{1}{d}=\frac{1}{2}$ and $\|\cdot\|_q$ denotes the usual norm in $L^q$.  And  $|\hat{b}\sigma|\in L^d(\mathbb{R}^d;dx)+L^\infty(\mathbb{R}^d;m)$, $c\in L^{d/2}(\mathbb{R}^d;dx)+L^\infty(\mathbb{R}^d;dx)$.
\vskip.10in

If (A2'), (A2) are satisfied, for $u,v\in bF$, we have
$$\mathcal{E}^\rho(u,v)=\int\langle D_\sigma u,D_\sigma v\rangle dm+\int cuv dm+\int\langle (b\sigma+\hat{b}\sigma,D_\sigma u\rangle v dm.$$
If $\rho=0$, we additionally have $\mathcal{E}^{a,\hat{b}}(u,u)\leq C\mathcal{E}_1^a(u,u),$
and that $F=F^a$. We also need the following condition:

\vskip.10in
\no (H6). $\mathcal{E}^a(u)<\infty,u\in L^2\Rightarrow u\in F.$
\vskip.10in

The Sobolev inequality and (H6) are satisfied if $a$ is uniformly elliptic.
By [21, Lemma 4.20] we have:

\vskip.10in
\th{Lemma 3.10} Assume conditions (A2') and (H6) hold. Let $\rho>0$ and suppose $\sigma$ is bounded. Then it holds
$\mathcal{E}^\rho(u,\varphi)=\mathcal{E}(u,\varphi\exp(-\theta \rho))+(M_\rho u,\varphi)_\rho,$
for $u\in F_\rho,\varphi\in bF_\rho,$ where $M_\rho u=\rho\langle D_\sigma \theta, D_\sigma u\rangle$.
\vskip.10in

\vskip.10in
\th{Theorem 3.11} Suppose that $m(dx)=dx$, $\sigma$ is bounded and that the conditions (A1)(A2')(A3)(A4), (H1)-(H5), (H6) hold. Then there exists a unique generalized solution of equation (3.1) and it satisfies the following estimates with constants $K_1$ and $K_2$ independent of $u, \phi,f$
$$\|u\|_T^2\leq K_1(\|\phi\|_2^2+\int_0^T\|f_t^0\|_2^2dt).$$
$$ \|u\|_\infty\leq K_2(\|\phi\|_\infty+\|f^0\|_\infty).$$

\proof Set for $\rho>0$ $f^\rho(t,x,y,z):=f(t,x,y,z)+\rho\sum_{l=1}^k\sum_{i=1}^d\sigma_l^i(x)\partial_i\theta(x)z_l(x),$
and consider
\begin{equation}(\partial_t+L_\rho)u+f^\rho(u,D_\sigma u)=0,\qquad u_T=\phi.\end{equation}
The associated weak equation has the form $\forall \varphi\in b\mathcal{C}^\rho_T$
\begin{equation}\int_0^T\mathcal{E}^\rho(u_t,\varphi_t)+(u_t,\partial_t\varphi_t)_\rho dt=\int_0^T(f_t^\rho(u_t,D_\sigma u_t),\varphi_t)_\rho dt+(u_T,\varphi_T)_\rho-(u_0,\varphi_0)_\rho.\end{equation}
As $f^\rho$ satisfies conditions (H1)-(H5), we have a generalized solution $u^\rho$ of (3.8).

Fix $\rho>0$ and take $f_n\in C_0^\infty(\mathbb{R}^d)$ such that $f_n(x)=1$ for $x\in B_n(0)$, $f_n(x)=0$ for $x\in B_{2n}^c(0)$, $\partial_{x_i}f_n(x)$ are uniformly bounded and $\partial_{x_i}f_n(x)\rightarrow0$ as $n\rightarrow\infty$.
If $\varphi\in b\mathcal{C}_T$, then $\varphi f_n\exp(\theta\rho)\in b\mathcal{C}_T^\rho$.
As
$$\int_0^T\mathcal{E}^\rho(u_t^\rho,\varphi_tf_n\exp(\theta\rho))+(u_t^\rho,\partial_t\varphi_tf_n) dt=\int_0^T(f_t^\rho(u_t^\rho,D_\sigma u_t^\rho),f_n\varphi_t) dt+(u_T^\rho,f_n\varphi_T)-(u_0^\rho,f_n\varphi_0),$$
by Lemma 3.10 we have
\begin{equation}\int_0^T\mathcal{E}(u_t^\rho,\varphi_tf_n)+(u_t^\rho,\partial_t\varphi_tf_n) dt=\int_0^T(f_t(u_t^\rho,D_\sigma u_t^\rho),f_n\varphi_t) dt+(u_T^\rho,f_n\varphi_T)-(u_0^\rho,f_n\varphi_0).\end{equation}
If $u\in \hat{F}_{\tilde{\rho}}$ satisfies (3.10) for fixed $\tilde{\rho}$ with test function $\varphi\in b\mathcal{C}_T$, then
$u$ satisfies (3.9) for  $\rho \geq \tilde{\rho}$, with test functions $\varphi $ where $\varphi\in b\mathcal{C}_T^\rho$.

Now fix $\rho_1>0$. Then there exists a solution $u^{\rho_1}$ of (3.8) associated to $\rho_1$. We  conclude $u^{\rho_1}$ satisfies the weak equation (3.9) for all $\rho>\rho_1$ with $\varphi\in b\mathcal{C}_T^\rho$. Then by Lemma 3.9 and the same arguments as in the uniqueness proof of Theorem 3.8 we have  $u^{\rho_1}=u^\rho$ for all $\rho>\rho_1$.

Finally, we deduce that a solution $u^{\tilde{\rho}}$ of (3.8) associated to $\tilde{\rho}$ is a solution of (3.8) for all $\rho>0$. Then by Theorem 3.8, we have
$\|u^{\tilde{\rho}}\|_{T,\rho}^2\leq K_1(\|\phi\|_{2,\rho}^2+\int_0^T\|f_t^0\|_{2,\rho}^2dt).$
Letting $\rho\rightarrow0$, by Fatou's Lemma, we obtain
$\aligned \|u^{\tilde{\rho}}\|_T^2&=\lim_{\rho\rightarrow0}\|u^{\tilde{\rho}}\|_{T,\rho}^2\leq \lim_{\rho\rightarrow0}K_1(\|\phi\|_{2,\rho}^2+\int_0^T\|f_t^0\|_{2,\rho}^2dt)
=K_1(\|\phi\|_{2}^2+\int_0^T\|f_t^0\|_{2}^2dt),\endaligned$
and
$ \|u^{\tilde{\rho}}\|_\infty\leq K_2(\|\phi\|_\infty+\|f^0\|_\infty).$
By (H6),  we have  $u^{\tilde{\rho}}\in L^2((0,T),F)$. For $u^{\tilde{\rho}}\in \hat{F}^\rho$ for $\rho>0$, we obtain
$\|u_{t+h_n}^{\tilde{\rho}}-u_t^{\tilde{\rho}}\|_{2,\rho}\rightarrow0.$
Then there exists a subsequence such that $u_{t+h_{n_k}}^{\tilde{\rho}}\rightarrow u_t^{\tilde{\rho}}$ for $m$-almost every $x$. Hence, $u_{t+h_{n_k}}^{\tilde{\rho}}\rightarrow u_t^{\tilde{\rho}}$ for $dx$-almost every $x$. Then by the same arguments as Lemma 3.9, we have
$$\aligned|\|u_t^{\tilde{\rho}}\|_{2,\rho}^2-\|u_{t+h}^{\tilde{\rho}}\|_{2,\rho}^2|\leq &2[|\int_t^{t+h}(u_s^{\tilde{\rho}},f_s^\rho)_\rho ds|+|\int_t^{t+h}\mathcal{E}^{a,\hat{b}}(u_s^{\tilde{\rho}})ds|+\alpha\int_t^{t+h}\|u_s^{\tilde{\rho}}\|^2_{2,\rho}ds]\\\leq &M\int_t^{t+h}\|f_s\|_{2,\rho} ds+M\int_t^{t+h}\mathcal{E}_{c_2+1}^{a,\hat{b}}(u_s^{\tilde{\rho}})ds
\endaligned$$
Letting $\rho\rightarrow0$, we get
$\aligned|\|u_t^{\tilde{\rho}}\|_{2}^2-\|u_{t+h}^{\tilde{\rho}}\|_{2}^2|\leq  M\int_t^{t+h}\|f_s\|_{2} ds+M\int_t^{t+h}\mathcal{E}_{c_2+1}^{a,\hat{b}}(u_s^{\tilde{\rho}})ds.
\endaligned$
Hence we have $u_{t+h_{n_k}}^{\tilde{\rho}}\rightarrow u_t^{\tilde{\rho}}$ in $L^2(\mathbb{R}^d,dx)$. Since this reasoning holds for every sequence $h_n\rightarrow0$, we have
$u^{\tilde{\rho}}\in \mathcal{C}([0,T],L^2)$, hence $u^{\tilde{\rho}}\in\hat{F}$.
By above arguments, we  deduce that
$$\int_0^T\mathcal{E}(u_t^{\tilde{\rho}},\varphi_tf_n)+(u_t^{\tilde{\rho}},\partial_t\varphi_tf_n) dt=\int_0^T(f_t(u_t^{\tilde{\rho}},D_\sigma u_t^{\tilde{\rho}}),f_n\varphi_t) dt+(u_T^{\tilde{\rho}},f_n\varphi_T)-(u_0^{\tilde{\rho}},f_n\varphi_0).$$
Letting $n\rightarrow\infty$, we conclude that
$\int_0^T\mathcal{E}(u_t^{\tilde{\rho}},\varphi_t)+(u_t^{\tilde{\rho}},\partial_t\varphi_t) dt=\int_0^T(f_t(u_t^{\tilde{\rho}},D_\sigma u_t^{\tilde{\rho}}),\varphi_t) dt+(u_T^{\tilde{\rho}},\varphi_T)-(u_0^{\tilde{\rho}},\varphi_0).$

Since $f_t(u_t^{\tilde{\rho}},D_\sigma u_t^{\tilde{\rho}})\in L^1([0,T];L^2)$, we can choose  $(f^n)_{n\in N}\subset C_0^\infty ([0,T]\times \mathbb{R}^d)$ such that $\int_0^T\|f_t^n-f_t(u_t^{\tilde{\rho}},D_\sigma u_t^{\tilde{\rho}})\|_2dt\rightarrow0$. Let $v^n$ be the generalized solution associated with $(f^n,\phi)$. Then $v^n$ is bounded. For $v_t:=P_{T-t}\phi(x)+\int_t^TP_{s-t}f(s,\cdot,u_s^{\tilde{\rho}},D_\sigma u_s^{\tilde{\rho}})(x)ds,$ we have $\|v^n-v\|_T\rightarrow0$. On the other hand, by Lemma 3.9
we have
$$\aligned\|u_t^{\tilde{\rho}}-v^n_t\|_2^2+2\int_t^T\mathcal{E}^{a,\hat{b}}(u_s^{\tilde{\rho}}-v^n_s)ds\leq &2\int_t^T(f_s(u_s^{\tilde{\rho}},D_\sigma u_s^{\tilde{\rho}})-f^n_s,u_s^{\tilde{\rho}}-v^n_s)ds+2\alpha\int _t^T\|u_s^{\tilde{\rho}}-v^n_s\|_2^2ds\\\leq&2M\int_t^T\|f_s(u_s^{\tilde{\rho}},D_\sigma u_s^{\tilde{\rho}})-f^n_s\|_2ds+2\alpha\int _t^T\|u_s^{\tilde{\rho}}-v^n_s\|_2^2ds.
\endaligned$$
By Gronwall's lemma we obtain $\|v^n-u^{\tilde{\rho}}\|_T\rightarrow0$, as $n\rightarrow\infty$. Therefore,
 we have $u_t^{\tilde{\rho}}=v_t$.
$\hfill\Box$
\vskip.10in
\section{Martingale representation for the processes}
\subsection{Representation under $P^x$}
In order to obtain the results for the probabilistic part, we need $\mathcal{E}$ to be a generalized Dirichlet form in the sense of Remark 2.1 (iii) with $c_2,\hat{c}\equiv0$ and $c\equiv0$. The Markov process $X=(\Omega,\mathcal{F}_\infty,\mathcal{F}_t,X_t,P^x)$ is properly associated in the resolvent sense with $\mathcal{E}$, i.e. $R_\alpha f:=E^x\int_0^\infty e^{-\alpha t}f(X_t)dt$ is an $\mathcal{E}$-quasi-continuous $m$-version of $G_\alpha f$, where $G_\alpha,\alpha>0$ is the resolvent  of $\mathcal{E}$ and $f\in \mathcal{B}_b(\mathbb{R}^d)\cap L^2(\mathbb{R}^d;m)$. The coform $\hat{\mathcal{E}}$ introduced in Section 2 is a generalized Dirichlet form with the associated resolvent $(\hat{G}_\alpha)_{\alpha>0}$ and there exists an $m$-tight special standard process property associated in the resolvent sense with $\hat{\mathcal{E}}$. From now on, we obtain all the results under the above assumption.

As mentioned in Remark 2.1 (viii), such a process can be constructed by quasi-regularity ([23, IV. 1. Definition 1.7]) and a structural condition ([23, IV. 2. D3] on the domain $\mathcal{F}$ of the generalized Dirichlet form.

We will now introduce the spaces which will be relevant for our further investigations. Define
$$\aligned \mathcal{M}:=\{&M|M \textrm{ is a finite additive functional, } E^z[M_t^2]<\infty, E^z[M_t]=0\\& \textrm{ for }\mathcal{E}-q.e. z\in E \textrm{ and all } t\geq0\}.\endaligned$$
$M\in \mathcal{M}$ is called a \emph{martingale additive functional}(MAF). Furthermore, define
$\dot{\mathcal{M}}:=\{M\in\mathcal{M}|e(M)<\infty\}.$
Here $e(M)=\frac{1}{2}\lim_{\alpha\rightarrow\infty}\alpha^2E^m[\int_0^\infty e^{-\alpha t}M_t^2dt]$. The elements of $\dot{\mathcal{M}}$ are called martingale additive functional's (MAF) of finite energy. Let $M\in \mathcal{M}$. There exists an $\mathcal{E}$-exceptional set $N$, such that $(M_t,\mathcal{F}_t,P_z)_{t\geq 0}$ is a square integrable martingale for all $z\in E\backslash N$. By [26, Theorem 2.10] $\dot{\mathcal{M}}$ is a real Hilbert space with inner product $e$.

We consider the following condition:
\vskip.10in

(A5) $X$ is a continuous conservative Hunt process in the state space $\mathbb{R}^d\cup\{\partial\}$.   $\hat{G}_\alpha$ is  strongly continuous on
 $\mathcal{V}$ and $\hat{\mathcal{E}}$ is quasi-regular.   $C_0^\infty(\mathbb{R}^d)\subset \mathcal{F}$ and  for $u\in \mathcal{F}$, there exists  a sequence $\{u_n\}\subset C_0^\infty(\mathbb{R}^d)$ such that $\mathcal{E}(u_n-u,u_n-u)\rightarrow0, n\rightarrow\infty$. $F_k:=\{x\in \mathbb{R}^d,|x|\leq k\}$  is an $\mathcal{E}$-nest.

\vskip.10in

\th{Remark 4.1} The last two conditions in (A5) are satisfied if  $C_0^\infty(\mathbb{R}^d)$ is dense in $\mathcal{F}$. If $\mathcal{E}$ satisfies the weak sector condition, (A5) can be easily verified. They are also satisfied in our next two examples and Example 5.6  which do not satisfy the weak sector condition.

\vskip.10in
\th{Example 4.2} Consider  $b=(b^i):\mathbb{R}^d\rightarrow \mathbb{R}^d$ be a Borel-measurable vector field. Let us define
$$Lu=\Delta u+\langle b, \nabla u \rangle, \qquad \forall u\in C_b^\infty(\mathbb{R}^d).$$
Assume that $\lim_{|x|\rightarrow\infty}\langle b(x),x\rangle=-\infty,$
and that there exist $C_1,C_2,m\in[0,\infty)$ such that
$|b(x)|\leq C_1+C_2|x|^m, x\in \mathbb{R}^d$

Then by [8, Theorem 5.3], there exists a probability measure $\mu$ on $\mathbb{R}^d$ such that
$\int_{\mathbb{R}^d}Lu d\mu=0, \forall u\in C_b^\infty(\mathbb{R}^d).$
And $b\in L^2(\mu).$ By [8, Theorem 3.1]
 we have $d\mu\ll dx$ and  the density admits a representation $\varphi^2$, where $\varphi\in H^{1,2}(\mathbb{R}^d,dx)$. The closure of
$\mathcal{E}^0(u,v)=\frac{1}{2}\int\langle \nabla u,\nabla v\rangle d\mu; u,v\in C_0^\infty(\mathbb{R}^d),$
on $L^2(\mathbb{R}^d,\mu)$ is a Dirichlet form.
Denote $b^0:=2\nabla\varphi/\varphi$ and $\beta:=b-b^0$. Then we have $\beta\in L^2(\mathbb{R}^d;\mathbb{R}^d,\mu)$.
Then by [22, Proposition 1.10 and Proposition 2.4]  $(L,C_0^\infty(\mathbb{R}^d))$ is $L^1$-unique. Then by the proof of [22, Proposition 2.4] for $u\in b\mathcal{F}$  there exists
a sequence $\{u_n\}\subset C_0^\infty(\mathbb{R}^d)$ such that $\mathcal{E}(u_n-u,u_n-u)\rightarrow0, n\rightarrow\infty.$

Consider the bilinear form
$$\mathcal{E}(u,v)=\frac{1}{2}\int\langle \nabla u,\nabla v\rangle d\mu-\int\langle \frac{1}{2}\beta,\nabla u\rangle vd\mu \textrm{ } u,v\in C_0^\infty(\mathbb{R}^d).$$
Then by the computation in [26, Section 4d] we have that  conditions (A1)-(A5) hold for the bilinear form $\mathcal{E}$.
\vskip.10in

\th{Example 4.3} Consider $d\geq 2$, $A=(a^{ij})$ a Borel-measurable mapping on $\mathbb{R}^d$ with values in the non-negative symmetric matrices on $\mathbb{R}^d$, and let $b=(b^i):\mathbb{R}^d\rightarrow \mathbb{R}^d$ be a Borel-measurable vector field. Consider the operator
$L_{A,b}\psi=\partial_i(a^{ij}\partial_j\psi)+b^i\partial_i\psi,  \forall \psi\in C_0^\infty(\mathbb{R}^d),$
 where we use the standard summation rule for repeated indices. By $H^{1,p}(\mathbb{R}^d,dx)$ we denote the standard Sobolev space of functions on $\mathbb{R}^d$ whose first order derivatives are in $L^p(\mathbb{R}^d,dx)$. Assume that for $p>d$

 (C1)$a^{ij}\in H^{p,1}_{\rm{loc}}(\mathbb{R}^d,dx), (a^{ij})$ is uniformly strictly elliptic in $\mathbb{R}^d$.

 (C2)$b^i\in L_{\rm{loc}}^p(\mathbb{R}^d,dx)$.

 Here by $H_{loc}^{1,p}(\mathbb{R}^d,dx)$ we denote the class of all functions $f$ on $\mathbb{R}^d$ such that $f\chi\in H^{1,p}(\mathbb{R}^d,dx)$ for all $\chi\in C_0^\infty(\mathbb{R}^d)$. And $L_{\rm{loc}}^p(\mathbb{R}^d,dx)$ denotes the class of all functions $f$ on $\mathbb{R}^d$ such that $f\chi\in L^p(\mathbb{R}^d)$ for all $\chi\in C_0^\infty(\mathbb{R}^d)$.
 Assume that there exists $V\in C^2(\mathbb{R}^d)$ ("Lyapunov function") such that
 $\lim_{|x|\rightarrow\infty} V(x)=+\infty,\lim_{|x|\rightarrow\infty}L_{A,b}V(x)=-\infty.$
Examples of $V$ can be found in \cite{BRS} and the reference therein.

Then by [9, Theorem 2.2] there exists a probability measure $\mu$ on $\mathbb{R}^d$ such that
$\int_{\mathbb{R}^d}L_{A,b}\psi d\mu=0, \forall \psi\in C_0^\infty(\mathbb{R}^d).$
Then by [9, Theorem 2.1] we have $d\mu\ll dx$ and that the density admits a representation $\varphi^2$, where $\varphi^2\in H_{loc}^{1,p}(\mathbb{R}^d,dx)$. The closure of
$\mathcal{E}^0(u,v)=\frac{1}{2}\int\langle \nabla u a,\nabla v\rangle d\mu; u,v\in C_0^\infty(\mathbb{R}^d),$
on $L^2(\mathbb{R}^d,\mu)$ is a Dirichlet form.

If in addition, there is a positive Borel function $\theta$ on $[0,\infty)$ such that $\lim_{t\rightarrow\infty}\theta(t)=+\infty$ and
$L_{A,b}V(x)\leq c_1-c_2\theta(|A^{-\frac{1}{2}}b|)|A^{-\frac{1}{2}}b|^2, $
outside some ball, then by [7, Theorem 2.6] $b\in L^2(\mathbb{R}^d;\mathbb{R}^d,\mu)$.
Set $b^0=(b_1^0,...,b_d^0)$, where $b_i^0:=2\sum_{j=1}^da^{ij}\partial_j\varphi/\varphi,i=1,...,d$. And $\beta:=b-b^0$.  By [7, Theorem 2.6] $\beta\in L^2(\mathbb{R}^d;\mathbb{R}^d,\mu)$.
Then  by [22, Proposition 1.10 and Proposition 2.4]  $(L,C_0^\infty(\mathbb{R}^d))$ is $L^1$-unique. Then by the proof of [22, Proposition 2.4] for $u\in b\mathcal{F}$  there exists
a sequence $\{u_n\}\subset C_0^\infty(\mathbb{R}^d)$ such that $\mathcal{E}(u_n-u,u_n-u)\rightarrow0, n\rightarrow\infty.$

Consider the bilinear form
$$\mathcal{E}(u,v)=\frac{1}{2}\int\langle \nabla u a,\nabla v\rangle d\mu-\int\langle \frac{1}{2}\beta,\nabla u\rangle vd\mu \textrm{ } u,v\in C_0^\infty(\mathbb{R}^d).$$
Then by the computation in [26, Section 4d] we have that  conditions (A1)-(A5) hold for the bilinear form $\mathcal{E}$.
\vskip.10in

For  an initial distribution $\mu\in \mathcal{P}(\mathbb{R}^d)$, here $\mathcal{P}(\mathbb{R}^d)$ denotes all the probabilities on $\mathbb{R}^d$,  we will prove the \emph{Fukushima reprensentation property} mentioned in \cite{QY} holds for $X$, i.e.
there is an algebra $K(\mathbb{R}^d)\subset\mathcal{B}_b(\mathbb{R}^d)$ which generates the Borel $\sigma$-algebra $\mathcal{B}(\mathbb{R}^d)$ and is invariant under $R_\alpha$ for $\alpha>0$, and there are finitely many continuous martingales $M^1,...,M^d$ over $(\Omega,\mathcal{F}^\mu,\mathcal{F}_t^\mu,P^\mu)$ such that for any potential $u=R_\alpha f$, where $\alpha>0$ and $f\in K(\mathbb{R}^d)$, the martingale part $M^{[u]}$ of the semimartingale $u(X_t)-u(X_0)$ has the martingale representation in terms of $(M^1,...,M^d)$, that is, there are predictable processes $F_1,...,F_d$ on $(\Omega,\mathcal{F}^\mu,\mathcal{F}^\mu_t)$ such that $M^{[u]}_t=\sum_{j=1}^d\int_0^tF_s^jdM_s^j, P^\mu-a.s..$

\vskip.10in
By [25, Theorem 4.5], if $\hat{G}_\alpha$ is sub-Markovian and strongly continuous on $\mathcal{V}$, the Fukushima decomposition holds for $u\in\mathcal{F}$.

Let us first calculate the energy measure related to $\langle M^{[u]}\rangle, u\in C_0^\infty(\mathbb{R}^d)$. By  [26, (23)], for bounded $g\in L^1(\mathbb{R}^d,m)$, we have
$$\aligned\int \hat{G}_\gamma gd\mu_{\langle M^{[u]}\rangle}=&\lim_{\alpha\rightarrow\infty}\alpha(U_{\langle M^{[u]}\rangle}^{\alpha+\gamma}1,\hat{G}_\gamma g)\\=&\lim_{\alpha\rightarrow\infty}\lim_{t\rightarrow\infty}E_{\hat{G}_\gamma g\cdot m}(\alpha e^{-(\gamma+\alpha)t}\langle M^{[u]}\rangle_t)+\lim_{\alpha\rightarrow\infty}E_{\hat{G}_\gamma g\cdot m}(\int_0^\infty\langle M^{[u]}\rangle_t\alpha(\gamma+\alpha)e^{-(\gamma+\alpha)t}dt)\\=&\lim_{\alpha\rightarrow\infty}\lim_{t\rightarrow\infty}\alpha\langle \mu_{\langle M^{[u]}\rangle},e^{-(\gamma+\alpha)t}\int_0^t\hat{P}_s \hat{G}_\gamma gds\rangle\\&+ \lim_{\alpha\rightarrow\infty}\alpha(\gamma+\alpha)(\int_0^\infty e^{-(\gamma+\alpha)t}E_{\hat{G}_\gamma g\cdot m}((u(X_t)-u(X_0)-N_t^{[u]})^2)dt)\\=&\lim_{\alpha\rightarrow\infty}\alpha(\gamma+\alpha)(\int_0^\infty e^{-(\gamma+\alpha)t}E_{\hat{G}_\gamma g\cdot m}((u(X_t)-u(X_0))^2)dt)
\\=&\lim_{\alpha\rightarrow\infty}2\alpha(u-\alpha G_\alpha u, u\hat{G}_\gamma g)-\alpha(u^2, \hat{G}_\gamma g-\alpha\hat{G}_\alpha \hat{G}_\gamma g)
\\=&2(-Lu,u \hat{G}_\gamma g)-(-Lu^2,\hat{G}_\gamma g)\\
=&2\mathcal{E}(u,u\hat{G}_\gamma g)-\mathcal{E}(u^2,\hat{G}_\gamma g)\\
=&2\mathcal{E}^a(u,u\hat{G}_\gamma g)-\mathcal{E}^a(u^2,\hat{G}_\gamma g)+2\int \langle b\sigma, D_\sigma u\rangle u\hat{G}_\gamma g m(dx)-\int \langle b\sigma, D_\sigma u^2\rangle \hat{G}_\gamma g m(dx)\\
=&2\mathcal{E}^a(u,u\hat{G}_\gamma g)-\mathcal{E}^a(u^2,\hat{G}_\gamma g)\\
=&2\int\langle D_\sigma u,D_\sigma (u\hat{G}_\gamma g)\rangle dm-\int\langle D_\sigma u^2,D_\sigma (\hat{G}_\gamma g)\rangle dm\\
=&2\int\langle D_\sigma u,D_\sigma u\rangle \hat{G}_\gamma gdm.\endaligned$$
Thus, by [26, Theorem 2.5] we obtain
$\mu_{\langle M^{[u]}\rangle}=2\langle D_\sigma u,D_\sigma u\rangle\cdot dm.$
So, for $u,v\in C_0^\infty(\mathbb{R}^d)$,
  under $P^x$ for quasi every point $x$,
 \begin{equation}\langle M^{[u]},M^{[v]}\rangle_t=2\int_0^t(\sum_{i,j=1}^d\sum_{l=1}^k\sigma_l^i\sigma_l^j\frac{\partial u}{\partial x^i}\frac{\partial v}{\partial x^j})(X_s)ds.\end{equation}
Then by (A5) and [25, Theorem 4.4], we deduce (4.1) for every $u,v\in\mathcal{F}$.

By (A5)  $F_k=\{x\in \mathbb{R}^d,|x|\leq k\}, k\in \mathbb{N}$  is an $\mathcal{E}$-nest. By [26, Theorem 3.6], for $u_i(x)=x_i$, we have the Fukushima decomposition for $A^{[u_i]}$, and let $M^{(i)}\in\dot{\mathcal{M}}_{loc,(F_k)_{k\in N}}$ be the associated local martingale additive functional. We  define the stochastic integral $f\cdot M^{(i)}$ for $f\in L^2(\mathbb{R}^d;\mu_{\langle M^{(i)}\rangle})$ as in [12, p243], and  for $L\in \dot{\mathcal{M}}$ we have
$\langle f\cdot M^{(i)},L\rangle=f\cdot\langle  M^{(i)}, L\rangle.$

\vskip.10in
  \th{Theorem 4.4} Assume (A5) holds.  Let $u\in C_0^1(\mathbb{R}^d)$, then
$$M^{[u]}=\sum_{i=1}^d u_{x_i}\cdot M^{(i)}.$$
\proof By [26, Theorem 3.6], we have
$$\aligned&\langle M^{[u]}-\sum_{i=1}^d u_{x_i}\cdot M^{(i)}\rangle_t=\sum_{i,j=1}^n\int_0^tu_{x_i}(X_s)u_{x_j}(X_s)d\langle M^{(i)},M^{(j)}\rangle_s\\&-2\sum_{i,j=1}^n\int_0^tu_{x_i}(X_s)u_{x_j}(X_s)d\langle M^{(i)},M^{(j)}\rangle_s+\sum_{i,j=1}^n\int_0^tu_{x_i}(X_s)u_{x_j}(X_s)d\langle M^{(i)},M^{(j)}\rangle_s=0.\endaligned$$
$\hfill\Box$
\vskip.10in

 Then by [26, Lemma 2.4, Lemma 1.18], we have
 \begin{equation}\langle M^{(i)},M^{(j)}\rangle_t=2\int_0^t\sum_{l=1}^k\sigma_l^i(X_s)\sigma_l^j(X_s)ds.\end{equation}
\vskip.10in

\th{Lemma 4.5} Assume (A5) holds. Let $\mathcal{C}_1$  be a uniformly dense subset of $C_0(\mathbb{R}^d)$. Here $C_0(\mathbb{R}^d)$ denotes the continuous function with compact support.
Then the family $\{f\cdot M^{[u]}:f\in \mathcal{C}_1,u\in C_0^\infty(\mathbb{R}^d)\}$ of stochastic integrals is dense in $(\dot{\mathcal{M}},e)$.

\proof Suppose that an MAF $M\in\dot{\mathcal{M}}$ is $e$-orthogonal to the above family, namely, $\int_Xfd\mu_{\langle M,M^{[u]}\rangle}=0, \forall f\in \mathcal{C}_1, u\in C_0^\infty(\mathbb{R}^d)$.
 This identity extends to all $u\in\mathcal{F}$ by [25, (13)] and (A5).  Hence, $$\langle M,M^{[u]}\rangle=0 \qquad \forall u\in\mathcal{F}.$$ In particular, this holds for $u=G_\alpha g, \alpha>0, \forall g\in C_0(\mathbb{R}^d)$.
By [12, Theorem A.3.20] we deduce that $M=0$.$\hfill\Box$
\vskip.10in

\th{Theorem 4.6} Assume (A5) holds. Then the space $\dot{\mathcal{M}}$ can be represented by stochastic integrals based on  $M^{(i)}=M^{[x_i]},1\leq i\leq d$:
\begin{equation}\dot{\mathcal{M}}=\{\sum_{i=1}^df_i\cdot M^{(i)}: \sum_{i,j=1}^d\sum_{l=1}^k\int_{\mathbb{R}^d}(f_if_j\sigma_i^l\sigma_j^l)(x)m(dx)<\infty\},\end{equation}
and
$$e(\sum_{i=1}^df_i\cdot M^{(i)})=\sum_{i,j=1}^d\sum_{l=1}^k\int_{\mathbb{R}^d}(f_if_j\sigma_i^l\sigma_j^l)(x)m(dx).$$
\proof The space of the right hand side of (4.3) is dense in $(\dot{\mathcal{M}},e)$, since it contains the set
$\{f\cdot M^{[u]}=\sum_{i=1}^d(fu_{x_i})\cdot M^{(i)}; f\in C_0(\mathbb{R}^d), u\in C_0^1(\mathbb{R}^d)\},$
which is dense in $(\dot{\mathcal{M}},e)$ by Lemma 4.5. Hence, it is enough to show that the right hand side of (4.3) is closed in $(\dot{\mathcal{M}},e)$.

Suppose that $\lim_{n\rightarrow\infty}e(M_n-M)=0$ for
$M_n=\sum_{i=1}^d f_i^{(n)}\cdot M^{(i)}, \sum_{i,j=1}^d\int_{\mathbb{R}^d}a_{ij}f_i^{(n)}f_j^{(n)}dm<\infty, M\in \dot{\mathcal{M}}.$
Set ${f}^n:=(f_1^n,...,f_d^n).$ Since
$e(M_n-M_m)=\int_{\mathbb{R}^d}|({f}^n-{f}^m)\sigma|^2dm,$
we deduce that $f^n\sigma$ converges in $L^2(\mathbb{R}^d;m)$ to some $h_i$ for each $i=1,...,d.$ Let
$h=(h_1,...,h_d)$ and $f=h\tau$. Set $M'=\sum_{i=1}^df_i\cdot M^{(i)}$, then
$\aligned e(M_n-M')=\int_D|({f}^n-{f})\sigma|^2dm
=\int_D|{f}^n\sigma-{h}|^2dm,\endaligned$
which converges to zero as $n\rightarrow\infty$. Therefore, we have $M=M'$ and $$e(M)=\sum_{i,j=1}^d\sum_{l=1}^k\int_{\mathbb{R}^d}(f_if_j\sigma_i^l\sigma_j^l)(x)m(dx)<\infty.$$
$\hfill\Box$
\vskip.10in

As a consequence, $X$ satisfies  Assumption (1) in \cite{QY}. Hence by [19, Theorem 3.1], we have the martingale representation theorem for $X$:
\vskip.10in

\th{Theorem 4.7} Assume (A5) holds. Let $\mu\in \mathcal{P}(\mathbb{R}^d)$ charging no set of zero capacity. Then for any square-integrable martingale $N=(N_t)_{t\geq0}$ on $(\Omega, \mathcal{F}^\mu,\mathcal{F}^\mu_t, P^\mu)$, there are unique predictable processes $(F_t^i)$ such that
$$N_t-N_0=\sum_{i=0}^d\int_0^tF_s^idM_s^{(i)}\qquad P^\mu-a.s..$$
\vskip.10in

Moreover, by  an analogous method to [19, Theorem 3.1] we have the martingale representation theorem for $X$ which is similar to \cite{BPS05}.
\vskip.10in

\th{Theorem 4.8} Assume (A5) holds. Then there exists some exceptional set $\mathcal{N}$ such that the following representation result holds. For  every bounded $\mathcal{F}_\infty$-measurable random variable $\xi$, there  exists an predictable process $(\phi_1,...,\phi_d):[0,\infty)\times \Omega\rightarrow \mathbb{R}^d$, such that for each probability measure $\nu$, supported by $\mathbb{R}^d\setminus\mathcal{N}$, one has
$$\xi=E^\nu(\xi|\mathcal{F}_0)+\sum_{i=0}^d\int_0^\infty\phi_s^idM_s^{(i)}\qquad P^\nu-a.s.,$$
and
$$E^\nu\int_0^\infty|\phi_s\sigma(X_s)|^2ds\leq \frac{1}{2}E^\nu\xi^2.$$
If another predictable process $\phi'=(\phi_1',...,\phi_d')$ satisfies the same relations under a certain measure $P^\nu$, then one has $\phi_t'\sigma(X_t)=\phi_t\sigma(X_t), dt\times dP^\nu-a.s.$.

\proof Suppose that $\mathcal{N}$ is some fixed exceptional set. By $\mathcal{K}$ we denote the class of bounded random variables for which the statement holds outside this set. By the same arguments as in the proof of  [2, Theorem 4.7]  we have that if $(\xi_n)\subset \mathcal{K}$ is a uniformly bounded increasing sequence and $\xi=\lim_{n\rightarrow\infty}\xi_n$ then $\xi\in\mathcal{K}$.

Let $K(\mathbb{R}^d)\subset \mathcal{B}_b(\mathbb{R}^d)$ be a countable set which is closed under multiplication, generates the Borel $\sigma$-algebra $\mathcal{B}(\mathbb{R}^d)$ and $R_\alpha (K(\mathbb{R}^d))\subset K(\mathbb{R}^d)$ for each $\alpha \in \mathbb{Q}^+$. Such $K(\mathbb{R}^d)$ can be constructed as follows. Choose $N_0\subset \mathcal{B}_b(\mathbb{R}^d)$ to be a countable set which generates the Borel $\sigma$-algebra $\mathcal{B}(\mathbb{R}^d)$. For $l\geq 1$  define $N_{l+1}=\{g_1...g_k, R_\alpha fg_1...g_k,f,g_i\in N_l, k\in \mathbb{N}\cup \{0\},\alpha\in Q^+\}$ and $K(\mathbb{R}^d):=\cup_{l=0}^\infty N_l$.

Let $\mathcal{C}_0$ be all $\xi=\xi_1\cdot\cdot\cdot\xi_n$, $n\in N$,
$\xi_j=\int_0^\infty e^{-\alpha_jt}f_j(X_t)dt$, where $\alpha_j\in \mathbb{Q}^+$, $f_j\in K(\mathbb{R}^d), j=1,...,n.$ Following the same arguments as the proof in [19, Lemma 2.2], we have that the completion of the $\sigma$-algebra generated by $\mathcal{C}_0$ is $\mathcal{F}_\infty$ . By the first part of our proof a monotone class argument reduces the proof to the representation of a random variable in $\mathcal{C}_0$.

Let $\xi\in\mathcal{C}_0$. Following the same arguments as   the proof of [19, Theorem 3.1], we have
$N_t=E^x(\xi|\mathcal{F}_t)=\sum_m Z_t^m,$ where the sum is a finite one, and for each $m$, $Z^m=Z_t$ has the following form
$Z_t=V_tu(X_t),$
(the superscript $m$ will be dropped if no confusion may arise), where $V_t=\prod_{i=1}^{k'}\int_0^te^{-\beta_i s}g_i(X_s)ds$ and $u(x)=R_{\beta_1+...+\beta_k}(h_1(R_{\beta_2+...+\beta_k}h_2...(R_{\beta_k}h_k)...)$ for $\beta_i\in \mathbb{Q}^+, g_i, h_i\in K(\mathbb{R}^d)$. We have $u\in K(\mathbb{R}^d)$. Hence, by the Fukushima decomposition and the Fukushima representation we obtain
\begin{equation}u(X_t)-u(X_0)=M_t^{[u]}+A_t^{[u]}=\sum_{j=1}^d\int_0^t G_s^jdM_s^{(j)}+A_t^{[u]}\qquad P^x-a.s..\end{equation} for some predictable processes $G^j$.
Then by the same arguments as the proof of  [19, Theorem 3.1], we deduce that
$N_t=\sum_{i=1}^d\int_0^t\sum_mV_s^m\cdot G_s^{m,i}dM_s^{(i)}, P^x-a.s..$
As (4.4) holds for every $x$ outside a set of zero capacity. Then we take the exceptional set $\mathcal{N}$ in the assertion to be the union of all these exceptional sets corresponding to $u\in K(\mathbb{R}^d)$.$\hfill\Box$
\vskip.10in
One may represent separately the positive and the negative parts and then we have the following corollary.
\vskip.10in
\th{Corollary 4.9} Assume (A5) holds.  Let $\mathcal{N}$ be the set obtained in the preceding theorem. Then for any $\mathcal{F}_\infty$-measurable nonnegative random variable $\xi\geq0$ there exists a predictable process $\phi=(\phi_1,...,\phi_d):[0,\infty)\times \Omega\rightarrow \mathbb{R}^d$ such that the following holds
$$\xi=E^x(\xi|\mathcal{F}_0)+\sum_{i=0}^d\int_0^\infty\phi_s^idM_s^{(i)}\qquad P^x-a.s.,$$
and
$$E^\nu\int_0^\infty|\phi_s\sigma(X_s)|^2ds\leq \frac{1}{2}E^\nu\xi^2,$$
for each point $x\in\mathcal{N}^c$ such that $E^x\xi<\infty$.

If another predictable process $\phi'=(\phi_1',...,\phi_d')$ satisfies the same relations under a certain measure $P^x$, then one has $\phi_t'\sigma(X_t)=\phi_t\sigma(X_t), dt\times dP^x-a.s.$
\vskip.10in
\subsection{Representation under $P^m$}
In the following, we use the notation $\int_0^t\psi(s,X_s).dM_s:=\sum_{i=1}^d\int_0^t\psi_i(s,X_s)dM_s^{(i)}$.

\vskip.10in
\th{Lemma 4.10} Assume (A1)-(A5) hold. If $u\in\mathcal{D}(L)$ and $\psi\in\tilde{\nabla}u$, then
$$u(X_t)-u(X_0)=\int_0^t\psi(X_s).dM_s+\int_0^tLu(X_s)ds\qquad P^m-a.s..$$
\proof  The assertion  follows by the Fukushima decomposition,  (4.2) and Theorem 4.6.$\hfill\Box$
\vskip.10in

The aim of the rest of this section is to extend this representation to time dependent functions $u(t,x)$.
\vskip.10in

\th{Lemma 4.11} Assume (A1)-(A5) hold. Let $u:[0,T]\times \mathbb{R}^d\rightarrow \mathbb{R}$ be such that

(i) $\forall s,u_s\in \mathcal{D}(L)$ and $s\rightarrow Lu_s$ is continuous in $L^2$.

(ii) $u\in C_1([0,T];L^2)$.

Then clearly $u\in\mathcal{C}_T$. Moreover, for any $\psi\in \tilde{\nabla }u$ and any $s,t>0$ such that $s+t<T$, the following relation holds $P^m$-a.s.
$$u(s+t,X_t)-u(s,X_0)=\int_0^t\psi(s+r,X_r).dM_r+\int_0^t(\partial_s+L)u(s+r,X_r)dr.$$
\proof We prove the above relation with $s=0$, the general case being similar. Let $0=t_0<t_1<...<t_p=t$ be a partition of the interval $[0,t]$ and write
$u(t,X_t)-u(0,X_0)=\sum_{n=0}^{p-1}(u(t_{n+1},X_{t_{n+1}})-u(t_n,X_{t_n})).$
Then, on account of the preceding lemma, each term of the sum is expressed  as
$$\aligned &u(t_{n+1},X_{t_{n+1}})-u(t_n,X_{t_n})\\=&u(t_{n+1},X_{t_{n+1}})-u(t_{n+1},X_{t_{n}})+u(t_{n+1},X_{t_n})-u(t_n,X_{t_n})\\
=&\int_{t_n}^{t_{n+1}}\psi^{n+1}(X_s).dM_s+\int_{t_n}^{t_{n+1}}Lu_{t_{n+1}}(X_s)ds+\int_{t_n}^{t_{n+1}}\partial_su_s(X_{t_n})ds,\endaligned$$
where $\psi^{n+1}=(\psi^{n+1}_1,...,\psi^{n+1}_d)\in\tilde{\nabla}u_{t_{n+1}}$ and the last integral is obtained by using the Leibnitz-Newton formula for the $L^2$-valued function $s\rightarrow u_s$. Below we estimate in $L^2$ the differences between each term in the last expression and the similar terms corresponding to the formula we have to prove. Here we use $mP_t\leq m$ i.e. $\int P_tfdm\leq \int fdm$ for $f\in\mathcal{B}^+$. This holds since $\hat{P}_t$ is sub-Markovian.
Then we have
$$\aligned &E^m(\int_{t_n}^{t_{n+1}}\psi^{n+1}(X_s).dM_s-\int_{t_n}^{t_{n+1}}\psi(s,X_s).dM_s)^2\\
=&E^m\int_{t_n}^{t_{n+1}}|(\psi^{n+1}(X_s)-\psi(s,X_s))\sigma(X_s)|^2ds
\leq\int_{t_n}^{t_{n+1}}\mathcal{E}^a(u_{t_{n+1}}-u_s)ds.\endaligned$$
Since $s\rightarrow Lu_s$ is continuous in $L^2$, it follows that $s\rightarrow u_s$ is continuous w.r.t. $\mathcal{E}^a_1$-norm. Hence the difference appearing in the last integral $\mathcal{E}^a(u_{t_{n+1}}-u_s)$ is uniformly small, provided the partition is fine enough. From this one deduces that
$\sum_{n=0}^{p-1}\int_{t_n}^{t_{n+1}}\psi^{n+1}(X_s).dM_s\rightarrow \int_0^t\psi(s+r,X_r).dM_r.$
The next difference is estimated by using Minkowski's inequality
$$\aligned &(E^m(\sum_{n=0}^{p-1}\int_{t_n}^{t_{n+1}}(Lu_{t_{n+1}}-Lu_s)(X_s)ds)^2)^{1/2}\\
\leq&\sum_{n=0}^{p-1}\int_{t_n}^{t_{n+1}}(E^m(Lu_{t_{n+1}}-Lu_s)^2(X_s))^{1/2}ds
\leq\sum_{n=0}^{p-1}\int_{t_n}^{t_{n+1}}\|Lu_{t_{n+1}}-Lu_s\|_2ds,\endaligned$$
so that it is similarly expressed as in integral of a uniformly small quantity.

For the last difference we write
$$\aligned&(E^m(\sum_{n=0}^{p-1}\int_{t_n}^{t_{n+1}}(\partial_su_s(X_{t_n})-\partial_su_s(X_s))ds)^2)^{1/2}
\leq \sum_{n=0}^{p-1}\int_{t_n}^{t_{n+1}}(E^m(\partial_su_s(X_{t_n})-\partial_su_s(X_s))^2)^{1/2}ds\\
=&\sum_{n=0}^{p-1}\int_{t_n}^{t_{n+1}}(E^m(\partial_su_s(X_{t_n})^2+P_{s-t_n}(\partial_su_s)^2(X_{t_n})-2\partial_su_s(X_{t_n})(P_{s-t_n}
\partial_su)(X_{t_n})))^{1/2}ds\\
=&\sum_{n=0}^{p-1}\int_{t_n}^{t_{n+1}}(E^m((\partial_su_s(X_{t_n})-(P_{s-t_n}\partial_su_s)(X_{t_n}))^2+(P_{s-t_n}(\partial_su_s)^2(X_{t_n})-
((P_{s-t_n}\partial_su_s)(X_{t_n}))^2)))^{1/2}ds\\
\leq&\sum_{n=0}^{p-1}(\int_{t_n}^{t_{n+1}}\int(\partial_su_s-P_{s-t_n}\partial_su_s)^2+P_{s-t_n}(\partial_su_s)^2-(P_{s-t_n}\partial_su_s)^2dm))^{1/2}ds.\endaligned$$
From the hypotheses it follows that this will tend also to zero if the partition is fine enough. Hence the assertions follow$\hfill\Box$
\vskip.10in
\th{Theorem 4.12} Assume (A1)-(A5) hold. Let $f\in L^1([0,T];L^2)$ and $\phi\in L^2(\mathbb{R}^d)$ and define
$$u_t:=P_{T-t}\phi+\int_t^TP_{s-t}f_sds.$$
Then for each $\psi\in\tilde{\nabla}u$ and for each $s\in[0,T]$, the following relation holds $P^m$-a.s.
$$u(s+t,X_t)-u(s,X_0)=\int_0^t\psi(s+r,X_r).dM_r-\int_0^tf(s+r,X_r)dr.$$
In particular, if $u$ is a generalized solution of PDE (3.1), the following BSDE holds $P^m$-a.s.
$$u(t,X_{t-s})=\phi(X_{T-s})+\int_t^Tf(r,X_{r-s},u(r,X_{r-s}),D_\sigma u(r,X_{r-s}))dr-\int_{t-s}^{T-s}\psi(s+r,X_r).dM_r.$$
\proof Assume first that $\phi$ and $f$ satisfy the conditions in Proposition 2.6 (ii). Then we have $u$ satisfies the conditions in  Lemma 4.11. Then by Lemma 4.11, the assertion follows. For the general case we choose $u^n$ associated $(f^n,\phi^n)$ as in Proposition 2.9. Then we have if $n\rightarrow\infty$, $\|u^n-u\|_T\rightarrow0$.
For $u^n$ we have
\begin{equation}u^n(s+t,X_t)-u^n(s,X_0)=\int_0^t\psi^n(s+r,X_r).dM_r-\int_0^tf^n(s+r,X_r)dr.\end{equation}
As
$$\aligned E^m|\int_0^t(\psi^n(s+r,X_r)-\psi^p(s+r,X_r).dM_r|^2&\leq E^m\int_0^t|(\psi^n(s+r,X_r)-\psi^p(s+r,X_r))\sigma(X_r)|^2dr
\\&\leq \int_0^t\mathcal{E}^a(u^n_{s+r}-u^p_{s+r})dr.\endaligned$$
Letting $n\rightarrow\infty$ in (4.5), we obtain the assertion.
$\hfill\Box$
\vskip.10in

\section{BSDE's and Generalized Solutions}
The set $\mathcal{N}$ obtained in Theorem 4.8 will be fixed throughout this Section. By Theorem 4.8 we can solve BSDE's under all measures $P^x$, $x\in\mathcal{N}^c$, at the same time. We will treat systems of $l$ equations, $l\in \mathbb{N}$, associated to $\mathbb{R}^l$-valued functions $f:[0,T]\times \Omega\times \mathbb{R}^l\times \mathbb{R}^l\otimes \mathbb{R}^k\mapsto \mathbb{R}^l$. These functions are assumed to depend on the past in general and it turns out that a good theory is developed assuming that they are predictable. This means that we consider the map $(s,\omega)\mapsto f(s,\omega,\cdot,\cdot)$ as a  predictable process with respect to the canonical filtration of our process $(\mathcal{F}_t)$.

\vskip.10in
\th{Lemma 5.1} Assume (A5) holds. Let $\xi$ be an $\mathcal{F}_T$-measurable random variable and $f:[0,T]\times \Omega\mapsto \mathbb{R}$ an $(\mathcal{F}_t)_{t\geq0}$-predictable process. Let $A$ be the set of all points $x\in \mathcal{N}^c$ for which 
$E^x(|\xi|+\int_0^T|f(s,\omega)|ds)^2<\infty$ holds.
Then there exists a pair $(Y_t,Z_t)_{0\leq t\leq T}$ of predictable processes $Y:[0,T)\times \Omega\mapsto \mathbb{R}, Z:[0,T)\times \Omega\mapsto \mathbb{R}^d$, such that under all measures $P^x$, $x\in A$, they have the following properties:

\no (i) $Y$ is continuous,

\no (ii) $Z$ satisfy the integrability condition
$\int_0^T|Z_t\sigma(X_t)|^2dt<\infty, \qquad P^x-a.s..$

\no(iii) The local martingale $\int_0^tZ_s.dM_s$, obtained by integrating $Z$ against the coordinate martingales,  is a uniformly integrable martingale,

\no(iv) they satisfy the equation
$Y_t=\xi+\int_t^Tf(s,\omega)ds-\int_t^TZ_s.dM_s, P^x-a.s., 0\leq t\leq T.$
If another pair $(Y_t',Z_t')$ of predictable processes satisfies the above conditions (i),(ii),(iii),(iv), under a certain measure $P^\nu$ with the initial distribution $\nu$ supported by $A$, then one has $Y.=Y.', P^\nu-a.s.$ and $Z_t\sigma(X_t)=Z_t'\sigma(X_t), dt\times P^\nu-a.s..$

\proof The representation of the positive and negative parts of the random variable $\xi+\int_0^Tf_sds$ give us the predictable process $Z$ such that
$\xi+\int_0^Tf_sds=E^{X_0}(\xi+\int_0^Tf_sds)+\int_0^TZ_s.dM_s.$
Then we get the process $Y$ by the formula
$Y_t=E^{X_0}(\xi+\int_0^Tf_sds)+\int_0^tZ_s.dM_s-\int_0^tf_sds.$
$\hfill\Box$
\vskip.10in

\th{Definition 5.2} Let $\xi$ be an $\mathbb{R}^l$-valued, $\mathcal{F}_T$-measurable, random variable and $f:[0,T]\times \Omega\times \mathbb{R}^l\times \mathbb{R}^l\otimes \mathbb{R}^k\mapsto \mathbb{R}^l$ a measurable $\mathbb{R}^l$-valued function such that $(s,\omega)\mapsto f(s,\omega,\cdot,\cdot)$ as a process is predictable. Let $p>1$ and $\nu$ be a probability measure supported by $\mathcal{N}^c$ such that $E^\nu|\xi|^p<\infty$. We say that a pair $(Y_t,Z_t)_{0\leq t\leq T}$ of predictable processes $Y:[0,T)\times \Omega\mapsto \mathbb{R}^l$, $Z:[0,T)\times\Omega\mapsto \mathbb{R}^l\otimes \mathbb{R}^d$ is a solution of the BSDE (5.1) in $L^p(P^\nu)$ with data $(\xi,f)$ provided that $Y$ is continuous under $P^\nu$ and satisfies both the integrability conditions
$\int_0^T|f(t,\cdot,Y_t,Z_t\sigma(X_t))|dt<\infty, P^\nu-a.s.,$
$E^\nu(\int_0^T|Z_t\sigma(X_t)|^2dt)^{p/2}<\infty,$
and the following equation, with $0\leq t\leq T$,
\begin{equation}Y_t=\xi+\int_t^Tf(s,\omega, Y_s,Z_s\sigma(X_s))ds-\int_t^TZ_s.dM_s,\qquad P^\nu-a.s..\end{equation}
\vskip.10in

Let $f:[0,T]\times \Omega\times \mathbb{R}^l\times \mathbb{R}^l\otimes \mathbb{R}^k\mapsto \mathbb{R}^l$ be a measurable $\mathbb{R}^l$-valued function such that $(s,\omega)\mapsto f(s,\omega,\cdot,\cdot)$ is predictable and satisfies the following conditions:

\no($\Omega1$) (\emph{Lipschitz condition in $z$}) There exists a  constant $C>0$ such that for all $t,\omega,y,z,z'$,
$|f(t,\omega,y,z)-f(t,\omega,y,z')|\leq C|z-z'|.$

\no($\Omega2$) (\emph{Monotonicity condition in $y$})
 There exists a function $\mu_t\in L^1([0,T],\mathbb{R})$ such that for all $\omega,y,y',z$,
$\langle y-y',f(t,\omega,y,z)-f(t,\omega,y',z)\rangle
\leq \mu_t|y-y'|^2,$
and   $\alpha_t:=\int_0^t\mu_sds<\infty.$

\no($\Omega3$) (\emph{Continuity condition in $y$})
For $t,\omega$ and $z$ fixed, the map
$y\mapsto f(t,\omega,y,z),$ is continuous.

We need the following notation
$f^0(t,\omega):=f(t,\omega,0,0), f'(t,\omega,y):=f(t,\omega,y,0)-f(t,\omega,0,0)$, $f^{',r}(t,\omega):=\sup_{|y|\leq r}|f'(t,\omega,y)|.$

Let $\xi$ be an $\mathbb{R}^l$-valued, $\mathcal{F}_T$-measurable, random variable and, for each $p>0$ denote by $A_p$ the set of all points $x\in\mathcal{N}^c$ for which the following integrability conditions hold,
\begin{equation}E^x\int_0^Tf_t^{',r}dt<\infty, \qquad \forall r\geq 0,\end{equation}
$$E^x(|\xi|^p+(\int_0^T|f^0(s,\omega)|ds)^p)<\infty.$$
Denote by $A_\infty$ the set of points $x\in\mathcal{N}^c$ for which (5.2) holds and with the property that $|\xi|,|f^0|\in L^\infty(P^x)$.

\vskip.10in
\th{Proposition 5.3} Under the conditions (A5), $(\Omega1),(\Omega2),(\Omega3)$, there exists a pair $(Y_t,Z_t)_{0\leq t\leq T}$ of predictable processes $Y:[0,T)\times \Omega\mapsto \mathbb{R}^l, Z:[0,T)\times\Omega\mapsto \mathbb{R}^l\otimes \mathbb{R}^d$ that forms a solution of the BSDE (5.1) in $L^p(P^x)$ with data $(\xi,f)$ for each point $x\in A_p$. Moreover, the following estimate holds with some constant $K$ that depends only on $C,\mu$ and $T$,
$$E^x(\sup_{t\in[0,T]}|Y_t|^p+(\int_0^T|Z_t\sigma(X_t)|^2dt)^{p/2})\leq KE^x(|\xi|^p+(\int_0^T|f^0(s,\omega)|ds)^p), \qquad x\in A_p.$$
If $x\in A_\infty$, then $\sup_{t\in[0,T]}|Y_t|\in L^\infty(P^x)$.

If $(Y_t',Z_t')$ is another solution in $L^p(P^x)$, for some point $x\in A_p$, then one has $Y_t=Y_t'$ and $Z_t\sigma(X_t)=Z_t'\sigma(X_t),dt\times P^x-a.s.$.

\vskip.10in

The proof is based on more or less standard methods. Therefore, we include it not here, but in the Appendix below.

We shall now look at the connection between the solutions of BSDE's introduced in this Section and PDE's studied in Section 3. In order to do this we have to consider BSDE's over time intervals like $[s,T]$, with $0\leq s\leq T$. Since the present approach is based on the theory of Markov processes, which is a time homogeneous theory, we have to discuss solutions over the interval $[s,T]$, while the process and the coordinate martingales are indexed by a parameter in the interval $[0,T-s]$.

Let us give a formal definition for the natural notion of solution over a time interval $[s,T]$. Let $\xi$ be an $\mathcal{F}_{T-s}$-measurable, $\mathbb{R}^l$-valued, random variable and $f:[s,T]\times \Omega\times \mathbb{R}^l\times \mathbb{R}^l\otimes \mathbb{R}^k\rightarrow \mathbb{R}^l$ an $\mathbb{R}^l$-valued, measurable map such that $(f(s+l,\omega,\cdot,\cdot))_{l\in[0,T-s]}$ is predictable with respect to $(\mathcal{F}_l)_{l\in[0,T-s]}$. Let $\nu$ be a probability measure supported by $\mathcal{N}^c$ such that $E^\nu|\xi|^p<\infty$. We say a pair $(Y_t,Z_t)_{s\leq t\leq T}$ of processes $Y:[s,T]\times\Omega\rightarrow \mathbb{R}^l, Z:[s,T]\times\Omega\rightarrow \mathbb{R}^l\otimes \mathbb{R}^d$ is a solution in $L^p(P^\nu)$ of the BSDE (5.3) over the interval $[s,T]$ with data $(\xi,f)$, provided that they have the property that reindexed as $(Y_{s+l},Z_{s+l})_{l\in[0,T-s]}$ these processes are $(\mathcal{F}_l)_{l\in[0,T-s]}$-predictable, $Y$ is continuous and together they satisfy the integrability conditions
$\int_s^T|f(t,\cdot,Y_t,Z_t\sigma(X_{t-s}))|dt<\infty, P^\nu-a.s.,$
$E^\nu(\int_s^T|Z_t\sigma(X_{t-s})|^2dt)^{p/2}<\infty.$
and the following equation under $P^\nu$,
\begin{equation}Y_t=\xi+\int_t^Tf(r,Y_r,Z_r\sigma(X_{r-s}))dr-\int_{t-s}^{T-s}Z_{s+l}.dM_l,\qquad s\leq t\leq T.\end{equation}
The next result gives a probabilistic interpretation of Theorem 3.8. Let us assume that $f:[0,T]\times \mathbb{R}^d\times \mathbb{R}^l\times \mathbb{R}^l\otimes \mathbb{R}^k\rightarrow \mathbb{R}^l$
is the measurable function appearing in the basic equation (3.1).
Let $\phi:\mathbb{R}^d\rightarrow \mathbb{R}^l$ be measurable and for each $p>1$, denote by $A_p$ the set of points $(s,x)\in[0,T)\times \mathcal{N}^c$ with the following properties
\begin{equation}E^x\int_s^Tf^{',r}(t,X_{t-s})dt<\infty, \qquad \forall r\geq 0.\end{equation}
$$E^x(|\phi|^p(X_{T-s})+(\int_s^T|f^0(t,X_{t-s})|ds)^p)<\infty.$$
Set $A=\cup_{p>1}A_p, A_{p,s}=\{x\in\mathcal{N}^c,(s,x)\in A_p\}$, and $A_s=\cup_{p>1}A_{p,s},s\in [0,T)$. By the same arguments as in [2, Theorem 5.4], we have the following results. In particular, we can reconstruct  solutions to PDE (3.1) using Proposition 5.3.
\vskip.10in

\th{Theorem 5.4} Assume that (A5) holds and $f$ satisfies conditions (H1),(H2),(H3). Then there exist nearly Borel measurable functions $(u,\psi), u:A\rightarrow \mathbb{R}^l,\psi:A\rightarrow \mathbb{R}^l\otimes \mathbb{R}^d$, such that, for each $s\in[0,T)$ and each $x\in A_{p,s}$, the pair $(u(t,X_{t-s}),\psi(t,X_{t-s}))_{s\leq t\leq T}$ solves the BSDE (5.3) in $L^p(P^x)$ with data $(\phi(X_{T-s}),f(t,X_{t-s},y,z))$ over the interval $[s,T]$.

In particular, the functions $u,\psi$ satisfy the following estimates, for $(s,x)\in A_p$,
$$E^x(\sup_{t\in[s,T]}|u(t,X_{t-s})|^p+(\int_s^T|\psi\sigma(t,X_{t-s})|^2dt)^{p/2})\leq KE^x(|\phi(X_{T-s})|^p+(\int_s^T|f^0(t,X_{t-s})|dt)^p).$$
Moreover,  suppose (A1)-(A4) hold, and the conditions in Theorem 3.11 hold when $m(dx)=dx$. If $f$ and $\phi$ satisfy the conditions (H4) and (H5)
then the complement of $A_{2.s}$ is $m$-negligible (i.e. $m(A_{2,s}^c)=0$) for each $s\in[0,T)$, the class of $u1_{A_2}$ is an element of $\hat{F}^l$ which is a generalized solution of PDE (3.1), $\psi\sigma$ represents a version of $D_\sigma u$ and the following relations hold for each $(s,x)\in A$ and $1\leq i\leq l,$
\begin{equation}u^i(s,x)=E^x(\phi^i(X_{T-s}))+\int_s^TE^xf^i(t,X_{t-s},u(t,X_{t-s}),D_\sigma u(t,X_{t-s}))dt.\end{equation}

\vskip.10in
\th{Remark 5.5} In the above theorem, we need the analytic results, i.e. the existence of a generalized solution of nonlinear equation (3.1), to obtain the above results. In the following example, we drop the conditions (A1)-(A4), in particular, we don't need $|b\sigma|\in L^2(\mathbb{R}^d;m)$  and use the results that the existence of the solution of BSDE (5.3) to obtain the existence of a generalized solution of nonlinear equation (3.1), which is not covered by our analytic results in Section 2.
\vskip.10in

\th{Example 5.6} Consider $d\geq 2$, $A=(a^{ij})$ a Borel-measurable mapping on $\mathbb{R}^d$ with values in the non-negative symmetric matrices on $\mathbb{R}^d$, and let $b=(b^i):\mathbb{R}^d\rightarrow \mathbb{R}^d$ be a Borel-measurable vector field. Consider the operator
$$L_{A,b}\psi=a^{ij}\partial_i\partial_j\psi+b^i\partial_i\psi, \qquad \forall \psi\in C_0^\infty(\mathbb{R}^d),$$
 where we use the standard summation rule for repeated indices. By $H^{1,p}(\mathbb{R}^d,dx)$ we denote the standard Sobolev space of functions on $\mathbb{R}^d$ whose first order derivatives are in $L^p(\mathbb{R}^d,dx)$. Assume that for $p>d$

 (C1)$a^{ij}=a^{ji}\in H^{p,1}_{\rm{loc}}(\mathbb{R}^d,dx), 1\leq i,j\leq d$.

 (C2)$b^i\in L_{\rm{loc}}^p(\mathbb{R}^d,dx)$.

 (C3) for all $K$ relatively compact in $\mathbb{R}^d$ there exist $\nu_K>0$ such that
 $$\nu_K^{-1}|h|^2\leq \langle ha,h\rangle\leq \nu_K|h|^2 \textrm{ for all }h\in \mathbb{R}^d, x\in K.$$

 Here by $H_{loc}^{1,p}(\mathbb{R}^d,dx)$ we denote the class of all functions $f$ on $\mathbb{R}^d$ such that $f\chi\in H^{1,p}(\mathbb{R}^d,dx)$ for all $\chi\in C_0^\infty(\mathbb{R}^d)$. And $L_{\rm{loc}}^p(\mathbb{R}^d,dx)$ denotes the class of all functions $f$ on $\mathbb{R}^d$ such that $f\chi\in L^p(\mathbb{R}^d)$ for all $\chi\in C_0^\infty(\mathbb{R}^d)$.
 Assume that there exists $V\in C^2(\mathbb{R}^d)$ ("Lyapunov function") such that
 $\lim_{|x|\rightarrow\infty} V(x)=+\infty,\lim_{|x|\rightarrow\infty}L_{A,b}V(x)=-\infty.$
Examples of $V$ can be found in \cite{BRS} and the reference therein.

Then by [9, Theorem 2.2] there exists a probability measure $\mu$ on $\mathbb{R}^d$ such that $\int_{\mathbb{R}^d}L_{A,b}\psi d\mu=0, \forall \psi\in C_0^\infty(\mathbb{R}^d).$

Then by [9, Theorem 2.1] we have $d\mu\ll dx$ and that the density admits a representation $\varphi^2$, where $\varphi^2\in H_{loc}^{1,p}(\mathbb{R}^d,dx)$. The closure of
$\mathcal{E}^0(u,v)=\frac{1}{2}\int\langle \nabla u a,\nabla v\rangle d\mu;  u,v\in C_0^\infty(\mathbb{R}^d),$
on $L^2(\mathbb{R}^d,\mu)$ is a Dirichlet form.

Set $b^0=(b_1^0,...,b_d^0)$, where $b_i^0:=\sum_{j=1}^d(\partial_j a_{ij}+2a^{ij}\partial_j\varphi/\varphi),i=1,...,d$. And $\beta:=b-b^0$. Then, $\beta\in L^2_{\rm{loc}}(\mathbb{R}^d;\mathbb{R}^d,\mu)$.
Then  by [22, Proposition 1.10 and Proposition 2.4]  $(L,C_0^\infty(\mathbb{R}^d))$ is $L^1$-unique. Then by the proof of [22, Proposition 2.4] for $u\in b\mathcal{F}$  there exists
a sequence $\{u_n\}\subset C_0^\infty(\mathbb{R}^d)$ such that $\mathcal{E}(u_n-u,u_n-u)\rightarrow0, n\rightarrow\infty.$

Consider the bilinear form
$$\mathcal{E}(u,v)=\frac{1}{2}\int\langle \nabla u a,\nabla v\rangle d\mu-\int\langle \frac{1}{2}\beta,\nabla u\rangle vd\mu \textrm{ } u,v\in C_0^\infty(\mathbb{R}^d).$$
Then by the computation in [26, Section 4d] we have that  conditions (A5) hold for the bilinear form $\mathcal{E}$. Then we can use the first part of Theorem 5.4 to obtain the following results.
\vskip.10in

\th{Theorem 5.7} Consider the bilinear form obtained in Example 5.6. If $f$ satisfies conditions (H1),(H2),(H3). Then there exist nearly Borel measurable functions $(u,\psi), u:A\rightarrow \mathbb{R}^l,\psi:A\rightarrow \mathbb{R}^l\otimes \mathbb{R}^d$, such that, for each $s\in[0,T)$ and each $x\in A_{p,s}$, the pair $(u(t,X_{t-s}),\psi(t,X_{t-s}))_{s\leq t\leq T}$ solves the BSDE (5.3) in $L^p(P^x)$ with data $(\phi(X_{T-s}),f(t,X_{t-s},y,z))$ over the interval $[s,T]$.

In particular, the functions $u,\psi$ satisfy the following estimates, for $(s,x)\in A_p$,
$$E^x(\sup_{t\in[s,T]}|u(t,X_{t-s})|^p+(\int_s^T|\psi\sigma(t,X_{t-s})|^2dt)^{p/2})\leq KE^x(|\phi(X_{T-s})|^p+(\int_s^T|f^0(t,X_{t-s})|dt)^p).$$ Moreover, suppose $f$ and $\phi$ satisfy the conditions (H4) and (H5)
then the complement of $A_{2.s}$ is $\mu$-negligible (i.e. $\mu(A_{2,s}^c)=0$) for each $s\in[0,T)$, the class of $u1_{A_2}$ is an element of $\hat{F}^l$ which is a generalized solution of (3.1), $\psi\sigma$ represents a version of $D_\sigma u$ and the following relations hold for each $(s,x)\in A$ and $1\leq i\leq l,$
$$u^i(s,x)=E^x(\phi^i(X_{T-s}))+\int_s^TE^xf^i(t,X_{t-s},u(t,X_{t-s}),D_\sigma u(t,X_{t-s}))dt.$$

\proof By [22, Lemma 3.1], we have for $u\in D(L_{A,b})$, $u\in D(\mathcal{E}^0)$ and $\mathcal{E}^0(u,u)\leq -\int Luud\mu$. By this, the first part of proof in Proposition 2.9 hold in this case i.e. the mild solution is equivalent to the generalized solution. Then the results in Theorem 4.12 hold. By the same arguments as in the proof of Theorem 5.4 and using $P_t\mu=\mu$, the assertion follows.$\hfill\Box$
\vskip.10in

\section{Further Examples}

The following two examples  discuss the case where PDE satisfies some boundary conditions.

\th{Example 6.1} Let $D\subset \mathbb{R}^d$ be a bounded domain satisfying the cone condition. We choose $m(dx)=1_D(x)dx$. If $\mathcal{E}$ is a sectorial Dirichlet form, it is associated to a reflecting diffusion $X$ in the state space $\overline{D}$.   Then by Theorem 3.8 there exists a solution to the non-linear parabolic equation
$$\aligned (\partial_t+L)u+f(t,x,u,D_\sigma u)&=0,\qquad  0\leq t\leq T,\\
u_T(x)&=\phi(x), \qquad x\in \mathbb{R}^d,\\ \frac{\partial u(t,\cdot)}{\partial \nu}|_{\partial D}&=0, \qquad t>0,\endaligned$$
where $\frac{\partial}{\partial \nu}$ denotes the normal derivative. Then Theorem 5.5 provides a probabilistic interpretation for this equation.

\vskip.10in

\th{Example 6.2} Let $D\subset \mathbb{R}^d$ be a bounded domain satisfying the cone condition. We choose $m(dx)=1_D(x)m(dx)$ and replace $C_0^\infty(\mathbb{R}^d)$ by $C_0^\infty(D)$. Then the results in Theorem 3.8 apply and there exists a solution $u_1\in F=H^1_0(D)$ to the following non-linear parabolic equation:
$$\aligned (\partial_t+L)u+f(t,x,u,D_\sigma u)&=0,\qquad  0\leq t\leq T,\\
u_T(x)&=\phi(x), \qquad x\in \mathbb{R}^d.\endaligned$$

Assume $\mathcal{E}$ satisfies the weak sector condition. Let $X^0$ denote the diffusion associated with $\mathcal{E}^R$, where $\mathcal{E}^R$ denotes the Dirichlet form which has the same form as $\mathcal{E}$ with the reference measure  $m(dx)$ replaced by $dx$.  Then define
 $$X_t:=\left\{\begin{array}{ll}X_t^0,&\ \ \ \ \textrm{ if }t<\tau,\\ \Delta &\ \ \ \ \textrm{ otherwise, } \end{array}\right.$$
where $\tau=\inf\{t\geq 0, X_t^0\in D^c\cup \Delta\}$. Assume (A5) holds for $X^0$. We use Theorem 5.5 for $X^0$ with the data
$(\phi(X_{T-s}^0)1_{\{T-s<\tau\}},1_{[0,\tau+s]}(r)f(r,X_{r-s}^0,Y_r,Z_r\sigma(X_{r-s}^0))$. Then there exist nearly Borel measurable functions
$(u,\psi), u:A\rightarrow \mathbb{R}^l,\psi:A\rightarrow \mathbb{R}^l\otimes \mathbb{R}^d$, such that, for each $s\in[0,T)$ and each $x\in A_{p,s}$, the pair $(u(t,X_{t-s}^0),\psi(t,X_{t-s}^0))_{s\leq t\leq T}$ solves the BSDE
$$Y_t=\phi(X_{T-s}^0)1_{\{T-s<\tau\}}+\int_{t\wedge(\tau+s)}^{T\wedge (\tau+s)}f(r,X_{r-s}^0,Y_r,Z_r\sigma(X_{r-s}^0))dr-\int_{t-s}^{T-s}Z_{s+l}.dM_l,\qquad s\leq t\leq T.$$
Then by [17, Proposition 2.6] we have $Y_t=0,Z_t=0 \textrm{ when } t\in [\tau+s, T],$
and the pair $(u(t,X_{t-s}),\psi(t,X_{t-s}))_{s\leq t\leq T}$ solves the BSDE
$$Y_t=\phi(X_{T-s})+\int_{t}^{T}f(r,X_{r-s},Y_r,Z_r\sigma(X_{r-s}))dr-\int_{t-s}^{T-s}Z_{s+l}.dM_l,\qquad s\leq t\leq T.$$

In particular, the functions $u,\psi$ satisfy the following estimates, for $(s,x)\in A_p$,
$$E^x(\sup_{t\in[s,T]}|u(t,X_{t-s})|^p+(\int_s^T|\psi\sigma(t,X_{t-s})|^2dt)^{p/2})\leq KE^x(|\phi(X_{T-s})|^p+(\int_s^T|f^0(t,X_{t-s})|dt)^p).$$
 The class of $u1_{A_2}$ is an element in $\hat{F}^l$ which is an $m$-version of $u_1$, $\psi\sigma$ represents a version of $D_\sigma u$ and the following relations hold for each $(s,x)\in A$ and $1\leq i\leq l,$
\begin{equation}u^i(s,x)=E^x(\phi^i(X_{T-s}))+\int_s^TE^xf^i(t,X_{t-s},u(t,X_{t-s}),\psi(t,X_{t-s})\sigma(X_{t-s}))dt.\end{equation}

\th{Acknowledgement.} I would like to thank Professor M. R\"{o}ckner for valuable discussions and for suggesting me to associate BSDE's to generalized Dirichlet forms which was one motivation for this paper . I would also like to thank Professor Ma Zhiming and Zhu Xiangchan for their helpful discussions.

\vskip.10in
\th{Appendix A. Proof of Theorem 3.8}

[Uniqueness]

Let $u_1$ and $u_2$ be two solutions of equation (3.1). By using (2.7) for the difference $u_1-u_2$ we get
$$\aligned&\|u_{1,t}-u_{2,t}\|_2^2+2\int_t^T\mathcal{E}^{a,\hat{b}}(u_{1,s}-u_{2,s})ds\\\leq& 2\int_t^T(f(s,\cdot,u_{1,s},D_\sigma u_{1,s})-f(s,\cdot,u_{2,s},D_\sigma u_{2,s}),u_{1,s}-u_{2,s})ds+2\alpha\int _t^T\|u_{1,s}-u_{2,s}\|_2^2ds\\
\leq &2\int_t^T C(|D_\sigma u_{1,s}-D_\sigma u_{2,s}|,|u_{1,s}-u_{2,s}|)ds+2\alpha\int _t^T\|u_{1,s}-u_{2,s}\|_2^2ds
\\ \leq & (\frac{C^2}{c_1}+c_2+2\alpha) \int _t^T\|u_{1,s}-u_{2,s}\|_2^2ds+\int_t^T\mathcal{E}^{a,\hat{b}}(u_{1,s}-u_{2,s})ds.\endaligned$$
By Gronwall's lemma it follows that
 $u_1=u_2$.

[Existence] The existence will be proved in four steps.

\textbf{Step 1}: Suppose there exists  $r\in \mathbb{R}$ such that $r\geq 1+K(\|\phi\|_\infty+\|f^0\|_\infty+\|f^{',1}\|_\infty),$
where $K$ is the constant appearing in Lemma 3.7 (3.5), and  $f$ is uniformly bounded on the set $A_r=[0,T]\times \mathbb{R}^d\times \{|y|\leq r\}\times \mathbb{R}^l\otimes \mathbb{R}^k.$
Define $M:=\sup\{|f(t,x,y,z)|:(t,x,y,z)\in A_r\}<\infty.$
Next we regularize $f$ with respect to the variable $y$ by convolution
$f_n(t,x,y,z)=n^l\int_{R^l}f(t,x,y',z)\varphi(n(y-y'))dy',$
where $\varphi$ is a smooth nonnegative function with support contained in the ball $\{|y|\leq 1\}$ such that $\int \varphi =1$.
Then $f=\lim_{n\rightarrow\infty}f_n$ and for each $n$, $\partial_{y_i}f_n$ are uniformly bounded on $A_{r-1}$.  Set $h_n(t,x,y,z):=f_n(t,x,\frac{r-1}{|y|\vee(r-1)}y,z).$
Then each $h_n$ satisfies the Lipschitz condition with respect to both $y$ and $z$. Thus by Proposition 3.4 each $h_n$ determines  a solution $u_n\in \hat{F}^l$  of (3.1) with data $(\phi,h_n)$.
By the same arguments as in [21, Theorem 4.19], we have that $h_n$ satisfies conditions (H1) and (H2') with the same constants ($C>0$ and $\mu=0$).
As $m$ is a finite measure and $f^{',1}\in L^\infty([0,T]\times \mathbb{R}^d)$, we have $f^{',1}\in L^2([0,T];L^2)$.
 Since
$$\aligned |h_n(t,x,0,0) \leq &n^l\int_{R^l}|f(t,x,y')-f^0(t,x)+f^0(t,x)||\varphi(n(-y'))|dy'\\
\leq &|f^0(t,x)|+f^{',1}(t,x),\endaligned$$
one deduces from Lemma 3.7 that $\|u_n\|_\infty\leq r-1$ and $\|u_n\|_T\leq K_T$.
Since $h_n=f_n$ on $A_{r-1}$, it follows that $u_n$ satisfies (3.1) with data $(\phi,f_n)$.

Now for $b>0$, set
$d_{n,b}(t,x):=\sup_{|y|\leq r-1,|z|\leq b}|f(t,x,y,z)-f_n(t,x,y,z)|.$
Obviously one has $|d_{n,b}|\leq 2M$. Moreover, on account of the $y$-continuity and of the uniform $z$-continuity, one sees that for fixed $t,x,b$, the family of functions
$\{f(t,x,\cdot,z)||z|\leq b\},$ is equicontinuous and then compact in $\mathcal{C}(\{|y|\leq r-1\})$. Since the convolution operators approach the identity uniformly on such a compact set, we get
$\lim_{n\rightarrow\infty}d_{n,b}(t,x)=0,$
which implies $\lim_{n\rightarrow\infty}d_{n,b}(t,x)=0$ in $L^2(dt\times m)$ because of our assumption that $m(\mathbb{R}^d)<\infty$.
Moreover, for $u\in \hat{F}^l,|u|\leq r-1$
$$\aligned|f(u,D_\sigma u)-f_n(u,D_\sigma u)|\leq &1_{\{|D_\sigma u|\leq b\}}d_{n,b}+2M1_{\{|D_\sigma u|>b\}}\\
\leq &d_{n,b}+\frac{2M}{b}|D_\sigma u|.\endaligned$$
Next we will show that $(u_n)_{n\in N}$ is a $\|\cdot\|_T$-Cauchy sequence. By (2.7)  for the difference $u_l-u_n$, we have
$$\aligned&\|u_{l,t}-u_{n,t}\|_2^2+2\int_t^T\mathcal{E}^{a,\hat{b}}(u_{l,s}-u_{n,s})ds\\\leq& 2\int_t^T(f_l(s,\cdot,u_{l,s},D_\sigma u_{l,s})-f_n(s,\cdot,u_{n,s},D_\sigma u_{n,s}),u_{l,s}-u_{n,s})ds+2\alpha\int _t^T\|u_{l,s}-u_{n,s}\|_2^2ds
\\\leq
&2\int_t^T(|f_l(s,\cdot,u_{l,s},D_\sigma u_{l,s})-f(s,\cdot,u_{l,s},D_\sigma u_{l,s})|,|u_{l,s}-u_{n,s}|)ds\\&+
2\int_t^T(|f_n(s,\cdot,u_{n,s},D_\sigma u_{n,s})-f(s,\cdot,u_{n,s},D_\sigma u_{n,s})|,|u_{l,s}-u_{n,s}|)ds\\&+
 2\int_t^T C(|D_\sigma u_{l,s}-D_\sigma u_{n,s}|,|u_{l,s}-u_{n,s}|)ds+2\alpha\int _t^T\|u_{l,s}-u_{n,s}\|_2^2ds\endaligned$$
$$\aligned
\leq &2\int_t^T(d_{l,b}(s,\cdot)+d_{n,b}(s,\cdot),|u_{l,s}-u_{n,s}|)ds+2\int_t^T\frac{2M}{b}(|D_\sigma u_{l,s}|+|D_\sigma u_{n,s}|,|u_{l,s}-u_{n,s}|)ds\\&+
 2\int_t^T C(|D_\sigma u_{l,s}-D_\sigma u_{n,s}|,|u_{l,s}-u_{n,s}|)ds +2\alpha\int _t^T\|u_{l,s}-u_{n,s}\|_2^2ds\\
 \leq & \int_t^T\|d_{l,b}(s,\cdot)\|_2^2ds+\int_t^T\|d_{n,b}(s,\cdot)\|_2^2ds+\frac{1}{b^2}\int_t^T(\|D_\sigma u_{l,s}\|_2^2+\|D_\sigma u_{n,s}\|_2^2)ds\\&+(1+4M^2+\frac{C^2}{c_1}+2\alpha+c_2)\int_t^T\|u_{l,s}-u_{n,s}\|_2^2ds+\int_t^T\mathcal{E}^{a,\hat{b}}(u_{l,s}-u_{n,s})ds.\endaligned$$
Since  $\|u_n\|_T\leq K_T$, we obtain $\int_0^T\|D_\sigma u_{l,s}\|_2^2ds\leq \frac{K_T}{c_1},$
where the $K_T$ is independent of $l$ and $b$. Thus, for $b,l,n$ large enough,  for  arbitrary $\varepsilon>0$ we get
$$\|u_{l,t}-u_{n,t}\|_2^2+\int_t^T\mathcal{E}^{a,\hat{b}}(u_{l,s}-u_{n,s})ds\leq \varepsilon+\tilde{K}\int_t^T\|u_{l,t}-u_{n,t}\|_2^2ds,$$
where $\tilde{K}$ depends on $C,M,\mu,\alpha$. It is easy to see that Gronwall's lemma implies that $(u_n)_{n\in N}$ is a Cauchy-sequence in $\hat{F}$. Define $u:=\lim_{n\rightarrow\infty} u_n$ and take a subsequence $(n_k)_{k\in N}$ such that $u_{n_k}\rightarrow u$ a.e. We have $f(\cdot,\cdot,u_{n_k},D_\sigma u)\rightarrow f(\cdot,\cdot,u, D_\sigma u) \textrm{ in } L^2(dt\times m).$
Since $\|u_{n_k}-u\|_T\rightarrow0$, we obtain
$\|D_\sigma u-D_\sigma u_{n_k}\|_{L^2(dt\times m)}\rightarrow0.$
Then by (H1), it follows that
$$\aligned &\lim_{k\rightarrow\infty}\|f(\cdot,\cdot,u_{n_k},D_\sigma u)-f(\cdot,\cdot,u_{n_k},D_\sigma u_{n_k})\|_{L^2(dt\times m)}\\ \leq &\lim_{k\rightarrow\infty} C\|D_\sigma u-D_\sigma u_{n_k}\|_{L^2(dt\times m)}\\=&0.\endaligned$$
We also have $$\aligned &\|f(\cdot,\cdot,u_{n_k},D_\sigma u_{n_k})-f_{n_k}(\cdot,\cdot,u_{n_k},D_\sigma u_{n_k})\|_{L^2(dt\times m)}\\ \leq &||d_{n_k,b}\|_{L^2(dt\times m)}+\frac{2M}{b}\|D_\sigma u_{n_k}\|_{L^2(dt\times m)}.\endaligned$$
Letting $k\rightarrow\infty$ and then $b\rightarrow\infty$ the above equality converges to zero.
Finally, we conclude
$$\aligned&\lim_{k\rightarrow\infty}\|f_{n_k}(u_{n_k},D_\sigma u_{n_k})-f(u,D_\sigma u)\|_{L^2(dt\times m)}\\
\leq &\lim_{k\rightarrow\infty}\|f_{n_k}(u_{n_k},D_\sigma u_{n_k})-f(u_{n_k},D_\sigma u_{n_k})\|_{L^2(dt\times m)}\\&+\lim_{k\rightarrow\infty}\|f(u_{n_k},D_\sigma u_{n_k})-f(u_{n_k},D_\sigma u)\|_{L^2(dt\times m)}\\&+\lim_{k\rightarrow\infty}\|f(u_{n_k},D_\sigma u)-f(u,D_\sigma u)\|_{L^2(dt\times m)}\\=&0.\endaligned$$
By passing to the limit in the mild equation associated to $u_{n_k}$ with data $(\phi,f_{n_k})$, it follows that $u$ is the solution associated to $(\phi,f)$.

\textbf{Step 2}: In this Step we will prove the assertion under the assumption that there exists some constant $r$ such that $f^{',r}$ is uniformly bounded and $r\geq 1+K(\|\phi\|_\infty+\|f^0\|_\infty+\|f^{',1}\|_\infty),$
where $K$ is the constant appearing in Lemma 3.7 (3.5). Define
$f_n(t,x,y,z):=f(t,x,y,\frac{n}{|z|\vee n}z).$
$f_n\leq Cn+\|f^{',r}\|_\infty+\|f^0\|_\infty$ on $A_r$. Each of the functions $f_n$ satisfies the same conditions as $f$ and by  Step 1, there exists a solution $u_n$ associated to the data $(\phi,f_n)$. One has $\|u_n\|_\infty\leq r-1$, $\|u_n\|_T\leq K_T$. Conditions (H1) and (H2) yield
$$\aligned&|(f_l(u_l,D_\sigma u_l)-f_n(u_n,D_\sigma u_n),u_l-u_n)|\\
\leq &C(|D_\sigma u_l-D_\sigma u_n|,|u_l-u_n|)+|(f_l(u_n,D_\sigma u_n)-f_n(u_n,D_\sigma u_n),u_l-u_n)|.\endaligned$$
Since $f_n(t,x,y,z)1_{|z|\leq n}=f(t,x,y,z)1_{|z|\leq n},$ and for $n\leq l$,
$|f_l-f_n|1_{|z|\geq n}\leq 2C|z|1_{|z|\geq n},$
we have $$|(f_l(u_n,D_\sigma u_n)-f_n(u_n,D_\sigma u_n),u_l-u_n)|\leq |(2C|D_\sigma u_n|1_{\{|D_\sigma u_n|\geq n\}},|u_l-u_n|)|.$$
Then,
$$\aligned&\|u_{l,t}-u_{n,t}\|_2^2+2\int_t^T\mathcal{E}^{a,\hat{b}}(u_{l,s}-u_{n,s})ds\\\leq& 2\int_t^T(f_l(u_{l,s},D_\sigma u_{l,s})-f_n(u_{n,s},D_\sigma u_{n,s}),u_{l,s}-u_{n,s})ds+2\alpha\int _t^T\|u_{l,s}-u_{n,s}\|_2^2ds\\
\leq& 2\int_t^TC(|D_\sigma u_l-D_\sigma u_n|,|u_l-u_n|)ds+2\int_t^T|(2C|D_\sigma u_n|1_{\{|D_\sigma u_n|\geq n\}},|u_l-u_n|)|ds\\&+2\alpha\int _t^T\|u_{l,s}-u_{n,s}\|_2^2ds\\\leq &(\frac{C^2}{c_1}+2\alpha+c_2)\int_t^T\|u_l-u_n\|_2^2ds+\int_t^T\mathcal{E}^{a,\hat{b}}(u_l-u_n)ds\\&+8C(r-1)(\int_t^T\|1_{\{|D_\sigma u_n|\geq n\}}\|_2^2ds)^{\frac{1}{2}}(\int_t^T\|D_\sigma u_n\|_2^2ds)^{\frac{1}{2}}.\endaligned$$
As $\|u_n\|_T^2\leq K_T$, we have $\int_0^T\|D_\sigma u_n\|_2^2ds\leq \frac{K_T}{c_1}$. Hence,
$$n^2\int_t^T\|1_{\{|D_\sigma u_n|\geq n\}}\|_2^2ds\leq \int_t^T\|D_\sigma u_n1_{\{|D_\sigma u_n|\geq n\}}\|_2^2ds\leq \frac{K_T}{c_1}.$$
Therefore,  for $n$ big enough
$$\|u_{l,t}-u_{n,t}\|_2^2+\int_t^T\mathcal{E}^{a,\hat{b}}(u_{l,s}-u_{n,s})ds\\\leq (\frac{C^2}{c_1}+2\alpha+c_2)\int_t^T\|u_l-u_n\|_2^2ds+\varepsilon.$$
By Gronwalls' lemma it follows that $(u_n)_{n\in N}$ is a Cauchy sequence in $\hat{F}^l$. Hence, $u:=\lim_{n\rightarrow\infty}u_n$ is well defined. We can find a subsequence such that $(u_{n_k},D_\sigma u_{n_k})\rightarrow (u,D_\sigma u)$ a.e.
and conclude $$|f_{n_k}(u_{n_k},D_\sigma u_{n_k})-f(u,D_\sigma u)|\leq C|\frac{n_k}{|D_\sigma u_{n_k}|\vee n_k}D_\sigma u_{n_k}-D_\sigma u|+|f(u_{n_k},D_\sigma u)-f(u,D_\sigma u)|\rightarrow0.$$
Since
$$\aligned &|f_{n_k}(u_{n_k},D_\sigma u_{n_k})-f(u,D_\sigma u)|\\\leq &|f(u,0)-f(u,D_\sigma u)|+|f_{n_k}(u_{n_k},D_\sigma u_{n_k})-f_{n_k}(u_{n_k},0)|+|f_{n_k}(u_{n_k},0)-f^0|+|f^0-f(u,0)|\\\leq &C(|D_\sigma u|+|D_\sigma u_{n_k}|)+2f^{',r},\endaligned$$
we have $$f_{n_k}(u_{n_k},D_\sigma u_{n_k})\rightarrow f(u,D_\sigma u) \textrm{ in } L^1([0,T],L^2).$$ We conclude $u$ is a solution of (3.1) associated to the data $(\phi,f)$.

\textbf{Step 3}: Now we only suppose that $f^{',1}$ is bounded. Hence, we can choose a constant $r$ such that
 $r\geq 1+K(\|\phi\|_\infty+\|f^0\|_\infty+\|f^{',1}\|_\infty),$
where $K$ is the constant appearing in Lemma 3.7 (3.5). Let us define
$f_n:=\frac{n}{f^{',r}\vee n}(f-f^0)+f^0.$
Easily we see that the $f_n$ have the same properties as $f$. Since $f_n(t,x,y,z)=f(t,x,y,z)$ for $f^{',r}\leq n$, we have
$\lim_{n\rightarrow\infty}f_n=f.$
We introduce the following notation:
$f_n^{',r}(t,x):=\sup_{|y|\leq r}|f_n'(t,x,y)|, \textrm{ and } f_n'(t,x,y):=f_n(t,x,y,0)-f^0(t,x).$
By the same arguments as in [21, Theorem 4.19] we have $|f_n^{',r}|\leq n\wedge |f^{',r}|.$
Hence, by  Step 2 we obtain that there exists a solution $u_n$ associated to the data $(\phi,f^n)$ such that $\|u_n\|_\infty\leq r-1,\|u_n\|_T\leq M$,  where $M$ is a constant. For $n\leq l$, we have
$$|f_l-f_n|\leq (C|z|+|f'|)|\frac{l}{f^{',r}\vee l}-\frac{n}{f^{',r}\vee n}|\leq (C|z|+|f'|)1_{\{f^{',r}>n\}}.$$
Hence  $$\int_t^T|(f_l(u_n,D_\sigma u_n)-f_n(u_n,D_\sigma u_n),u_l-u_n)|ds\leq 2(r-1)\int_t^T\int_{\{f^{',r}>n\}}(C|D_\sigma u_n|+f^{',r})dmds.$$
We start as in the preceding steps:
$$\aligned&\|u_{l,t}-u_{n,t}\|_2^2+2\int_t^T\mathcal{E}^{a,\hat{b}}(u_{l,s}-u_{n,s})ds\\\leq& 2\int_t^T(f_l(u_{l,s},D_\sigma u_{l,s})-f_n(u_{n,s},D_\sigma u_{n,s}),u_{l,s}-u_{n,s})ds+2\alpha\int _t^T\|u_{l,s}-u_{n,s}\|_2^2ds\\
\leq &(\frac{C^2}{c_1}+2\alpha+c_2)\int_t^T\|u_l-u_n\|_2^2ds+\int_t^T\mathcal{E}^{a,\hat{b}}(u_l-u_n)ds+4(r-1)\int_t^T\int_{\{f^{',r}>n\}}(C|D_\sigma u_n|+f^{',r})dmds.\endaligned$$
As $\lim_{n\rightarrow\infty}\int_t^T\int_{\{f^{',r}>n\}}f^{',r}dmds=0,$
and
$$\int_t^T\int_{\{f^{',r}>n\}}|D_\sigma u_n|dmdt\leq \|1_{\{f^{',r}>n\}}\|_{L^2(dt\times m)}\|D_\sigma u_n\|_{L^2(dt\times m)}\rightarrow0,$$
we have as above that $(u_n)_{n\in N}$ is a Cauchy sequence in $\hat{F}^l$. Hence, $u:=\lim_{n\rightarrow\infty}u_n$ exists in $\hat{F}^l$. We can find a subsequence such that $(u_{n_k},D_\sigma u_{n_k})\rightarrow (u,D_\sigma u)$ a.e. and we have that
$$\aligned &|f_{n_k}(u_{n_k},D_\sigma u_{n_k})-f(u,D_\sigma u)|\\\leq&
1_{\{f^{',r}\leq n_k\}}|f(u,D_\sigma u)-f(u_{n_k},D_\sigma u_{n_k})|+1_{\{f^{',r}>n_k\}}[|f(u,D_\sigma u)-f^0|+|f(u,D_\sigma u)-f(u_{n_k},D_\sigma u_{n_k})|]\\ \leq &|f(u,D_\sigma u)-f(u_{n_k},D_\sigma u_{n_k})|+1_{\{f^{',r}>n_k\}}|f(u,D_\sigma u)-f^0|\\ \leq &|f(u_{n_k},D_\sigma u)-f(u_{n_k},D_\sigma u_{n_k})|+|f(u_{n_k},D_\sigma u)-f(u,D_\sigma u)|+1_{\{f^{',r}>n_k\}}|f(u,D_\sigma u)-f^0|.\endaligned$$
As in the above proof we have
 $f_{n_k}(u_{n_k},D_\sigma u_{n_k})\rightarrow f(u,D_\sigma u),$ in $L^1([0,T],L^2)$. We conclude $u$ is a solution of (3.1) associated to the data $(\phi,f)$.

\textbf{Step 4}: Now we prove the theorem without additional conditions. Define
$f_n:=\frac{n}{f^{',1}\vee n}(f-f^0)+f^0.$
Since $f_n(t,x,y,z)=f(t,x,y,z)$ for $f^{',1}\leq n$, we have
$\lim_{n\rightarrow\infty}f_n=f.$
Introduce the following notation:
$f_n^{',1}(t,x):=\sup_{|y|\leq 1}|f_n'(t,x,y)| \textrm{ and } f_n'(t,x,y):=f_n(t,x,y,0)-f^0(t,x).$
As in Step 3 we have $|f_n^{',1}|\leq n\wedge|f^{',1}|.$
Since $f_n^{',1}$ is uniformly bounded, we can apply Step 3. Then we get a solution $u_n$ for the data $(\phi,f_n)$. The convergence of $u_n$ can be shown analogously to Step 3.$\hfill\Box$

\vskip.10in

\th{Appendix B. Proof of Proposition 5.3}

Let $M^p_x(\mathbb{R}^l)$ denote the set of (equivalence classes of )predictable processes $\{\phi_t\}_{t\in[0,T]}$ with values in $\mathbb{R}^l$ such that
$\|\phi\|_{M^p_x}:=(E^x[(\int_0^T|\phi_r|^2dr)^{p/2}])^{1/p}<\infty.$

$M^p_{\sigma,x}(\mathbb{R}^l\otimes \mathbb{R}^d)$ denotes the set of (equivalence classes of )predictable processes $\{\phi_t\}_{t\in[0,T]}$ with values in $\mathbb{R}^l\otimes \mathbb{R}^d$ such that
$\|\phi\|_{M^p_{\sigma,x}}:=(E^x[(\int_0^T|\phi_r\sigma(X_r)|^2dr)^{p/2}])^{1/p}<\infty.$

Fix $x\in A_p$.

We note that $(Y,Z)$ solves the BSDE (5.1) with data $(\xi,f)$ iff
$(\bar{Y}_t,\bar{Z}_t):=(e^{\alpha_t}Y_t,e^{\alpha_t}Z_t),$
solve the BSDE (5.1) with data $(e^{\alpha_T}\xi,f')$, where
$f'(t,y,z):=e^{\alpha_t}f(t,e^{-\alpha_t}y, e^{-\alpha_t}z)-\mu_t y.$
Therefore, we may replace $(\Omega2)$ by
$$\langle y-y',f(t,\omega,y,z)-f(t,\omega,y',z)\rangle
\leq 0, \textrm{ for all } t,x,y,y',z.$$

\textbf{Step 1} Assume that $f$ is Lipschitz continuous with respect to both $y$ and $z$.
Define a mapping $\Phi$ from $B^2_x:=M^2_x(\mathbb{R}^l)\times M^2_{\sigma,x}(\mathbb{R}^l\otimes \mathbb{R}^d)$ into itself as follows. Given $(U,V)\in B^2_x$, we can set $\Phi(U,V):=(Y,Z)$, where $(Y,Z)$ is the solution of the BSDE (5.1)  associated with data $(\xi,f(U,V\sigma(X)))$ given by Lemma 5.1. Then by It\^{o}'s formula and BDG inequality, we obtain
$E^x[\sup_{t\in [0,T]}|Y_t|^2]<\infty.$
Let $(U,V),(U',V')\in B^2_x$, $(Y,Z)=\Phi(U,V)$, $(Y',Z')=\Phi(U',V')$, $(\bar{U},\bar{V})=(U-U',V-V')$, $(\bar{Y},\bar{Z})=(Y-Y',Z-Z')$. It follows from Ito's formula that for each $\gamma\in R$,
$$\aligned &e^{\gamma t}E^x|\bar{Y}_t|^2+E^x\int_t^Te^{\gamma s}(\gamma |\bar{Y}_s|^2+|\bar{Z}_s\sigma(X_s)|^2)ds\\
\leq &2KE^x\int_t^Te^{\gamma s}|\bar{Y}_s|(|\bar{U}_s|+|\bar{V}_s\sigma(X_s)|)ds\\
\leq &4K^2E^x\int_t^Te^{\gamma s}|\bar{Y}_s|^2+\frac{1}{2}E^x\int_t^Te^{\gamma s}(|\bar{U}_s|^2+|\bar{V}_s\sigma(X_s)|^2)ds,\endaligned$$
where $K$ is the Lipschitz constant of $f$. We choose $\gamma=1+4K^2$. Then
$$E^x\int_0^Te^{\gamma s}( |\bar{Y}_s|^2+|\bar{Z}_s\sigma(X_s)|^2)ds\leq \frac{1}{2}E^x\int_0^Te^{\gamma s}(|\bar{U}_s|^2+|\bar{V}_s\sigma(X_s)|^2)ds,$$
from which it follows that $\Phi$ is a strict contraction on $B^2_x$ equipped with the norm:
$$|||(Y,Z)|||_\gamma^x=(E^x\int_0^Te^{\gamma t}( |Y_t|^2+|Z_t\sigma(X_t)|^2)dt)^{1/2}.$$

Define a sequence $(Y^n,Z^n)$ by $(Y^{n+1},Z^{n+1}):=\Phi(Y^n,Z^n)$. We have for  $\gamma=1+4K^2$
$$E^x\int_0^Te^{\gamma s}( |Y_s^n-Y_s^{n+1}|^2+|(Z_s^n-Z_s^{n+1})\sigma(X_s)|^2)ds\leq\frac{1}{2^n}E^x\int_0^Te^{\gamma s}( |Y_s^0-Y_s^{1}|^2+|(Z_s^0-Z_s^{1})\sigma(X_s)|^2)ds.$$
Then  we have
the a.s. pointwise convergence of $(Y^n_s,Z^n_s\sigma(X_s))$  under each measure $P^x$, $x\in A^2$. Denote the limit by $(Y_s,Z_s\sigma(X_s))$. Then this is the fixed point of $\Phi$ under the norm $|||(Y,Z)|||_\gamma^x$. So we have $(Y_s,Z_s)$ is the solution of BSDE (5.1).

\textbf{Step 2} We assume $f,\xi$ are bounded.

We need the following proposition.

\vskip.10in
\th{Proposition B.1} Assume  condition (A5). Given $V\in \cap_x M^2_{\sigma,x}(\mathbb{R}^l\otimes \mathbb{R}^d)$,  there exists a unique pair of predictable processes $(Y_t,Z_t)\in M^2_x\times M^2_{\sigma,x}(\mathbb{R}^l\otimes \mathbb{R}^d), \forall x\in \mathcal{N}^c$ satisfying under all $P^x$, $ x\in \mathcal{N}^c$
$$Y_t=\xi+\int_t^Tf(s,Y_s,V_s)ds-\int_t^TZ_sdM_s,\qquad 0\leq t\leq T.$$
\vskip.10in
Using Proposition B.1, we can construct a mapping $\Phi$ from $B^2_x$ into itself as follows. For any $(U,V)\in B^2_x$, $(Y,Z)=\Phi(U,V)$ is the solution of the BSDE$$Y_t=\xi+\int_t^Tf(s,Y_s,V_s)ds-\int_t^TZ_sdM_s,\qquad 0\leq t\leq T.$$
 Then as in \textbf{Step 1}, we have
$$\aligned &e^{\gamma t}E^x|\bar{Y}_t|^2+E^x\int_t^Te^{\gamma s}(\gamma |\bar{Y}_s|^2+|\bar{Z}_s\sigma(X_s)|^2)ds\\
= &2E^x\int_t^Te^{\gamma s}\langle\bar{Y}_s,f(Y_s,V_s\sigma(X_s))-f(Y_s',V_s'\sigma(X_s))\rangle ds\\
\leq &2KE^x\int_t^Te^{\gamma s}|\bar{Y}_s|\times|\bar{V}_s\sigma(X_s)|ds\\
\leq & E^x\int_t^Te^{\gamma s}(2K^2|\bar{Y}_s|^2+\frac{1}{2}|\bar{V}_s\sigma(X_s)|^2)ds.\endaligned$$
Then by the same argument as in  \textbf{Step 1}, we obtain the assertion of Proposition 5.3 if $f,\xi$ are bounded.
\vskip.10in

\vspace{1mm}\noindent{\it Proof of Proposition B.1}\quad We write $f(s,y)$ for $f(s,y,V_s)$.

By $C$ we denote the constant satifying $|\xi|^2+\sup_t|f(t,0)|^2\leq C$ a.s..
Define $f^n(t,y):=(\rho_n*f(t,\cdot))(y),$
where $\rho_n:\mathbb{R}^l\mapsto \mathbb{R}^+$ is a sequence of smooth functions with compact support satisfying $\int\rho_n(z)dz=1$, which approximate the Dirac measure at $0$. Then each $f^n$ is locally Lipschitz in $y$, uniformly with respect to $s$ and $\omega$.

Define for each $m\in N$, $f^{n,m}(t,y):=f^n(t,\frac{\inf(m,|y|)}{|y|}y).$ Then $f^{n,m}$ is globally Lipschitz and bounded, uniformly w.r.t. $(t,\omega)$. As in \textbf{Step 1}, we have a unique pair $(Y_t^{n,m}, Z_t^{n,m})\in M_x^2\times M^2_{\sigma,x}(\mathbb{R}^l\otimes \mathbb{R}^d)$ such that
$$Y_t^{n,m}=\xi+\int_t^Tf^{n,m}(s,Y_s^{n,m})ds-\int_t^TZ_s^{n,m}dM_s,\qquad 0\leq t\leq T.$$
By It\^{o}'s Formula we have $|Y_t^{n,m}|^2\leq e^TC, 0\leq t\leq T.$
Consequently, for $m^2>e^TC$, $(Y_t^{n,m}, Z_t^{n,m})$ does not depend on $m$. Therefor, we denote it by $(Y_t^{n}, Z_t^{n})$.
Then by the same arguments as [10, Proposition 3.2] we have
$$E^x(\sup_{0\leq t\leq T}|Y_t^k-Y_t^l|^2)+E^x(\int_0^T|(Z_t^k-Z_t^l)\sigma(X_t)|^2dt)\leq KE^x[\int_0^T|f^k(t,Y_t^k)-f^l(t,Y_t^k)|^2dt].$$
We have for fixed $\omega$,$\sup_{k>l}\int_0^T|f^k(t,Y_t^k)-f^l(t,Y_t^k)|^2dt\rightarrow0, l\rightarrow\infty.$
Then we have
$$\sup_{k>l}E^x\int_0^T|f^k(t,Y_t^k)-f^l(t,Y_t^k)|^2dt\leq E^x\sup_{k>l}\int_0^T|f^k(t,Y_t^k)-f^l(t,Y_t^k)|^2dt\rightarrow0,\qquad l\rightarrow\infty.$$
and we can obtain a sequence of representable variables that converges rapidly enough under all measures $P^x, x\in\mathcal{N}^c$. For each $l=0,1,...$ set
$$n_l(x)=\inf\{n>n_{l-1}(x);\sup_{k\geq n} E^x[\int_0^T|f^k(t,Y_t^k)-f^n(t,Y_t^k)|^2dt]<\frac{1}{2^l}\},$$
$$\bar{Y}^l=Y^{n_l(X_0)},\bar{Z}^l=Z^{n_l(X_0)}.$$
With this sequence one may pass to the limit and define $Z_s'=\limsup_{l\rightarrow\infty}\bar{Z}^l_s\sigma(X_s)$ and $Z_s=Z_s'\tau(X_s)$. Then we obtain the claimed results.$\hfill\Box$
\vskip.10in

So far we have proved the assertion when $\xi, f$ are bounded. Then by the same arguments as in [10, Theorem 4.2], one proves the general case.$\hfill\Box$
\vskip.10in

\th{Appendix C. Basic Relations for the Linear Equation}
In this section we assume that (A1)-(A4) hold.

\vskip.10in

\th{Lemma C.1} If $u$ is a bounded generalized solution of equation (2.5), then $u^+$ satisfies the following relation for $0\leq t_1<t_2\leq T$
$$\|u_{t_1}^+\|_2^2\leq 2\int_{t_1}^{t_2}(f_s,u_s^+)ds+\|u_{t_2}^+\|_2^2.$$

\proof Choose  the approximation sequence $u^n$ for $u$ as in the existence proof of Proposition 2.9. Denote its related data by $f^n,\phi^n$ .

Suppose that  the following holds
$$\|(u_{t_1}^n)^+\|_2^2\leq 2\int_{t_1}^{t_2}(f_s^n,(u_s^n)^+)ds+\|(u_{t_2}^n)^+\|_2^2,\eqno(C.1)$$
where $0\leq t_1\leq t_2\leq T$. Since
$\|u^n\|_2$ are uniformly bounded, we have $\lim_{n\rightarrow\infty}\int_{t_1}^{t_2}(f_s^n,(u_s^n)^+)ds=\int_{t_1}^{t_2}(f_s,u_s^+)ds.$

By passing $n$ to the limit in equation (C.1) the assertion follows.

Therefore, the problem is reduced to the case where $u$ belongs to $b\mathcal{C}_T$ ; in the
remainder we assume $u\in b\mathcal{C}_T$ .
(2.9), written with $u^+\in bW^{1,2}([0,T];L^2)\cap L^2([0,T];F)$ as test function, takes the form
$$\aligned& \int_{t_1}^{t_2}(u_t,\partial_t(u_t^+))dt+\int_{t_1}^{t_2}\mathcal{E}^{a,\hat{b}}(u_t,u_t^+)dt+\int_{t_1}^{t_2}\int\langle b\sigma, D_\sigma u_t
\rangle u_t^+dmdt
\\=&\int_{t_1}^{t_2}(f_t,u_t^+)dt+(u_{t_2},u^+_{t_2}))-(u_{t_1},u^+_{t_1})),\endaligned\eqno(C.2)$$
By [3, Theorem 1.19], we  obtain
$\int_{t_1}^{t_2}(u_t,\partial_t(u_t^+))dt=\frac{1}{2}(\|u_{t_2}^+\|_2^2-\|u_{t_1}^+\|_2^2).$
Then
$$\aligned&\|u_{t_1}^+\|_2^2 +2\int_{t_1}^{t_2}\mathcal{E}^{a,\hat{b}}(u_t,u_t^+)dt+2\int_{t_1}^{t_2}\int\langle b\sigma, D_\sigma u_t
\rangle u_t^+dmdt
=2\int_{t_1}^{t_2}(f_t,u_t^+)dt+\|u_{t_2}^+\|_2^2.\endaligned\eqno(C.3)$$
Next we prove for $u\in bF$
$$\mathcal{E}(u,u^+)\geq0.\eqno(C.4)$$

We have the above relation for $u\in \mathcal{D}(L)$. For $u\in bF$, by (A4) we can choose a uniformly bounded sequence $\{u_n\}\subset \mathcal{D}(L)$ such that  $\mathcal{E}^{a,\hat{b}}_{c_2+1}(u_n-u)\rightarrow0$. Then we have
$$\aligned&|\int\langle b\sigma, D_\sigma u
\rangle u^+dm-\int\langle b\sigma, D_\sigma u_n
\rangle u_n^+dm|
\rightarrow&0. \endaligned$$
Because $\mathcal{E}^{a,\hat{b}}(u^+)\leq \mathcal{E}^{a,\hat{b}}(u)$, $\sup_n\mathcal{E}^{a,\hat{b}}(u^+_n)\leq \sup_n\mathcal{E}^{a,\hat{b}}(u_n)<\infty$, we also have
$$\aligned&|\mathcal{E}^{a,\hat{b}}(u_n,(u_n)^+)-\mathcal{E}^{a,\hat{b}}(u,u^+)|
\rightarrow0.\endaligned$$
As a result we have (C.4) for  $u\in bF$.
So we have $\|u_{t_1}^+\|_2^2 \leq2\int_{t_1}^{t_2}(f_t,u_t^+)dt+\|u_{t_2}^+\|_2^2.$
$\hfill\Box$
\vskip.10in
To extend the class of solutions we are working with  to allow $f\in L^1(dt\times dm)$, we need the following proposition. It is a modified version of the above lemma.
\vskip.10in
\th{Lemma C.2} Let $u\in b\hat{F}$  and $f\in L^1(dt\times dm)$ satisfying the weak relation (2.9)  with test functions in $b\mathcal{C}_T$ and some function $\phi\geq0$, $\phi\in L^2\cap L^\infty$. Then $u^+$ satisfies the following relation with
$0\leq t_1<t_2\leq T$
$$\|u_{t_1}^+\|_2^2\leq 2\int_{t_1}^{t_2}(f_s,u_s^+)ds+\|u_{t_2}^+\|_2^2.$$

\proof First note that we can prove analogously to the above proof that for each $u\in b\mathcal{C}_T$ satisfying the weak relation (2.9) with data $(\phi,f)$ over the interval $[t_1,t_2]$, where $\varepsilon\leq t_1\leq t_2\leq T$ for $\varepsilon>0$, the following holds
$$\|u_{t_1}^+\|_2^2 \leq2\int_{t_1}^{t_2}(f_t,u_t^+)dt+\|u_{t_2}^+\|_2^2.$$
For $u\in\hat{F}$ we take approximating functions $u^n$ and $(\phi^n,f^n)$ as in the last part of the proof of Proposition 2.9 . And  we have that $u^n$ satisfies the weak relation (2.9) for the data $\phi^n,f^n$ with  test functions in $b\mathcal{C}_T$ over the interval $[\varepsilon,t_2]$ and $\frac{1}{n}\leq\varepsilon\leq t_2\leq T$. Note $\lim_{n\rightarrow\infty}\int_\varepsilon^T\|f_t^n-f_t\|_1dt=0.$
Then we have
$\|(u_{t_1}^n)^+\|_2^2 \leq2\int_{t_1}^{t_2}(f_t^n,(u_t^n)^+)dt+\|(u_{t_2}^n)^+\|_2^2,$
 where $\varepsilon\leq t_1\leq t_2\leq T$ for $\varepsilon>0$. The convergence of all terms, which do not depend on $f$, follows by the same arguments as in the above proof. Since $u$ is bounded, it is easy to see that $u^n$ is uniformly bounded.
Then we have $\lim_{n\rightarrow\infty}| \int_{t_1}^{t_2}(f_s^n,(u_s^n)^+)ds-\int_{t_1}^{t_2}(f_s,u_s^+)ds|
=0.$
Finally, we obtain
$\|u_{t_1}^+\|_2^2 \leq2\int_{t_1}^{t_2}(f_t,u_t^+)dt+\|u_{t_2}^+\|_2^2,$
where $\varepsilon\leq t_1\leq t_2\leq T$ for $\varepsilon>0$. Letting $\varepsilon\rightarrow0$, the assertion follows.$\hfill\Box$
\vskip.10in
The next Proposition is a modification of [2, Proposition 2.9]. It represents a version of the maximum principle.
\vskip.10in
\th{Proposition C.3} Let $u\in b\hat{F}$ and $f\in L^1(dt\times dm),f\geq0$, satisfying the weak relation (2.9) with test functions in $b\mathcal{C}_T$ and some function $\phi\geq0$, $\phi\in L^2\cap L^\infty$. Then $u\geq 0$ and it is represented by the following relation:
$$ u_t=P_{T-t}\phi+\int_t^TP_{s-t}f_sds.$$
Here we use $P_t$ is a $C_0$-semigroup on $L^1(\mathbb{R}^d;m)$ to make $P_{s-t}f_s$ meaningful.

\proof Let $(f^n)_{n\in N}$ be a sequence of bounded functions such that
$0\leq f^n\leq f^{n+1}\leq f, \lim_{n\rightarrow\infty}f^n=f.$
Since $f^n$ is bounded, we have $f^n\in L^1([0,T];L^2)$. Define
$ u_t^n:=P_{T-t}\phi+\int_t^TP_{s-t}f_s^nds.$
Then by Proposition 2.9, $u^n\in\hat{F}$ is a unique generalized solution for the data $(\phi,f^n)$. Clearly $0\leq u^n\leq u^{n+1}$ for $n\in N$. Define $y:=u^n-u$ and $\tilde{f}:=f^n-f$. Then $\tilde{f}\leq 0$ and $y$ satisfies the weak relation (2.9) for the data $(0,\tilde{f})$. Therefore by Lemma C.2, we have for $t_1\in [0,T]$
$\|y_{t_1}^+\|_2^2\leq 2\int_{t_1}^{T}(\tilde{f}_s,y_s^+)ds\leq 0.$
We conclude that $\|y_{t_1}^+\|_2^2=0$. Therefore, $u\geq u^n\geq0$ for $n\in N$. Set $v:=\lim_{n\rightarrow\infty}u^n$. By  (2.7) we have
$$\|u_t^n\|_2^2+2\int_t^T\mathcal{E}^{a,\hat{b}}(u_s^n)ds\leq 2\int_t^T(f_s^n,u_s^n)ds+\|\phi\|_2^2+2\alpha\int_t^T\|u_s^n\|_2^2ds,$$
which implies that
$$\|u_t^n\|_2^2+2\int_t^T\mathcal{E}^{a,\hat{b}}(u_s^n)ds\leq 2M\int_t^T\int|f_s^n|dmds+\|\phi\|_2^2+2\alpha\int_t^T\|u_s^n\|_2^2ds.$$
By Gronwall's Lemma, we have $\sup_n\sup_{t\in[0,T]}\|u_t^n\|_2^2\leq$const. And we obtain
 $\lim_{n\rightarrow\infty}\|u_t^n-v_t\|_2^2=0$ and
$\lim_{n\rightarrow\infty}|\int_t^T\int(f_s^nu_s^n-f_sv_s)dmds|=0.$
By [14, Lemma 2.12] we have
$$\int_t^T\mathcal{E}^{a,\hat{b}}_{c_2+1}(v_s)ds\leq \int_t^T\liminf_{n\rightarrow\infty}\mathcal{E}^{a,\hat{b}}_{c_2+1}(u_s^n)ds\leq\liminf_{n\rightarrow\infty}\int_t^T\mathcal{E}^{a,\hat{b}}_{c_2+1}(u_s^n)ds.$$
Finally, for $t\in[0,T]$ we get
$$\aligned\|v_t\|_2^2+2\int_t^T\mathcal{E}^{a,\hat{b}}(v_s)ds\leq&\lim_{n\rightarrow\infty}\|u_t^n\|_2^2
+\liminf_{n\rightarrow\infty}\int_t^T\mathcal{E}^{a,\hat{b}}(u_s^n)ds\\
\leq &\lim_{n\rightarrow\infty}(2\int_t^T(f_s^n,u_s^n)ds+
\|\phi\|_2^2)+\lim_{n\rightarrow\infty}2\alpha\int_t^T\|u_s^n\|_2^2ds\\
=&2\int_t^T(f_s,v_s)ds+\|\phi\|_2^2+\int_t^T\|v_s\|_2^2ds.\endaligned$$
Since the right hand  side of this inequality is finite and $t\mapsto v_t$ is $L^2$-continuous, it follows that $v\in\hat{F}$.

Now we show that $v$ satisfies the weak relation (2.9) for the data $(\phi,f)$. As $\varphi^n(t):=\|u_t^n-v_t\|_2$ is continuous and decreasing, we conclude by Dini's theorem
$\lim_{n\rightarrow\infty}\sup_{t\in[0,T]}\|u_t^n-v_t\|_2=0,$
and therefore
$\lim_{n\rightarrow\infty}\int_0^T\|u_t^n-v_t\|_2^2=0.$
Furthermore, there exists $K\in \mathbb{R}_+$ and a subsequence $(n_k)_{k\in N}$ such that
$|\int_0^T\mathcal{E}^{a,\hat{b}}_{c_2+1}(u_s^{n_k})ds|\leq K, \forall k\in N.$
in particular $\int_0^T\int|D_\sigma u_s^{n_k}|^2dmds\leq \frac{ K }{c_1}, \forall k\in N.$
We obtain  $\lim_{k\rightarrow\infty}\int_0^T\mathcal{E}^{a,\hat{b}}(u_s^{n_k},\varphi_s)ds=\int_0^T\mathcal{E}^{a,\hat{b}}(v_s,
\varphi_s)ds,$
and
$\lim_{k\rightarrow\infty}\int_0^T\int\langle b\sigma,D_\sigma u_s^{n_k}\rangle\varphi_sdmds=\int_0^T\int\langle b\sigma,D_\sigma v_s\rangle\varphi_sdmds,$
which implies (2.9) for $v$ associated to $(\phi,f)$. Clearly $u-v$ satisfies (2.9) with data $(0,0)$ for $\varphi\in b\mathcal{C}_T$. By Proposition 2.9 we have $u-v=0$. Since
$$ v_t=P_{T-t}\phi+\int_t^TP_{s-t}f_sds,$$
the assertion follows.$\hfill\Box$
\vskip.10in
\th{Corollary C.4} Let $u\in b\hat{F}$  and $f\in L^1(dt\times dm)$ satisfy the weak relation (2.9)  with test functions in $b\mathcal{C}_T$ and some function  $\phi\in L^2\cap L^\infty$. Let $g\in L^1(dt\times dm)$ be a bounded function such that $f\leq g$. Then $u$  is represented by the following relation:
$$ u_t=P_{T-t}\phi+\int_t^TP_{s-t}f_sds.$$
\proof Define $f^n:=(f\vee(-n))\wedge g, n\in \mathbb{N}$. Then $(f^n)_{n\in \mathbb{N}}$ is a sequence of bounded functions such that $f^n\downarrow f$ and $f^n\leq g$ then by the same arguments as in Proposition C.3, the assertion follows.$\hfill\Box$
\vskip.10in

The following proposition is a modification of [2, Proposition 2.10] . It is essential for the analytic treatment of the  non-linear equation (1.1) which is done in the next Section.

\vskip.10in
\th{Proposition C.5} Let $u=(u^1,...,u^l)$ be a vector valued function where each component is a generalized solution of the linear equation (2.5) associated to certain data $f^i,\phi^i$, which are bounded and satisfy the conditions in Proposition 2.6 (ii)  for $i=1,...,l$. Denote by $\phi,f$ the vectors $\phi=(\phi^1,...,\phi^l),f=(f^1,...,f^l)$ and by $D_\sigma u$ the matrix whose rows consist of the row vectors $D_\sigma u^i$. Then the following relations hold $m$-almost everywhere
$$
 |u_t|^2+2\int_t^TP_{s-t}(|D_\sigma u_s|^2+\frac{1}{2}c|u_s|^2)ds=P_{T-t}|\phi|^2+2\int_t^TP_{s-t}\langle u_s,f_s\rangle ds.
\eqno(C.5)$$
$$
 |u_t|\leq P_{T-t}|\phi|+\int_t^TP_{s-t}\langle \hat{u}_s,f_s\rangle ds.
\eqno(C.6)$$
Here we write $\hat{x}=x/|x|$, for $x\in \mathbb{R}^l$, $x\neq0$ and $\hat{x}=0$, if $x=0$.

\proof By Proposition 2.6 (ii) we have $u\in b\mathcal{C}_T$.

First we assume $l=1$.
If we can check that $u^2$ satisfies (2.9) with data $(2uf-2|D_\sigma u|^2-cu^2,\phi^2)$ for $\varphi\in b\mathcal{C}_T$,  then (C.5) will follow by Corollary C.4.
We have the following relations:
$$\int_0^T(u_t^2,\partial_t\varphi_t)dt=2\int_0^T(u_t,\partial_t(u_t\varphi_t))dt+(u_0^2,\varphi_0)-(u_T^2,\varphi_T),$$
$$\mathcal{E}^{a,\hat{b}}(u_t^2,\varphi_t)=2\mathcal{E}^{a,\hat{b}}(u_t,u_t\varphi_t)-(2|D_\sigma u_t|^2+cu^2_t,\varphi_t),$$
and
$$\int\langle b\sigma, D_\sigma (u_t^2)\rangle\varphi_tdm=2\int\langle b\sigma, D_\sigma u_t\rangle u_t\varphi_tdm.$$
For the second relation, we use (2.2).
Since $u$ is a generalized solution of (2.5), we have
$$\aligned&\int_0^T(u_t,\partial_t(u_t\varphi_t))dt-(u_T,u_T\varphi_T)+(u_0,u_0\varphi_0)-\int_0^T(f_t,u_t\varphi_t)dt
\\=&-\int_0^T\mathcal{E}^{a,\hat{b}}(u_t,u_t\varphi_t)dt-\int_0^T\int\langle b\sigma, D_\sigma u_t\rangle u_t\varphi_tdmdt.\endaligned$$
By the above relations, we have
$$\aligned&\int_0^T(u_t^2,\partial_t\varphi_t)dt+(u_0^2,\varphi_0)-(u_T^2,\varphi_T)+\int_0^T(\mathcal{E}^{a,\hat{b}}(u_t^2,\varphi_t)+\int\langle b\sigma, D_\sigma (u_t^2)\rangle\varphi_tdm)dt\\
=&2\int_0^T(f_tu_t,\varphi_t)dt-\int_0^T(2|D_\sigma u_t|^2+c|u_t|^2,\varphi_t)dt.\endaligned\eqno(C.7)$$
Hence, by Corollary C.4 (C.5) holds in the case $l=1$. To deduce this relation in the case $l>1$, it suffices to add the relations corresponding to the components $|u_t^i|^2,i=1,...,l$.

For (C.6), define for $\varepsilon>0$,
$h_\varepsilon(t):=\sqrt{t+\varepsilon}-\sqrt{\varepsilon}$ for $t\geq 0$.
Then by integration by parts, we have
$$\aligned\mathcal{E}^{a,\hat{b}}(h_\varepsilon(|u|^2),\varphi)=&
\mathcal{E}^{a,\hat{b}}(|u|^2,h'_\varepsilon(|u|^2)\varphi)-(h''_\varepsilon(|u|^2)|D_\sigma(|u|^2)|^2,\varphi)
\\&+(c(h_\varepsilon(|u|^2)-|u|^2h'_\varepsilon(|u|^2)),\varphi),\endaligned$$
an
$$\aligned\int_0^T(h_\varepsilon(|u_t|^2),\partial_t\varphi_t)dt=&\int_0^T(|u_t|^2,\partial_t(\varphi_t h_\varepsilon'(|u_t|^2)))dt-(|u_T|^2,\varphi_Th_\varepsilon'(|u_T|^2))\\&+(|u_0|^2,\varphi_0h_\varepsilon'(|u_0|^2))+(h_\varepsilon(|u_T|^2),\varphi_T)-(h_\varepsilon(|u|_0^2),\varphi_0).\endaligned$$
 the above relations we obtain
$$\aligned&\int_0^T(h_\varepsilon(|u_t|^2),\partial_t\varphi_t)dt-(h_\varepsilon(|u_T|^2),\varphi_T)+(h_\varepsilon(|u_0|^2),\varphi_0)\\&
+\int_0^T(\mathcal{E}^{a,\hat{b}}(h_\varepsilon(|u_t|^2),\varphi_t)+\int\langle b\sigma, D_\sigma (h_\varepsilon(|u_t|^2))\rangle\varphi_tdm)dt\\=&\int_0^T-(h''_\varepsilon(|u_t|^2)|D_\sigma(|u_t|^2)|^2,\varphi_t)+(c(h_\varepsilon(|u|^2)-|u|^2h'_\varepsilon(|u|^2)),\varphi)dt\\&+2\int_0^T(\langle f_t,u_t\rangle h'_\varepsilon(|u_t|^2),\varphi_t)dt-\int_0^T(h'_\varepsilon(|u_t|^2)(2|D_\sigma u_t|^2+c|u_t|^2,\varphi_t)dt.\endaligned$$

As $ |D_\sigma(|u|^2)|^2=4\langle u, D_\sigma u (D_\sigma u)^*u\rangle,$
we deduce
$$\aligned&2\langle f,u\rangle h'_\varepsilon(|u|^2)-2h'_\varepsilon(|u|^2)|D_\sigma u|^2-h''_\varepsilon(|u|^2)|D_\sigma(|u|^2)|^2\\
=&\frac{\langle f,u\rangle-|D_\sigma u
|^2}{(|u|^2+\varepsilon)^{\frac{1}{2}}}+\frac{|u|^2\langle \hat{u}, D_\sigma u (D_\sigma u)^*\hat{u}\rangle}{(|u|^2+\varepsilon)^{\frac{3}{2}}
}\\ \leq &\frac{\langle f,u\rangle}{(|u|^2+\varepsilon)^{\frac{1}{2}}}.\endaligned$$
By Proposition C.5 and since $c(h_\varepsilon(|u_s|^2)-2|u_s|^2h'_\varepsilon(|u_s|^2))\leq 0$, we deduce
$$h_\varepsilon(|u_t|^2)\leq P_{T-t}h_\varepsilon(|\phi|^2)+\int_t^TP_{s-t}\frac{\langle f_s,u_s\rangle}{(|u_s|^2+\varepsilon)^{\frac{1}{2}}}ds.$$
Letting $\varepsilon\rightarrow0$ obtain the results.$\hfill\Box$
\vskip.10in
Proposition 2.10 is a version of the above proposition for general data. Here we use $P_t$ is a $C_0$-semigroup on $L^1$.
\vskip.10in

\vspace{1mm}\noindent{\it Proof of Proposition 2.10}\quad Analogously to the proof of Proposition C.5 it is enough to verify (2.12) for $l=1$. For $\phi\in L^2, f\in L^1([0,T],L^2)$, take $\phi_n,f_n$ as in  Proposition 2.9, then we have

(a). $ u_{n,t}:=P_{T-t}\phi_n+\int_t^TP_{s-t}f_{n,s}ds \textrm{ is a generalized solution },$

(b). $ \lim_{n\rightarrow\infty}\int_t^T\|f_{n,s}-f_s\|_2ds=0,$

(c). $ \lim_{n\rightarrow\infty}\|\phi_n-\phi\|_2=0,$

(d). $ \lim_{n\rightarrow\infty}\|u_n-u\|_T=0.$

\no By Proposition C.5 we have
$$
 |u_{n,t}|^2+2\int_t^TP_{s-t}(|D_\sigma u_{n,s}|^2+\frac{1}{2}c|u_{n,s}|^2)ds=P_{T-t}|\phi_n|^2+2\int_t^TP_{s-t}\langle u_{n,s},f_{n,s}\rangle ds.
\eqno(C.8)$$
By (b) and (d) we obtain
$$\aligned&\|\int_t^TP_{s-t}((u_{n,s},f_{n,s})-(u_s,f_s))ds\|_1\\
\leq &C(\sup_{s\in[0,T]}\|u_{n,s}\|_2\int_t^T\|f_{n,s}-f_s\|_2ds+\sup_{s\in[0,T]}\|u_{n,s}-u_s\|\int_t^T\|f_s\|_2ds)
\rightarrow0.\endaligned$$
Here we used that $P_t$ is a $C_0$-semigroup on $L^1(\mathbb{R}^d;m)$.
By (d) we conclude that
$$\aligned&\int_t^T\||D_\sigma u_{n,s}|^2-|D_\sigma u_{s}|^2\|_1ds\\ \leq &
((\int_t^T\|D_\sigma u_{n,s}\|_2^2ds)^\frac{1}{2}+(\int_t^T\|D_\sigma u_{s}\|_2^2ds)^\frac{1}{2})(\int_t^T\|D_\sigma u_{n,s}-D_\sigma u_{s}\|_2^2ds)^\frac{1}{2}
\rightarrow0,\endaligned$$
and that
$$\aligned&\int_t^T\||c u_{n,s}|^2-|c u_{s}|^2\|_1ds\\\leq&M((\int_t^T\mathcal{E}_{c_2+1}^{a,\hat{b}}(u_{n,s})ds)^\frac{1}{2}+(\int_t^T\mathcal{E}_{c_2+1}^{a,\hat{b}}(u_{s})ds)^\frac{1}{2})(\int_t^T\mathcal{E}^{a,\hat{b}}_{c_2+1}(u_{n,s}-u_s)ds)^\frac{1}{2}\rightarrow0.\endaligned$$
Thus, we obtain
$\lim_{n\rightarrow\infty}\int_t^TP_{s-t}(|D_\sigma u_{n,s}|^2)ds=\int_t^TP_{s-t}|D_\sigma u_{s}|^2ds,$
and
$\lim_{n\rightarrow\infty}\int_t^TP_{s-t}(|c u_{n,s}|^2)ds=\int_t^TP_{s-t}|c u_{s}|^2ds.$
Passing to the limit in equation (C.8) we get (2.12).

For proving (2.13), we use the same method.$\hfill\Box$
\vskip.10in
By the above results, we obtain the following lemma by the same arguments as the proof in [2, Lemma 2.12].
\vskip.10in
\th{Lemma C.6} If $f,g\in L^1([0,T];L^2)$ and $\phi\in L^2$, then:
$$\int_t^TP_{s-t}(f_sP_{T-s}\phi)ds\leq \frac{1}{2}P_{T-t}\phi^2+\int_t^T\int_s^TP_{s-t}(f_sP_{r-s}f_r)drds,\textrm{  } m-a.e.\eqno(C.9)$$

\end{document}